\documentclass[12pt]{elsarticle}
\usepackage[left=1in,right=1in,top=1in,bottom=1in]{geometry}%
\usepackage[utf8]{inputenc}
\usepackage[usenames,dvipsnames]{xcolor}

\usepackage{amsmath}
\usepackage{amssymb}
\usepackage{amsthm}
\usepackage{mathtools}

\usepackage{graphicx}
\usepackage{enumerate}
\usepackage{bbm}
\usepackage{mdwlist}

\usepackage{macros}
\usepackage{packages}
\usepackage{style}

\RequirePackage{xcolor}[2.11]
\colorlet{siaminlinkcolor}{green!50!black}
\colorlet{siamexlinkcolor}{red!50!black}
\colorlet{siamreviewcolor}{black!50}
\hypersetup{
  colorlinks = false,
  allcolors = siaminlinkcolor,
  urlcolor = siamexlinkcolor,
}
\hypersetup{hidelinks}

\theoremstyle{plain}
\newtheorem{theorem}{Theorem}
\newtheorem{corollary}[theorem]{Corollary}
\newtheorem{lemma}[theorem]{Lemma}

\theoremstyle{definition}
\newtheorem{definition}[theorem]{Definition}

\theoremstyle{remark}
\newtheorem{remark}[theorem]{Remark}
\newtheorem{hypothesis}[theorem]{Hypothesis} 

\numberwithin{theorem}{section}
\numberwithin{equation}{section}

\begin{document}

\begin{frontmatter}

\title{Periodic multi-pulses and spectral stability in Hamiltonian PDEs with symmetry}

\author[1]{Ross Parker}
    \ead{rhparker@smu.edu}
\author[2]{Bj\"{o}rn Sandstede}
    \ead{bjorn\_sandstede@brown.edu}

\address[1]{Department of Mathematics, Southern Methodist University, 
Dallas, TX 75275, USA}
\address[2]{Division of Applied Mathematics, Brown University, Providence, RI 02912, USA}

\begin{abstract}
We consider the existence and spectral stability of periodic multi-pulse solutions in Hamiltonian systems which are translation invariant and reversible, for which the fifth-order Korteweg-de Vries equation is a prototypical example. We use Lin's method to construct multi-pulses on a periodic domain, and in particular demonstrate a pitchfork bifurcation structure for periodic double pulses. We also use Lin's method to reduce the spectral problem for periodic multi-pulses to computing the determinant of a block matrix, which encodes both eigenvalues resulting from interactions between neighboring pulses and eigenvalues associated with the essential spectrum. We then use this matrix to compute the spectrum associated with periodic single and double pulses. Most notably, we prove that brief instability bubbles form when eigenvalues collide on the imaginary axis as the periodic domain size is altered. These analytical results are all in good agreement with numerical computations, and numerical timestepping experiments demonstrate that these instability bubbles correspond to oscillatory instabilities.
\end{abstract}


\end{frontmatter}

\section{Introduction}\label{sec:intro}

Solitary waves, localized disturbances that maintain their shape as they propagate at a constant velocity, have been an object of mathematical and experimental interest since the nineteenth century \cite{KdVoriginal} and have applications not only in fluid mechanics but also nonlinear optics \cite{Taylor1992}, molecular systems \cite{Davydov1985}, Bose-Einstein condensates \cite{Panos2008BEC}, and ferromagnetics \cite{Kosevich1998}. Of more recent interest are multi-pulses, which are multi-modal solitary waves resembling multiple, well-separated copies of a single solitary wave. 
The entire multi-pulse travels as a unit, and it maintains its shape unless perturbed.
The study of multi-pulses goes back to at least the early 1980s, where Evans, Fenichel, and Faroe proved the existence of a double pulse traveling wave in nerve axon equations \cite{Evans1982}. The stability of these double pulses was shown in \cite{Yanagida1989}, and the existence result was extended to arbitrary multi-pulses in \cite{Feroe1986}. The existence of multi-pulse traveling wave solutions to semilinear parabolic equations, which includes reaction-diffusion systems, was established in \cite{Alexander1994}, and the stability of these solutions was determined using the Evans function, an analytic function whose zeros coincide with the point spectrum of a linear operator \cite{KapitulaSandstede2004}.
Existence of multi-pulse solutions to a family of Hamiltonian equations was shown in \cite{Buffoni1996} using the dynamics on the Smale horseshoe set, and a spatial dynamics approach to the same problem is found in \cite{SandstedeStrut}. Multi-pulses have since been studied in diverse systems, including a pair of nonlinearly coupled Schr\"{o}dinger equations \cite{Yew2001,Yew2000}, coupled nonlinear Schr\"{o}dinger equations \cite{Pelinovsky2001,Pelinovsky2005}, the vector nonlinear Schr\"{o}dinger equation \cite{Kapitula2007}, and lattice systems such as the discrete nonlinear Schr\"{o}dinger equation \cite{parkerlattice} and the discrete sine-Gordon equation \cite{ParkerSineGordon}. In general, the spectrum of the linearization of the underlying PDE about a multi-pulse contains a finite set of eigenvalues close to 0 \cite{Alexander1990,Sandstede1998}. Since these result from nonlinear interactions between the tails of neighboring pulses, we call them interaction eigenvalues. Under the assumption that the essential spectrum lies in the left half plane, spectral stability of multi-pulses depends on these interaction eigenvalues. For semilinear parabolic equations, these eigenvalues are computed in \cite{Sandstede1998} by using Lin's method, an implementation of the Lyapunov-Schmidt technique, to reduce the eigenvalue problem to a matrix equation. An extension of this technique was used to study the existence and spectral stability of multi-pulses in systems with both reflection and phase symmetries, such as the complex cubic-quintic Ginzburg-Landau equation \cite{Manukian2009}. This was further adapted to multi-pulses in certain Hamiltonian systems with two continuous symmetries, such as a fourth order nonlinear Schr\"{o}dinger equation \cite{ParkerAceves2021}.

\begin{figure}
\begin{center}
\begin{tabular}{cc}
\includegraphics[width=7cm]{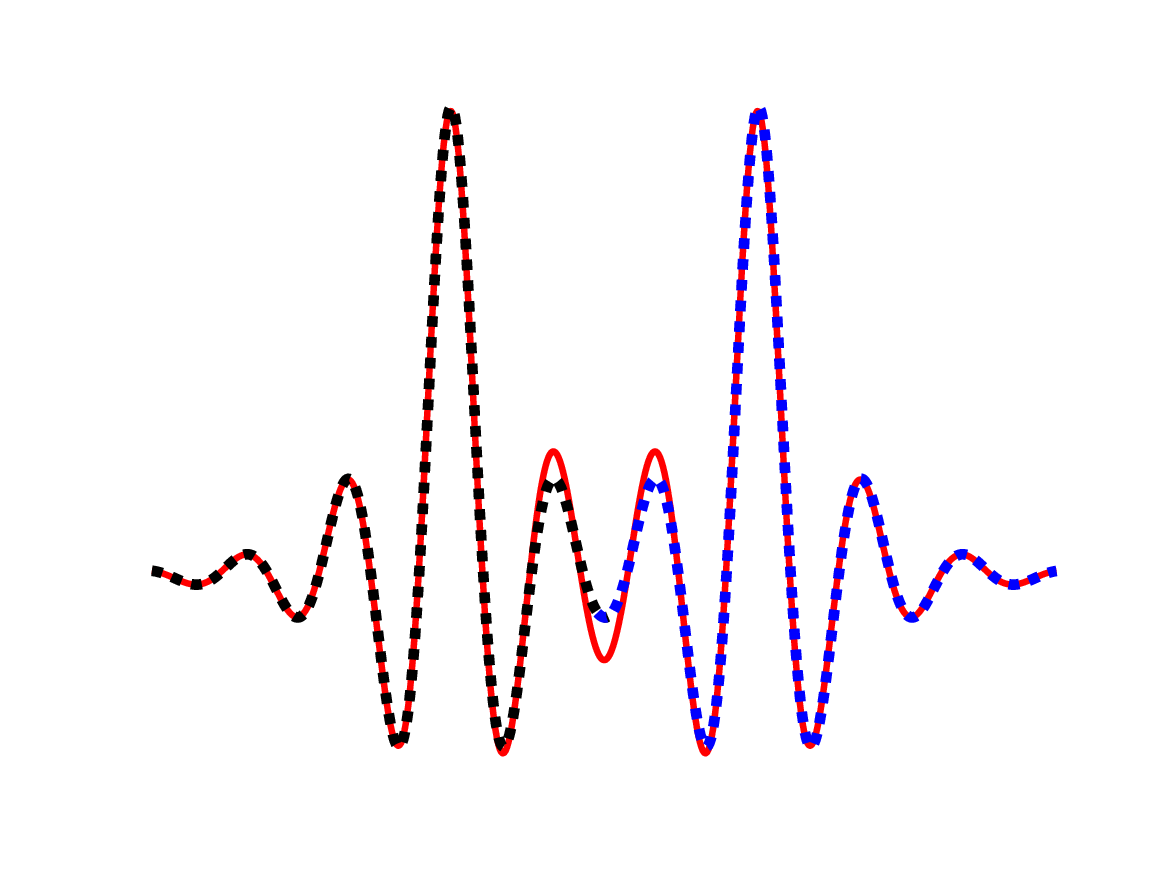} &
\includegraphics[width=7cm]{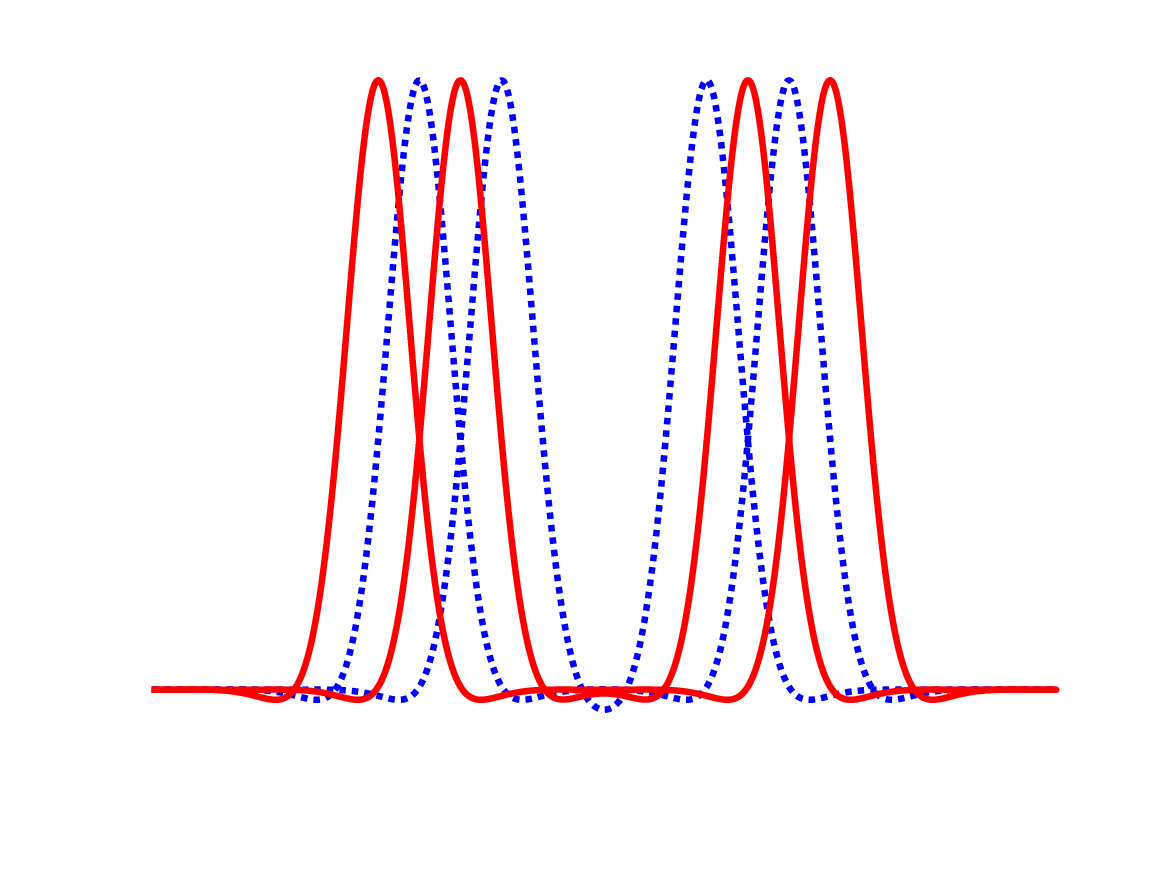}
\end{tabular}
\end{center}
\caption{Left panel is schematic of construction of double pulse solution (solid line) from two single pulses (dotted lines). Oscillatory tails are exaggerated. Right panel shows the first four double pulse solutions to \cref{KdV5example} for $c = 10$.}
\label{fig:KdV5double}
\end{figure}

\begin{figure}
\begin{center}
\begin{tabular}{cc}
\includegraphics[width=8cm]{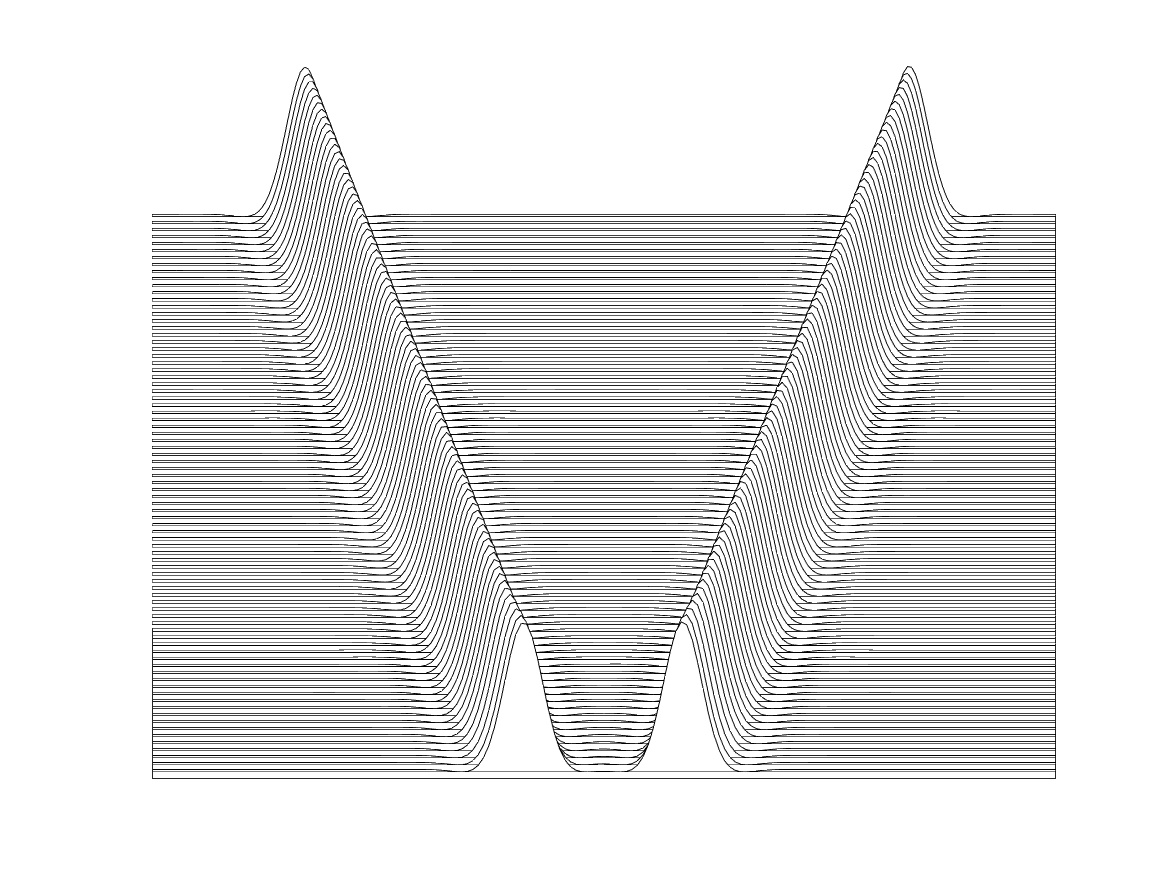} &
\includegraphics[width=8cm]{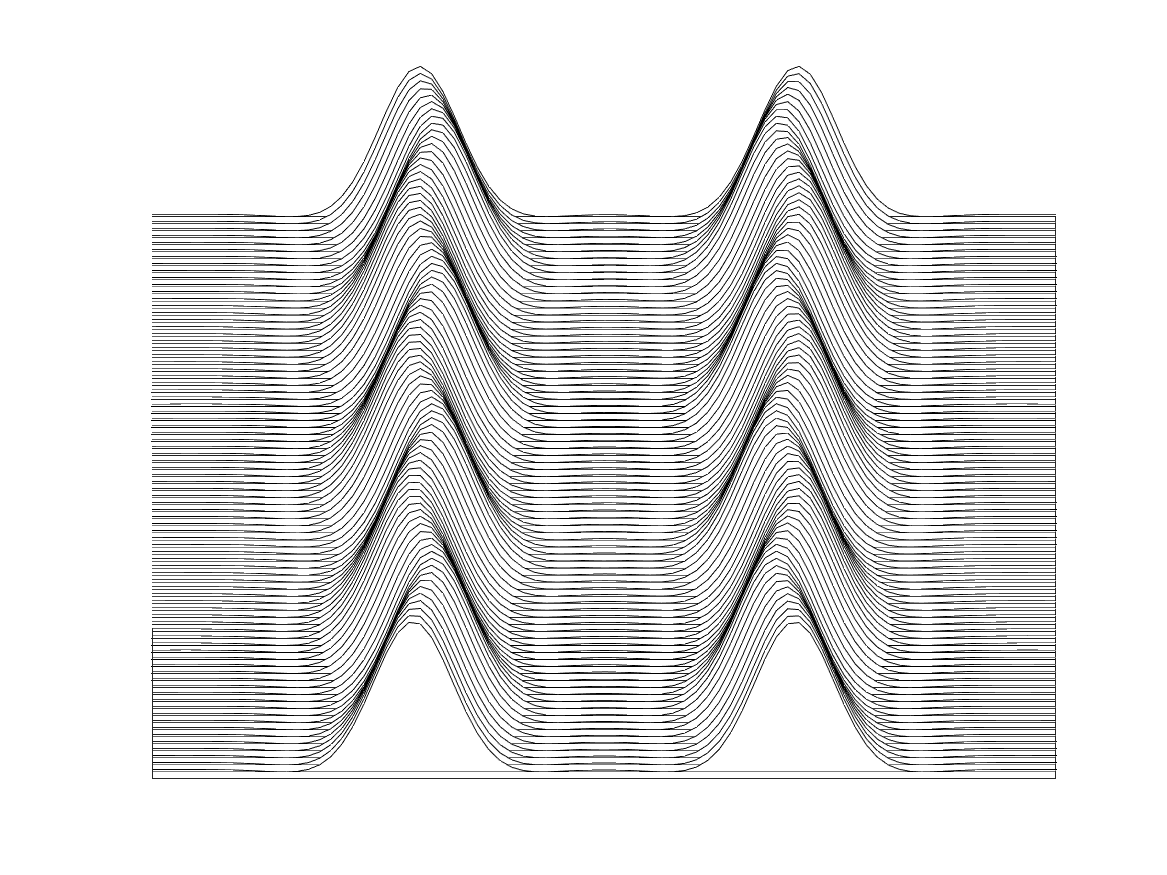}
\vspace{-0.25cm}
\end{tabular}
\end{center}
\caption{Results of numerical timestepping simulations for perturbations of double pulse solutions to \cref{KdV5example} with real eigenvalues (left panel) and imaginary eigenvalues (right panel). Crank-Nicolson/Adams-Bashforth 2 IMEX scheme in time with Chebyshev spectral discretization, Dirichlet and Neumann boundary conditions.}
\label{fig:KdV5waterfall}
\end{figure}

A much more difficult problem concerns the spectral stability of multi-pulses in Hamiltonian PDEs in the case where the essential spectrum consists of the entire imaginary axis. (If the essential spectrum is imaginary but bounded away from the origin, the spectral problem is considerably easier; see, for example, applications to a fourth-order beam equation \cite[Section 6]{Kapitula2020} and a fourth-order nonlinear Schr\"{o}dinger equation \cite{ParkerAceves2021}.) As a concrete example, we consider solitary wave solutions to the fifth-order Korteweg-de Vries equation (KdV5)
\begin{align}\label{KdV5example}
u_t &= \partial_x\left( u_{xxxx} - u_{xx} + c u - u^2\right) && c > 1/4,
\end{align}
which is the equation studied in \cite{Pelinovsky2007} written in a moving frame with speed $c$. 
For the remainder of this paper, we will consider the wavespeed $c$ to be a fixed parameter. For $c>1/4$, the solitary wave solutions will have oscillatory, exponentially decaying tails (see \cref{hyp:hypeq} and \cref{remark:hypsatisfied} below).
A double pulse solution resembles two well-separated copies of the primary solitary wave which are joined together in such a way that the tail oscillations ``match up'' (\cref{fig:KdV5double}, left panel). The distance between the two peaks takes values in a discrete set (\cref{fig:KdV5double}, right panel) \cite{Pelinovsky2007,SandstedeStrut}. This constraint is a consequence of a specific alignment of the stable and unstable manifolds which is a necessary condition for multi-pulses to occur, and these discrete values represent the number of twists made by the manifolds near the equilibrium at the origin \cite{Sandstede1998}. 
\noi For the spectral problem, the essential spectrum is the entire imaginary axis, and depends only on the background state. In addition, there is a pair of interaction eigenvalues which is symmetric about the origin and alternates between real (corresponding to double pulses with dashed lines in the right panel of \cref{fig:KdV5double}) and imaginary with negative Krein signature (corresponding to double pulses with solid lines in the right panel of \cref{fig:KdV5double}) \cite{Pelinovsky2007}. Numerical timestepping verifies that double pulses with real eigenvalues are unstable; when perturbed, the two peaks move away from each other with equal and opposite velocities (\cref{fig:KdV5waterfall}, left panel). For the remaining double pulses, numerical timestepping suggests that the two peaks exhibit oscillatory behavior when perturbed (\cref{fig:KdV5waterfall}, right panel). Similar timestepping results can be seen in \cite[Figure 9]{Pelinovsky2007}, as well as a reduction of the system to a two-dimensional phase plane \cite[Figure 10]{Pelinovsky2007}. 
A different argument that half of the double pulses are stable can be found in \cite{Buryak1997}, which uses the asymptotic method of \cite{Gorshkov1981}. Stability of double pulses in KdV5 is also discussed in terms of the Maslov index, an integer-valued topological invariant associated with homoclinic orbits in a finite-dimensional Hamiltonian system, in \cite[Section 15.1]{Chardard1}. The methods in \cite{Chardard1} are extended to Hamiltonian systems with phase space of dimension greater than four in \cite{Chardard2}; a specific example is the 7th order KdV model considered in \cite[Section 8]{Chardard2}.

Although these numerical results suggest that every other double pulse is neutrally stable, this remains an open question, since the imaginary eigenvalues are embedded in the essential spectrum. Furthermore, the timestepping simulation in \cref{fig:KdV5waterfall} was performed with separated boundary conditions, which shifts the essential spectrum into the left half plane, and thus could fundamentally alter the behavior of the system. As an alternative, we will look at multi-pulse solutions on a periodic domain subject to co-periodic perturbations. The advantage is that the essential spectrum becomes a discrete set of points on the imaginary axis; by analogy to the problem on the real line, we will refer to this set as essential spectrum eigenvalues, even though they are elements of the point spectrum. Purely imaginary interaction eigenvalues can then lie between essential spectrum eigenvalues. Periodic traveling waves were described by Korteweg and de Vries in their 1895 paper \cite{KdVoriginal}, and the stability of these cnoidal waves is shown in \cite{Pava2006,Bottman2009}. Since then, stability of periodic solutions has been investigated for many other systems, including the generalized KdV equation \cite{Johnson2009}, the generalized Kuramoto-Sivashinsky equation \cite{Barker2013}, the Boussinesq equation \cite{Hakkaev2014}, the Klein-Gordon equation \cite{Demirkaya2015}, a generalized class of nonlinear dispersive equations \cite{Hur2015}, the regularized short pulse and Ostrovsky equations \cite{Hakkaev2017}, and the Lugiato-Lefever model of optical fibers \cite{Delcey2018,Milena2018,Hakkaev2019}.

\begin{figure}
\begin{center}
\includegraphics[width=8cm]{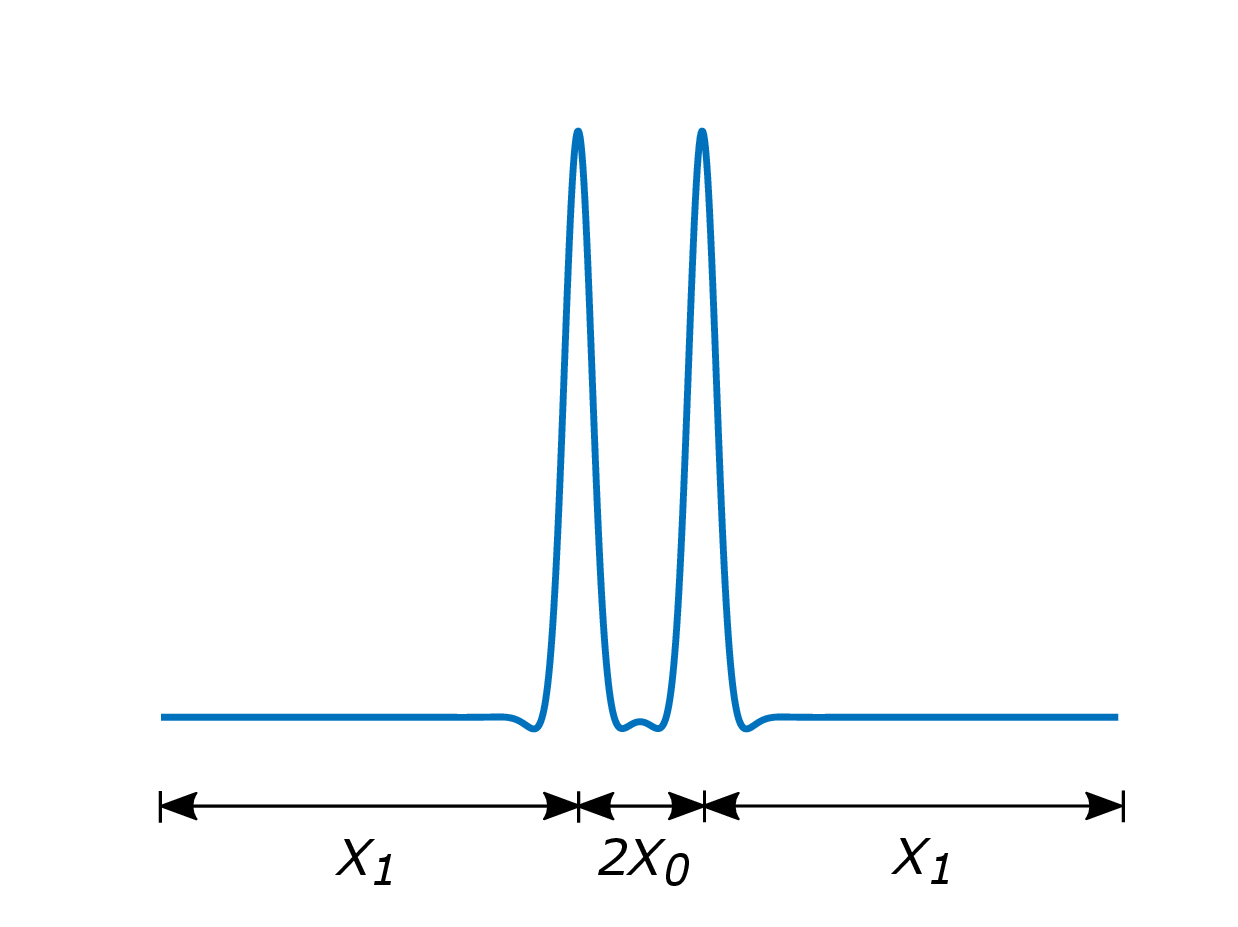} \hspace{-1cm}
\includegraphics[width=9cm]{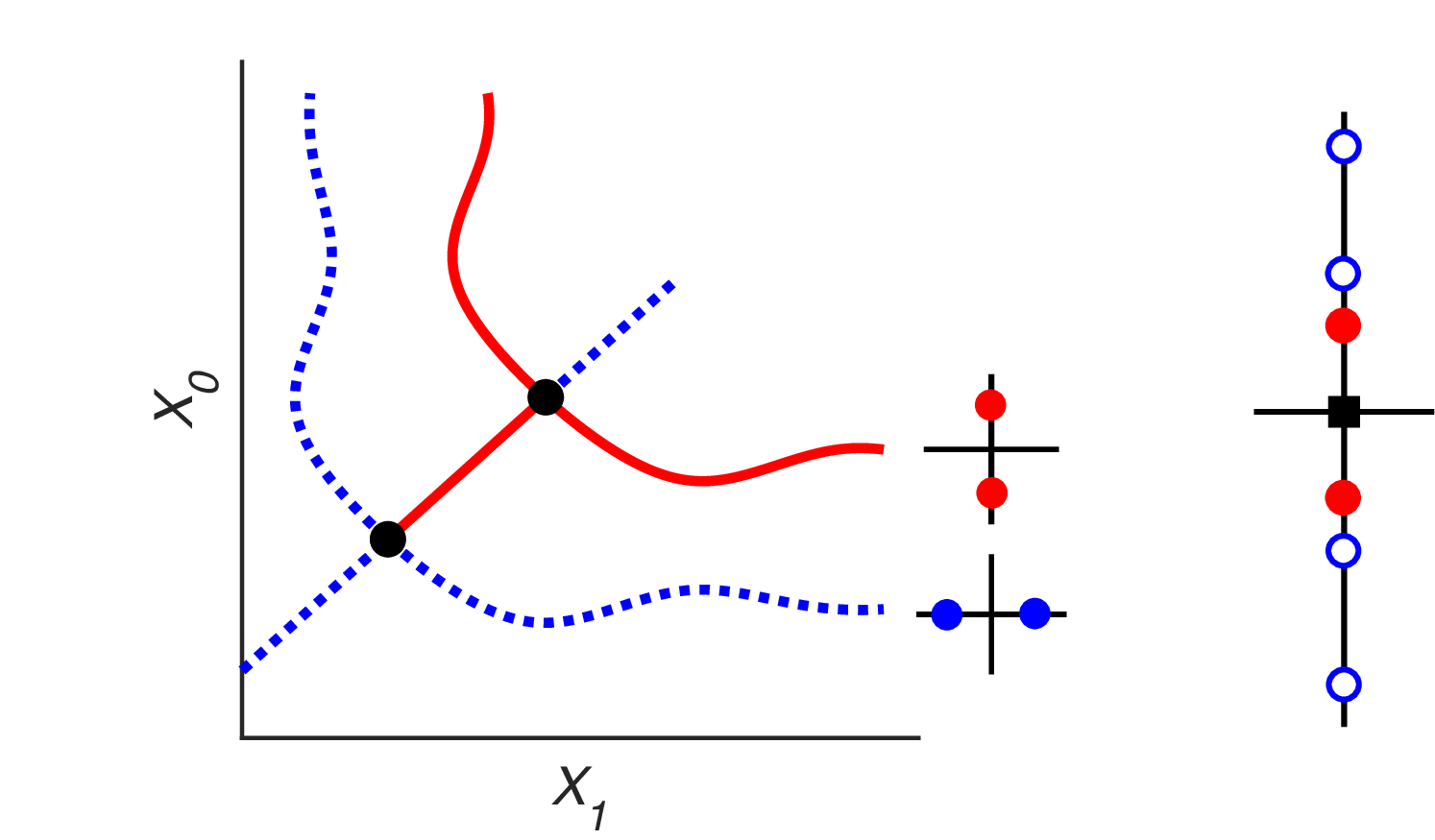}
\end{center}
\caption{Left panel is an illustration of a periodic 2-pulse showing the two length parameters $X_0$ and $X_1$. Center panel shows the pitchfork bifurcation structure for periodic 2-pulses as these length parameters are varied. Pitchfork bifurcation points are indicated with black dots. Blue dotted lines correspond to solutions with real interaction eigenvalues, and red solid lines correspond to solutions with imaginary interaction eigenvalues. Right panel is a schematic of the eigenvalue pattern for a periodic 2-pulse with a pair of imaginary interaction eigenvalues, corresponding to the solid red line in the center panel. The spectrum comprises an imaginary pair of interaction eigenvalues (red dots), a double eigenvalue at the origin from translation invariance (black square), and essential spectrum eigenvalues (blue open circles). The wavespeed $c$ is fixed.}
\label{fig:KdV5periodicMP}
\end{figure}

As in \cite{SandstedeStrut,Sandstede1998}, we will use a spatial dynamics approach. 
We note that since the wavespeed $c$ is a fixed parameter, all solutions obtained this way will be traveling waves with speed $c$. From this perspective, the primary solitary wave is a homoclinic orbit connecting the unstable and stable manifolds of a saddle equilibrium point at the origin. A multi-pulse is a multi-loop homoclinic orbit which remains close to the primary homoclinic orbit, and a periodic multi-pulse is a multi-loop periodic orbit. Unlike multi-pulses on the real line, which exist in discrete families (see \cref{fig:KdV5double}, right panel), periodic multi-pulses exist in continuous families, since there is an additional degree of freedom in their construction. Consider, for example, a 2-pulse. Whereas a 2-pulse on the real line can be described by a single length parameter representing the distance between the two peaks, the characterization of a periodic 2-pulse requires two length parameters $X_0$ and $X_1$ (\cref{fig:KdV5periodicMP}, left panel). The length of the periodic domain is $2X = 2X_0 + 2X_1$. 
Double pulse solutions on the real line correspond to the formal limit $X_1 \rightarrow \infty$. The length parameter $X_0$ (represented by the red solid and blue dotted horizontal lines in \cref{fig:KdV5periodicMP}) formally converges to the distance between the two peaks in the double pulse solution on the real line (these solutions are shown in the right panel of \cref{fig:KdV5double}). As a consequence of this additional degree of freedom and the reversibility of the system, symmetric periodic 2-pulses $(X_0 = X_1)$ exist for all sufficiently large $X_0$, and asymmetric periodic 2-pulses $(X_0 \neq X_1)$ bifurcate from these symmetric periodic 2-pulses in a series of pitchfork bifurcations (\cref{fig:KdV5periodicMP}, center panel).
The symmetric periodic 2-pulse solutions ($X_0 = X_1$) correspond to periodic single-pulse solutions with the period $2X_0$ repeated twice on the period $2X$. 

To compute the spectrum of the linearization of the underlying PDE about a periodic multi-pulse, we use Lin's method to reduce the eigenvalue problem to a block matrix equation; the block matrix encodes the interaction eigenvalues, the essential spectrum eigenvalues, and the translational eigenvalues (which are in the kernel of the linearization). For a periodic 2-pulse, the block matrix is $4\times4$, and the resulting equation can be solved. This yields the interaction eigenvalue pattern in the center panel of \cref{fig:KdV5periodicMP}, which corresponds exactly to the pitchfork bifurcation structure. The arms of the pitchforks alternate between solutions with a pair of real interaction eigenvalues and solutions with a pair of imaginary interaction eigenvalues; stability changes at the pitchfork bifurcation points, when the interaction eigenvalues collide at the origin. 

\begin{figure}
\begin{center}
\includegraphics[width=15cm]{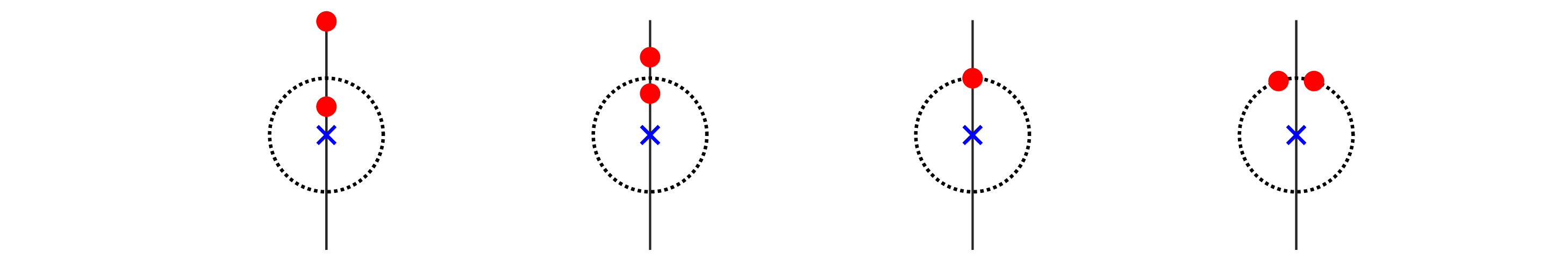}
\includegraphics[width=15cm]{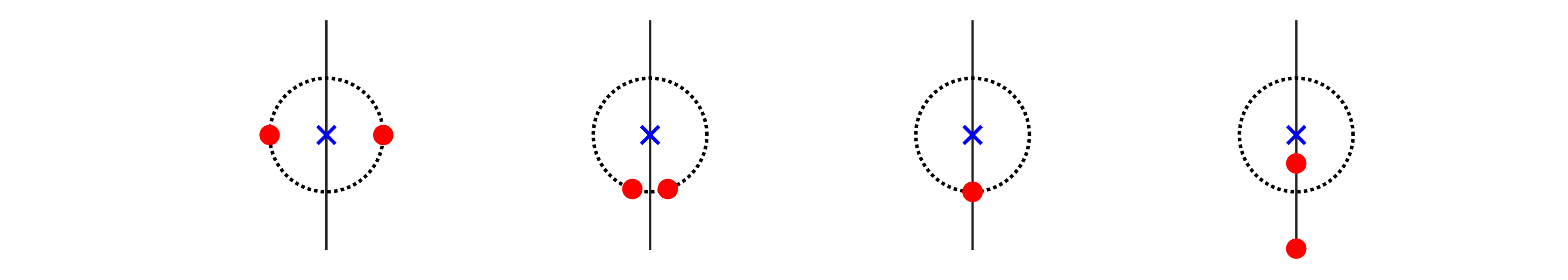}
\end{center}
\caption{Schematic showing brief instability bubble which forms as the imaginary interaction eigenvalue and the first essential spectrum eigenvalue collide on the imaginary axis. Domain size $X$ increases from left to right, and top to bottom.
Eigenvalues are shown in red dots, Krein bubble is dotted circle, and center of instability bubble is marked with blue X.}
\label{fig:KreinBubbleCartoon}
\end{figure}

The right panel of \cref{fig:KdV5periodicMP} is a schematic of the eigenvalue pattern of a periodic double pulse with imaginary interaction eigenvalues. The schematic also shows the first two essential spectrum eigenvalues. The essential spectrum eigenvalues are approximately equally spaced on the imaginary axis, and, to leading order, their location depends only on the domain length parameter $X$. As long as the interaction eigenvalues and the essential spectrum eigenvalues do not get too close, which we can guarantee by choosing the length parameters $X_0$ and $X_1$ sufficiently large, the interaction eigenvalue pattern is as shown in \cref{fig:KdV5periodicMP}. As the domain size $X$ is increased (moving to the right along the arm of the pitchfork corresponding to the red solid line in \cref{fig:KdV5periodicMP}), however, the essential spectrum eigenvalues move along the imaginary axis towards the origin. At a critical value of $X$, there is a collision between one of the essential spectrum eigenvalues and a purely imaginary interaction eigenvalue. Since the two eigenvalues have opposite Krein signatures, we expect them to leave the imaginary axis upon collision. In fact, what occurs is that a brief instability bubble is formed, where the two eigenvalues collide, move off the imaginary axis, trace an approximate circle in the complex plane, and recombine on the imaginary axis in a ``reverse'' Krein collision. This brief instability bubble, which we call a Krein bubble, is also a consequence of the block matrix reduction, and is shown in schematic form in \cref{fig:KreinBubbleCartoon}. The radius of the Krein bubble in the complex plane and the value of $X$ at which the Krein bubble occurs can be computed using the block matrix reduction. Similar instability bubbles have been observed in other systems. As one example, they are found for dark soliton solutions of the discrete nonlinear Schr{\"o}dinger equation on a finite lattice as the coupling parameter is increased \cite{Johansson1999}; in that case, however, the instability bubbles disappear after a critical value of the coupling parameter is reached (\cite[Figure 2]{Johansson1999}).

This paper is organized as follows. In \cref{sec:KdV5}, we introduce a generalization of KdV5 as our motivating example. In \cref{sec:setup}, we set up the problem of interest in general terms as a Hamiltonian system in $2m$ dimensions which is reversible and translation invariant,
for which KdV5 (corresponding to $m=2$) is a special case. We also comment on extensions to higher order models, for which $m>2$.
We then present the main results of this paper, which concern the existence (\cref{sec:perexist}) and spectrum (\cref{sec:perstab}) of periodic multi-pulse solutions. This is then applied to the periodic single pulse and the periodic double pulse. In particular, we prove that Krein bubbles occur, and we give a formula for their radius in terms of fundamental constants associated with the system. In \cref{sec:numerics}, we present numerical results which provide verification for our theoretical work,
including timestepping simulations to demonstrate the dynamical consequence of the Krein bubble. 
The next sections contain proofs of the main results, after which we discuss our findings in \cref{sec:conclusions} and offer some directions for future work. 

\section{Background and motivation}\label{sec:KdV5}

The Kawahara equation, also known as a fifth-order KdV-type equation, is used as a model for water waves, magneto-acoustic waves, plasma waves, and other dispersive phenomena. This equation takes the general form
\begin{equation}\label{kawahara}
u_t + \alpha u_{xxx} + \beta u_{xxxxx} = \frac{\partial}{\partial x} f(u, u_x, u_{xxx}),
\end{equation}
where $u(x, t)$ is a real-valued function, the parameters $\alpha$ and $\beta$ are real with $\beta \neq 0$, and $f$ is a smooth function \cite{Bridges2002,Bridges2002a}. If $f$ is a variational derivative, then \cref{kawahara} is the Hamiltonian system $\partial_t u = \calJ \calE'(u)$, where $\calJ = \partial_x$ is skew-Hermitian,
\begin{equation}\label{kawaharaE}
\calE(u) = -\frac{1}{2}\int_{-\infty}^\infty 
\left( \frac{1}{2}u_{xx}^2 - \frac{1}{2}\alpha u_x^2 + h(u, u_x, u_{xx})\right)dx
\end{equation}
is the energy, and $f$ in \cref{kawahara} is the variational derivative of the term involving $h$ in \cref{kawaharaE} \cite{Bridges2002}. A prototypical example is
\[
u_t + \frac{2}{15}u_{xxxxx} - b u_{xxx} + 3 u u_x + 2 u_x u_{xx} + u u_{xxx} = 0,
\]
which is a weakly nonlinear long-wave approximation for capillary-gravity water waves \cite{Champneys1998,Champneys1997}. We will consider instead the simpler equation 
\begin{equation}\label{KdV5}
u_t = u_{xxxxx} + p u_{xxx} - 2 u u_x,
\end{equation}
which is a general form of the equation studied in \cite{Pelinovsky2007}. Writing \cref{KdV5} in a co-moving frame with speed $c$ by letting $\xi = x - ct$, equation \cref{KdV5} becomes
\begin{equation}\label{KdV5c}
u_t = \partial_x(u_{xxxx} + p u_{xx} + cu - u^2),
\end{equation}
where we have renamed the independent variable back to $x$. Localized traveling pulse solutions satisfy the 4th order ODE
\begin{equation}\label{KdV5eq4}
u_{xxxx} + p u_{xx} + c u - u^2 = 0,
\end{equation}
which is obtained from \cref{KdV5c} by integrating once. Equation \cref{KdV5eq4} is Hamiltonian, with conserved quantity
\begin{equation}\label{KdV5ham}
H(u, u_x, u_{xx}, u_{xxx}) = u_x u_{xxx} - \frac{1}{2}u_{xx}^2 + \frac{p}{2}u_x^2 + \frac{c}{2}u^2 - \frac{1}{3}u^3,
\end{equation}
which is obtained by multiplying \cref{KdV5eq4} by $u_x$ and integrating once. Letting
\begin{equation}\label{KdV5U}
U = (q_1, q_2, p_1, p_2) = (u, u_x, -u_{xxx} + u_x, u_{xx}),
\end{equation}
we can also write \cref{KdV5eq4} in standard Hamiltonian form as the first order system
\begin{equation}\label{KdV5ham2}
U' = F(U) = J \nabla \tilde{H}(U),
\end{equation}
where $J$ is the standard $4 \times 4$ symplectic matrix and
\begin{equation}
\tilde{H}(q_1, q_2, p_1, p_2) = \frac{1}{3}q_1^3 - \frac{1}{2}c q_1^2 + p_1 q_2 + \frac{p}{2}q_2^2 + \frac{1}{2}p_2^2.
\end{equation}
We have the following theorem concerning the existence of localized solutions to \cref{KdV5eq4}, which is a direct consequence of \cite{Groves1998}, reversibility, and the stable manifold theorem.

\begin{theorem}\label{KdV1pulse}
If $p < 2 \sqrt{c}$, then there exists a single-pulse solution $q(x)$ to \cref{KdV5eq4} which is an even function and decays exponentially to 0 at $\pm \infty$.
\end{theorem}

Linearization of \cref{KdV5eq4} about a solution $u(x)$ is the self-adjoint linear operator
\begin{equation}\label{KdV5hessian}
\calE''(u) = \partial_x^4 + p \partial_x^2 + c - 2 u^*,
\end{equation}
where $\calE''(u)$ is the Hessian of the energy. The rest state $u = 0$ corresponds to the equilibrium point $U = 0$ of the first order system \cref{KdV5ham2}. When $p < 2 \sqrt{c}$, this equilibrium is a hyperbolic saddle with 2-dimensional stable and unstable manifolds. The single pulse $q(x)$ corresponds to a homoclinic orbit connecting the stable and unstable manifolds of this equilibrium. If $-2 \sqrt{c} < p < 2 \sqrt{c}$, the eigenvalues of $DF(0)$ are a complex quartet $\pm \alpha_0 \pm \beta i_0$, and multi-modal homoclinic and periodic orbits exist which lie close to the primary homoclinic orbit \cite{SandstedeStrut}. We adapt Lin's method as in \cite{Sandstede1993, SandstedeStrut} to construct periodic multi-pulses ($n$-periodic solutions) by gluing together consecutive copies of the primary pulse end-to-end in a loop using small remainder functions. This provides not only an existence result but also estimates for these small remainder functions. 
As opposed to $n$-homoclinic solutions, for which the pulse tails are spliced together at $n-1$ locations, these $n$-periodic solutions require $n$ splices at the pulse tails, which provides an additional degree of freedom. 
For spectral stability, as in \cite{Sandstede1998}, we reduce the computation of the spectrum of the linearization of the PDE \cref{KdV5c} about a periodic $n$-pulse to a matrix equation. In contrast to \cite{Sandstede1998}, we obtain a $2n \times 2n$ block matrix, which encodes both the interaction eigenvalues and the essential spectrum eigenvalues near the origin.

\section{Mathematical Setup}\label{sec:setup}

\subsection{Hamiltonian PDE}\label{sec:HamPDE}

First, we define a Hamiltonian PDE which is reversible and translation invariant. This analysis follows Grillakis, Shatah, and Strauss \cite{Grillakis1987}. Let $X = H^{2m}(\R)$ for $m \geq 2$, and $Y = L^2(\R)$, and consider the PDE
\begin{equation}\label{genPDE}
u_t = \partial_x \calE'(u),
\end{equation}
where $u \in X$ and $\calE: X \subset Y \rightarrow \R$ is a smooth functional representing the conserved energy of the system. We take the following hypothesis regarding the energy $\calE(u)$.

\begin{hypothesis}\label{hyp:E}
The energy $\calE(u)$ has the following properties:
\begin{enumerate}[(i)]
\item $\calE(0) = 0$ and $\calE'(0) = 0$.
\item $\calE(u) = \calE(\rho(u))$, where $\rho: X \rightarrow X$ is the reversor operator $[\rho(u)](x) = u(-x)$.
\item $\calE(T(s)u) = \calE(u)$ for all $s \in \R$, where $\{T(s) : s \in \R \}$ is the one parameter group of unitary translation operators on $X$ defined by $[T(s)]u(\cdot) = u(\cdot - s)$.
\item $\calE'(u): X \rightarrow X$ is a differential operator of the form
\begin{equation}\label{Eprimeuform}
\calE'(u) = \partial_x^{2m}u - f(u, \partial_x u, \dots, \partial_x^{2m-1} u),
\end{equation}
where $f: \R^{2m} \rightarrow \R$ is smooth.
\end{enumerate}
\end{hypothesis}

\noi \cref{hyp:E}(ii) is reversibility, and \cref{hyp:E}(iii) is translation invariance. \cref{hyp:E}(iv) holds in applications such as KdV5, and lets us write the $2m$-th order ODE $\calE'(u) = 0$ as a first order system in $\R^{2m}$. 

\begin{remark}\label{remark:largerm}
Although we are most interested in the case where $m=2$, for which the Kawahara equation \cref{kawahara} and the fifth-order KdV model \cref{KdV5example} are specific examples, the theory is developed for general $m\geq 2$ so that it applies to higher order models as well. An example for $m=3$ is the seventh-order KdV equation (\cite[Chapter 15.10]{Wazwaz2009} and \cite[Section 8]{Chardard2})
\begin{align}\label{KdV7}
u_t &= \partial_x\left( -a u_{xxxxxx} + u_{xxxx} - u_{xx} + c u - 3 u^2\right),
\end{align}
which was introduced to study the KdV equation under singular perturbations (see also equation (24) in \cite{Pomeau1988}). The ninth-order KdV equation \cite[Chapter 15.10]{Wazwaz2009} corresponds to $m=4$. There has also been recent interest in nonlinear Schr\"{o}dinger models incorporating higher order dispersion terms \cite{Runge2020}.
\end{remark}

Differentiating the reversibility relation $\calE(u) = \calE(\rho(u))$ with respect to $u$,
\[
\calE'(u) = \rho^*( \calE'(\rho(u) ) ) = \rho( \calE'(\rho(u) ) ),
\]
since $\rho$ is self-adjoint. Differentiating the symmetry relation $\calE(T(s)u) = \calE(u)$ with respect to $u$,
\begin{align}
\calE'(u) &= T(s)^* \calE'(T(s)u) \label{Eprimesymm} \\
\calE''(u) &= T(s)^* \calE''(T(s)u) T(s). \label{EHessiansymm}
\end{align}
Differentiating the symmetry relation $\calE(T(s)u) = \calE(u)$ with respect to $s$ at $s = 0$, 
\begin{align*}
0 = \langle \calE'(u), T'(s) u \rangle|_{s = 0}
= \langle \calE'(u), T'(0) u \rangle
= \langle \calE'(u), \partial_x u \rangle
\end{align*}
for all $u \in X$, since $T'(0) = \partial_x$ is the infinitesimal generator of the translation group $T(s)$. There is an additional conserved quantity $\calQ: L^2(\R) \rightarrow \R$, given by
\begin{equation}\label{defV}
\calQ(u) = -\frac{1}{2} \int_{-\infty}^\infty u^2 dx,
\end{equation}
which represents charge in some applications. Traveling waves are solutions of \cref{genPDE} of the form $u(x, t) = T(ct)\phi(x) = \phi(x - ct)$. If $\phi$ satisfies the equilibrium equation $\calE'(\phi) = c \calQ'(\phi)$, then $T(ct)\phi(x)$ is a traveling wave \cite{Grillakis1987}. Since $\calQ'(\phi) = -\phi$, the equilibrium equation becomes
\begin{equation}\label{eqODE}
\calE'(\phi) + c \phi = 0.
\end{equation}
Without loss of generality, we will assume that $\calE'(\phi)$ does not contain any terms of the form $b\phi$ for $b$ constant, since that is accounted for by the $c \phi$ term in \cref{eqODE}.

We take the following hypothesis concerning the existence of traveling waves, which is similar to \cite[Assumption 2]{Grillakis1987}. In the next section, we will give a condition under which this hypothesis is satisfied. 
\begin{hypothesis}\label{hyp:cinterval}
There exists an open interval $(c_1, c_2) \subset \R$ and a $C^1$ map $c \mapsto \phi_c$ such that for every $c \in (c_1, c_2)$, $\calE'(\phi_c) + c \phi_c = 0$, i.e. $\phi_c$ is a traveling wave solution to \cref{genPDE}.
\end{hypothesis}

\noi The linearization of the PDE \cref{genPDE} about a traveling wave solution $\phi_c$ is the linear operator $\partial_x \calL(\phi_c)$, where $\calL(\phi_c)$ is the self-adjoint operator
\begin{equation}\label{PDElinearization}
\calL(\phi_c) = \calE''(\phi_c) + c,
\end{equation}
and $\calE''(\phi_c)$ is the Hessian of the energy $\calE(\phi_c)$. Differentiating \cref{eqODE} with respect to $x$ and with respect to $c$,
\begin{equation}\label{Lkernel}
\begin{aligned}
\calL(\phi_c) \partial_x \phi_c &= 0 \\
\calL(\phi_c) (-\partial_c \phi_c)& = \phi_c.
\end{aligned}
\end{equation}
Differentiating again with respect to $x$,
\begin{equation}\label{Ekernel}
\begin{aligned}
[\partial_x \calL(\phi_c)] \partial_x \phi_c &= 0 \\
[\partial_x \calL(\phi_c)](-\partial_c \phi_c) &= \partial_x \phi_c,
\end{aligned}
\end{equation}
thus the kernel of $\partial_x \calE''(\phi_c)$ has algebraic multiplicity at least 2 and geometric multiplicity at least 1.

\subsection{Spatial dynamics formulation}\label{sec:spatdym}

We reformulate the equilibrium equation \cref{eqODE} using a spatial dynamics approach by rewriting it as a first-order dynamical system in $\R^{2m}$ evolving in the spatial variable $x$. From this viewpoint, an exponentially localized traveling wave is a homoclinic orbit connecting a saddle point equilibrium at the origin to itself. Let $U = (u, \partial_x u, \dots, \partial_x^{2m-1} u)^T \in \R^{2m}$. Using \cref{hyp:E}(iv), equation \cref{eqODE} is equivalent to the first order system
\begin{equation}\label{genODE}
U'(x) = F(U(x); c),
\end{equation}
where $F: \R^{2m} \times \R \rightarrow \R^{2m}$ is smooth and is given by
\begin{equation}\label{defF}
F(u_1, u_2, \dots, u_{2m}; c) = 
\begin{pmatrix}
u_2 \\ u_3 \\ \vdots \\ f(u_1, u_2, \dots, u_{2m}) - c u_1
\end{pmatrix}.
\end{equation}
By reversibility,
\begin{equation}\label{genODErev}
\begin{aligned}
F(RU; c) &= -RF(U; c) \\
DF(RU; c) &= -RDF(U; c)R,
\end{aligned}
\end{equation}
where $R:\R^{2m} \rightarrow \R^{2m}$ is the standard reversor operator on $\R^{2m}$
\begin{equation}\label{reverserR2m}
R(u_1, u_2, \dots, u_{2m-1}, u_{2m}) = (u_1, -u_2, \dots, u_{2m-1}, -u_{2m}).
\end{equation}
First, we assume that \cref{genODE} is a conservative system.
\begin{hypothesis}\label{hyp:H}
There exists a smooth function $H: \R^{2m} \times \R \rightarrow \R$ such that 
\begin{enumerate}[(i)]
\item $H(0; c) = 0$ for all $c$.
\item $\nabla_U H(U; c) = 0$ if and only if $F(U; c) = 0$.
\item For all $U \in \R^{2m}$ and all $c$, $\langle F(U; c), \nabla_U H(U; c) \rangle = 0$.
\end{enumerate}
\end{hypothesis}

\noi It follows from \cref{hyp:H} that $H$ is conserved along solutions to \cref{genODE}. Since $F(0; c) = 0$ for all $c$, the rest state $U = 0$ is an equilibrium of \cref{genODE} for all $c$. The next hypothesis addresses the hyperbolicity of this equilibrium. Although the eigenvalue pattern described in \cref{hyp:hypeq} is not necessary for the existence of a homoclinic orbit solution, it is a sufficient condition for the existence of multi-pulse and periodic multi-pulse solutions.

\begin{hypothesis}\label{hyp:hypeq}
For a specific $c_0 > 0$, $U = 0$ is a hyperbolic equilibrium of \cref{genODE}. Furthermore, the spectrum of $DF(0; c_0)$ contains a quartet of simple eigenvalues $\pm \alpha_0 \pm \beta_0 i$, where $\alpha_0, \beta_0 > 0$, and for any other eigenvalue $\nu$ of $DF(0; c_0)$, $|\text{Re }\nu| > \alpha_0$.
\end{hypothesis}

We note that localized pulse solutions will have tails which are exponentially decaying with approximate rate $\alpha_0$, and are oscillatory with approximate frequency $\beta_0$.

\begin{figure}
\begin{center}
\includegraphics[width=8cm]{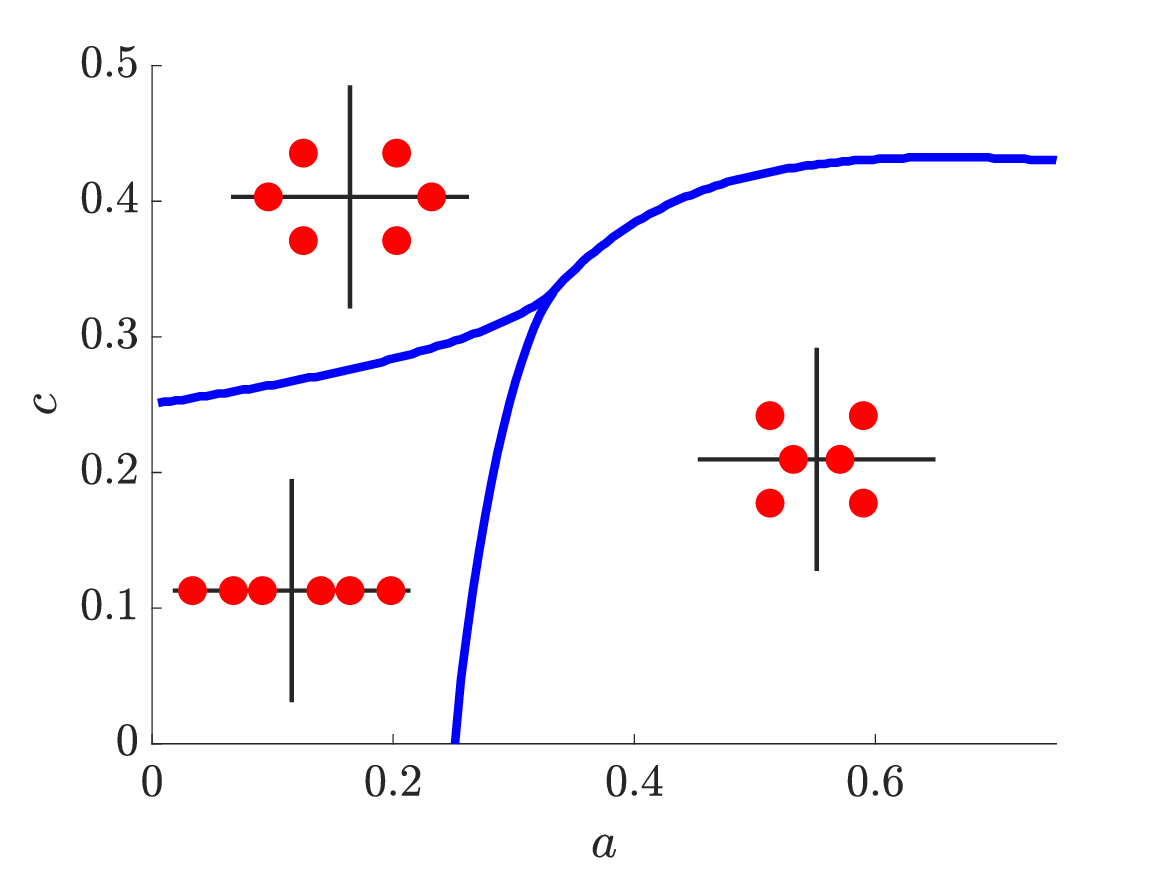}
\end{center}
\caption{Eigenvalue pattern of $DF(0; c)$ as the wavespeed $c$ and the parameter $a$ are varied for the 7th order KdV equation \cref{KdV7}.}
\label{fig:KdV7spec}
\end{figure}

\begin{remark}\label{remark:hypsatisfied}
For the 5th order KdV equation \cref{KdV5example}, corresponding to $m=2$, the spectrum of $DF(0; c)$ is the quartet of eigenvalues 
\[
\lambda = \pm \sqrt{ \frac{1 \pm \sqrt{1 - 4c } }{ 2} }.
\] 
For $c > 1/4$, this is a complex quartet $\pm \alpha_0 \pm \beta_0 i$, thus \cref{hyp:hypeq} is satisfied. For the 7th order KdV equation \cref{KdV7}, corresponding to $m=3$, the spectrum of $DF(0; c)$ comprises six eigenvalues, one pair of which is always real (see \cref{fig:KdV7spec}, as well as the first quadrant of \cite[Figure 3]{Chardard2}). \cref{hyp:hypeq} is satisfied in the upper region of \cref{fig:KdV7spec}. Note that in the lower right region of \cref{fig:KdV7spec}, there is a complex quartet of eigenvalues, but since the real pair of eigenvalues lies inside this complex quartet, \cref{hyp:hypeq} is not satisfied.
\end{remark}

We now address the existence of a primary pulse solution, which is a symmetric homoclinic orbit connecting the unstable manifold $\tilde{W}^u(0; c_0)$ and the stable manifold $\tilde{W}^s(0; c_0)$ of the rest state equilibrium $U = 0$. Both of these manifolds have dimension $m$ by reversibility. Since, in general, the existence of such a solution is unknown, we take the existence of a primary pulse solution for a specific wavespeed $c_0$ as a hypothesis. For specific systems, such as KdV5, the existence of a primary pulse solution has been proved (see, for example, \cref{KdV1pulse} above).

\begin{hypothesis}\label{Qexistshyp}
For the same $c_0$ as in \cref{hyp:hypeq}, there exists a homoclinic orbit solution $Q_1(x; c_0) = (q(x; c_0), \partial_x q(x; c_0), \dots, \partial_x^{2m-1}q(x; c_0))^T\in \tilde{W}^s(0; c_0) \cap \tilde{W}^u(0; c_0) \subset H^{-1}(0; c_0)$ to \cref{genODE}. In addition,
\begin{enumerate}[(i)]
\item $Q_1(0; c_0) \neq 0$.
\item $\nabla_U H(Q_1(0; c_0); c_0) \neq 0$.
\item $Q_1(x; c_0)$ is symmetric with respect to the reversor operator \cref{reverserR2m}, i.e. $Q_1(-x; c_0) = R Q_1(x; c_0)$.
\end{enumerate}
\end{hypothesis}

\noi It follows from \cref{Qexistshyp} that the first component $q(x; c_0)$ is a symmetric, exponentially localized traveling wave solution to \cref{genPDE}. In order to prove the  existence of homoclinic orbits $Q_1(x; c)$ for $c$ near $c_0$, we take the following additional hypothesis.

\begin{hypothesis}\label{hyp:transverse}
The stable manifold $\tilde{W}^s(0; c_0)$ and the unstable manifold $\tilde{W}^u(0; c_0)$ intersect transversely in $H^{-1}(0; c_0)$ at $Q_1(0; c_0)$.
\end{hypothesis}

\noi Using \cref{hyp:transverse} and a dimension-counting argument, we obtain the nondegeneracy condition
\begin{equation}\label{nondegencond}
T_{Q_1(0; c_0)}\tilde{W}^s(0; c_0) \cap T_{Q_1(0; c_0)}\tilde{W}^u(0; c_0) = \R Q_1'(0; c_0).
\end{equation}
We then have the following existence theorem. The proof is given in \cref{sec:transverseintproof}.

\begin{theorem}\label{transverseint}
Assume \cref{hyp:H}, \cref{hyp:hypeq}, and \cref{hyp:transverse}. Then there exists $\delta_0 > 0$ such that for $c \in (c_0 - \delta_0, c_0 + \delta_0)$, the stable and unstable manifolds $\tilde{W}^s(0; c)$ and $\tilde{W}^u(0; c)$ have a one-dimensional transverse intersection in $H^{-1}(0; c)$, which is a homoclinic orbit $Q_1(x; c)$. Furthermore, $Q_1(-x; c) = R Q_1(x; c)$, the map $c \rightarrow Q_1(x; c)$ is smooth, and $\partial_c Q_1(x; c)$ is exponentially localized, i.e. for any $\epsilon > 0$ there exists $\delta_1 > 0$ with $\delta_1 \leq \delta_0$ such that for $c \in (c_0 - \delta_1, c_0 + \delta_1)$,
\begin{equation}\label{Qcbound}
|\partial_c Q_1(x; c)| \leq C e^{-(\alpha_0 - \epsilon)|x|},
\end{equation}
where $\alpha_0$ is defined in \cref{hyp:hypeq}.
\end{theorem}

Finally, as in \cite{Grillakis1987}, we define the scalar 
\[
d(c) = \calE(q(x, c)) - \omega\calQ(q(x, c)).
\]
By \cite{Bona1987,Grillakis1987}, the traveling wave $q(x, c)$ is orbitally stable if $d''(c) > 0$, where
\begin{align}\label{ddoubleprime}
d''(c) = \langle \calQ'(q(x, c)), \partial_c q(x, c) \rangle
= \int_{-\infty}^\infty q(x, c) \partial_c q(x, c) dx.
\end{align}
This can be computed numerically, and we take this stability criterion as a hypothesis.

\begin{hypothesis}\label{hyp:dccpos}
For each $c \in (c_0 - \delta_0, c_0 + \delta_0)$, where $\delta_0$ is defined in \cref{transverseint}, $d''(c) > 0$.
\end{hypothesis}

\noi From this point on, we will fix a speed $c \in (c_0 - \delta_0, c_0 + \delta_0)$ and suppress the dependence on $c$ for simplicity of notation.

\subsection{Eigenvalue problem}\label{sec:EVP}

Let $U^*(x) = (u^*(x), \partial_x u^*(x), \dots, \partial_x^{2m-1}u^*(x) )^T$ be any solution to \cref{genODE}, so that $u^*(x)$ is a traveling wave solution to \cref{genPDE}.
Then $u^*(x)$ also solves the equation $\partial_x(\calE'(u) + cu) = 0$, which is equivalent to the system
\begin{equation}\label{eqsystem2}
\begin{aligned}
\calE'(u) + cu &= k \\
\partial_x k &= 0.
\end{aligned}
\end{equation}
Using a spatial dynamics approach, we rewrite \cref{eqsystem2} as the first order dynamical system in $\R^{2m+1}$ 
\begin{equation}\label{spsystem2}
\begin{pmatrix}
U \\ k
\end{pmatrix}'(x) =
\begin{pmatrix} 
F(U(x)) + k e_{2m} \\ 0
\end{pmatrix},
\end{equation}
where $e_{2m} = (0, \dots, 0, 1)^T \in \R^{2m}$ is the standard unit vector. We similarly reformulate the PDE eigenvalue problem $\partial_x \calL(u^*) v = \lambda v$ as the system 
\begin{equation}\label{genPDEeig2}
\begin{aligned}
\calL(u^*) v &= k \\
\partial_x k &= \lambda v.
\end{aligned}
\end{equation}
This is equivalent to the first order dynamical system in $\C^{2m+1}$
\begin{equation}\label{PDEeigsystem}
V'(x) = A(U^*(x))V(x) + \lambda B V(x),
\end{equation}
where $V(x) \in C^0(\R,\C^{2m+1})$, and $A(U^*(x))$ and $B$ are the $(2m+1) \times (2m+1)$ matrices
\begin{equation}\label{defAB}
A(U^*(x)) = 
\begin{pmatrix}
DF(U^*(x)) & e_{2m}\\
0 & 0
\end{pmatrix}, \qquad
B = \begin{pmatrix}0 & 0 & \cdots & 0 \\ & 
\vdots & \vdots & \\0 & 0 & \cdots & 0 \\ 1 & 0 & \cdots & 0 \end{pmatrix}.
\end{equation}
$A(0)$ has a one-dimensional kernel, which is characterized in the following lemma. 

\begin{lemma}\label{eigA0lemma}
The matrix $A(0)$ has a simple eigenvalue at 0 and a quartet of eigenvalues $\pm \alpha_0 \pm \beta_0 i$. For any other eigenvalue $\nu$ of $A(0)$, $|\text{Re }\nu| > \alpha_0$. The kernel of $A(0)$ is spanned by $V_0$ and the kernel of $A(0)^*$ is spanned by $W_0$, where
\begin{equation}\label{V0W0}
V_0 = \left(\frac{1}{c}, 0, \dots, 0, 1\right)^T, \quad
W_0 = (0, 0, \dots, 0, 1)^T,
\end{equation}
and $\langle V_0, W_0 \rangle = 1$. The projection on $\R V_0$ is given by $P_{V_0} = \langle W_0, \cdot \rangle$.
\begin{proof}
Let $p_1(\nu)$ and $p_2(\nu)$ be the characteristic polynomials of $DF(0)$ and $A(0)$. Since
\begin{align*}
p_2(\nu) &= \det(A(0) - \nu I) = -\nu \det(DF(0) - \nu I) = -\nu p_1(\nu),
\end{align*}
$A(0)$ has the same eigenvalues as $DF(0)$ as well as an additional eigenvalue at 0, thus part (i) follows from \cref{hyp:hypeq}. The kernel eigenvectors $V_0$ and $W_0$ and the projection $P_{V_0}$ can be verified directly.
\end{proof}
\end{lemma}

Since $A(0)$ is non-hyperbolic, the rest state at the origin $(U, k) = (0, 0)$ is a non-hyperbolic equilibrium of \cref{spsystem2}, and the results of \cite{Sandstede1998} do not apply. Let $W^s(0)$, $W^u(0)$, and $W^c(0)$ be the stable, unstable, and center manifolds of the equilibrium at the origin. By reversibility, $\dim W^s(0) = m$, $\dim W^u(0) = m$, and $\dim W^c(0) = 1$. Let $Q_1(x)$ be the primary pulse solution from \cref{transverseint}, and define
\begin{equation}\label{eq:defQx}
Q(x) = (Q_1(x), 0).
\end{equation}
The associated variational and adjoint variational equations are
\begin{align}
V'(x) &= A(Q(x)) V(x) \label{vareq2} \\
W'(x) &= -A(Q(x))^* W(x) \label{adjvareq2},
\end{align}
and $Q'(x)$ is an exponentially localized solution to \cref{vareq2}. Since $Q(x)$ is exponentially localized, $\R Q'(0) \subset T_{Q(0)}W^s(0) \cap T_{Q(0)}W^u(0)$. It follows from \cref{hyp:transverse} that these are in fact equal.

\begin{lemma}\label{nondegenlemma}
We have the nondegeneracy condition
\begin{equation}\label{nondegen2}
T_{Q(0)}W^s(0) \cap T_{Q(0)}W^u(0) = \R Q'(0).
\end{equation}
\begin{proof}
If the intersection were more than one-dimensional, there would exist another exponentially localized solution $V(x) = (v_1, \dots, v_{2m}, v_{2m+1})^T$ to \cref{vareq2}. By the definition of $A(Q(x))$, $v_{2m+1}$ is a constant, which must be 0 since $V(x)$ is exponentially localized. Then $(v_1, \dots, v_{2m})^T$ would be an exponentially localized solution to $V'(x) = DF(Q(x)) V(x)$, which contradicts the nondegeneracy condition \cref{nondegencond}.
\end{proof}
\end{lemma}

\noi Using \cref{nondegen2}, we can decompose the tangent spaces of the stable and unstable manifolds at $Q(0)$ as
\begin{equation}\label{TQ0decomp}
\begin{aligned}
T_{Q(0)}W^s(0) &= \R Q'(0) \oplus Y^+ \\
T_{Q(0)}W^u(0) &= \R Q'(0) \oplus Y^-.
\end{aligned}
\end{equation}
Since $\dim \R Q'(0) \oplus Y^+ \oplus Y^- = 2m-1$, we need two more directions to span $\R^{2m+1}$. We obtain these from the following lemma.

\begin{lemma}\label{varadjsolutions}
Let $Q(x)$ be defined by \cref{eq:defQx}. Then we have the following bounded solutions to the variational equation \cref{vareq2} and the adjoint variational equation \cref{adjvareq2}:
\begin{enumerate}[(i)]
	\item There are two linearly independent, bounded solutions to the variational equation \cref{vareq2}, which are given by $Q'(x)$ and $V^c(x)$. $V^c(x) \rightarrow V_0$ as $|x| \rightarrow \infty$, where $V_0$ is defined by \cref{V0W0}, and $V^c(-x) = R V^c(x)$, where $R$ is the standard reversor operator. Furthermore, $V^c = (\tilde{V}^c, 1)$, where $\tilde{V}^c$ solves the equation $\tilde{V}^c(x)' = DF(Q(x)) \tilde{V}^c(x) + e_{2m}$. Any other bounded solution to \cref{vareq2} is a linear combination of these.

	\item There are two linearly independent, bounded solutions to the adjoint variational equation \cref{adjvareq2}, which are given by $\Psi(x)$ and $W_0$. $\Psi(x)$ is the exponentially localized solution
	\begin{equation}\label{psicomponents}
	\Psi(x) = (\nabla H(Q(x)), q(x))^T, 
	\end{equation}
    where the conserved quantity $H$ is defined in \cref{hyp:H},
	and $\Psi(-x) = R \Psi(x)$. The constant solution $W_0$ is defined by \cref{V0W0}. Any other bounded solution to \cref{adjvareq2} is a linear combination of these.
\end{enumerate}
\begin{proof}
For part (i), the existence of $V^c(x)$ is a consequence of the geometry of the system, and will be proved below after \cref{lemma:Vpm}. The equation $V^c(x)' = A(Q(x))V^c(x)$ then reduces to $\tilde{V}^c(x)' = DF(Q(x)) \tilde{V}^c(x) + e_{2m}$. For part (ii), equation \cref{adjvareq2} can be written in block form as
\[
W'(x) = - 
\begin{pmatrix}DF(Q(x))^* & 0 \\ e_{2m}^T & 0 \end{pmatrix} W(x),
\]
for which $W_0 = (0, \dots, 0, 1)^T$ is a constant solution. Using \cref{psiform} below, $\Psi(x) = ( \nabla H(Q(x)), q(x) )^T$ is an exponentially localized solution to \cref{adjvareq2}. 
\end{proof}
\end{lemma}

\begin{remark}\label{remark:computeVc}
Let $v^c(x)$ be the first component of $V^c(x)$. Then $v^c$ is a formal solution to $\calL(q) v^c = 1$, which provides a convenient way of computing $V^c(x)$ numerically.
\end{remark}

\noi By \cref{eigadjoint} below, $\Psi(0)$ and $W_0$ are perpendicular to $\R Q'(0) \oplus Y^+ \oplus Y^-$ at $x = 0$, thus we can decompose $\R^{2m+1}$ as 
\begin{equation}\label{DSdecomp}
\R^{2m+1} = \R Q'(0) \oplus Y^+ \oplus Y^- \oplus \R \Psi(0) \oplus \R W_0.
\end{equation}

\section{Existence of periodic multi-pulses}\label{sec:perexist}

\begin{figure}
\begin{center}
\includegraphics[width=11cm]{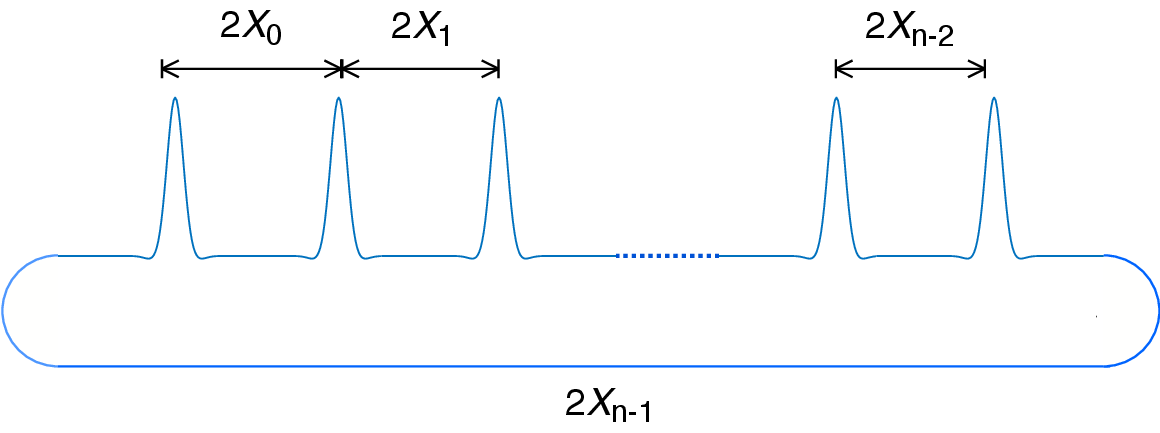}
\end{center}
\caption{Schematic showing the construction of a periodic $n-$pulse solution from the primary pulse.}
\label{fig:permultipulse}
\end{figure}

In this section, we prove the existence of periodic multi-pulse solutions to \cref{genODE}, which are multi-modal periodic orbits that remain close to the primary homoclinic orbit. Heuristically, we construct a periodic multi-pulse by gluing together multiple copies of the primary pulse end-to-end in a loop (\cref{fig:permultipulse}).

A periodic $n$-pulse can be described by the $n$ pulse distances $X_0, X_1, \dots, X_{n-1}$.
The distances between consecutive pulses are $2X_i$, as shown in \cref{fig:permultipulse}. 
The period of the orbit is $2X$, where $X = X_0 + \dots + X_{n-1}$. A periodic $n$-pulse requires one more length parameter than an $n$-pulse on the real line, since we need one more connection to ``close the loop''.  Rather than describing a periodic multi-pulse by the ``physical'' pulse distances $X_i$, we will use a parameterization which is more mathematically convenient and captures the underlying geometry necessary for a periodic $n$-pulse to exist. This parameterization is an adaptation of that in \cite{SandstedeStrut,Sandstede1998} to the periodic case. Let
\begin{equation}\label{defrho}
\rho = \frac{\beta_0}{\alpha_0}, \quad p^* = \arctan \rho,
\end{equation}
where $\alpha_0$ and $\beta_0$ are defined in \cref{hyp:hypeq}. Define the set
\begin{align}
\mathcal{R} &= \left\{ \exp\left(-\frac{2 m \pi}{\rho}\right) : m \in \N_0 \right\} \cup \{ 0 \},
\end{align}
which is a complete metric space. We will use $r \in \mathcal{R}$ as a scaling parameter. The parameterization is defined as follows.

\begin{figure}
\begin{center}
\includegraphics[width=10cm]{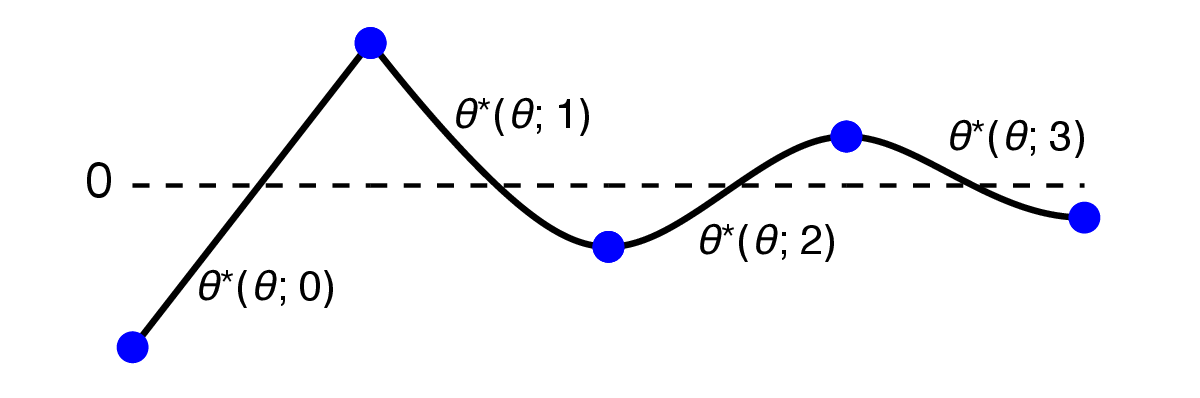}
\end{center}
\caption[Schematic for $\theta^*(\theta; m)$]{Schematic showing the first four functions $\theta^*(\theta; m)$ plotted consecutively. The region between each pair of blue dots corresponds to the domain $\theta = [-\pi + p^*, p^*]$.}
\label{fig:thetastarcartoon}
\end{figure}

\begin{definition}\label{def:perparam}
For $n \geq 2$, a \emph{periodic parameterization} of a periodic $n$-pulse is a sequence of $n+1$ parameters $(m_0, m_1, \dots, m_{n-1}, \theta)$, where $\theta \in (-\pi + p^*, p^*]$ and the $m_i$ are nonnegative integers which are chosen so that
\begin{enumerate}[(i)]
\item at least one of the $m_i \in \{0, 1\}$.
\item $m_{n-1} \geq m_i$ for $i = 0, \dots, n-2$.
\end{enumerate}
\end{definition}
\noi The selection of $m_{n-1}$ as the largest of the nonnegative integers $m_i$ is made for convenience of notation, and to allow the periodic parameterization to be unique. Since we are on a periodic domain, there is no loss of generality. The physical pulse distances $X_i$ are determined by the periodic parameterization and by the scaling parameter $r$. If $r = \exp\left(-\frac{2 m \pi}{\rho}\right)$, then
\begin{align*}
X_i &= \frac{1}{2 \beta_0}\big( (2 m + m_i)\pi + \theta^*(\theta; m_{n-1} - m_i)\big) + \tilde{L} \qquad\qquad i = 0, \dots, n-2  \\
X_{n-1} &= \frac{1}{2 \beta_0}\big( (2 m + m_{n-1})\pi + \theta \big) + \tilde{L},
\end{align*}
where $\tilde{L}$ is a constant. 
The functions $\theta^*(\theta; m): [-\pi + p^*, p^*] \rightarrow \R$ are defined for all nonnegative integers $m$, are continuous in $\theta$, and have the following properties:
\begin{enumerate}[(i)]
\item $\theta^*(0; m) = 0 \text{ for all } m$.
\item $|\theta^*(\theta; m)| \leq |\theta|$.
\item $|\theta^*(\theta; m)| \leq C \exp\left(-\frac{m \pi}{\rho} \right)$.
\item $\theta^*(\theta; 0) = \theta $.
\item $\theta^*(p^*; m) = \theta^*(-\pi+p^*; m+1)$.
\end{enumerate}
The last property is a matching condition which ``links up'' the parameterizations corresponding to adjacent $m_i$. \cref{fig:thetastarcartoon} shows a schematic of the first four functions $\theta^*(\theta; m)$ plotted consecutively to illustrate these properties. Together with the restriction of $\theta$ to the half-open interval $\theta \in (-\pi + p^*, p^*]$ in \cref{def:perparam}, these guarantee that each periodic parameterization corresponds to a unique periodic multi-pulse. The proof that the functions $\theta^*(\theta; m)$ exist and have these properties is given in \cref{thetaparamlemma} below.

We can now state the main theorem of this section, which gives conditions for the existence of periodic multi-pulses. The requirement that the scaling parameter $r$ be sufficiently small means that the individual pulses must be well-separated. The proof is given in \cref{sec:existproof}. 

\begin{theorem}[Existence of $n$-periodic solutions]\label{th:perexist}
Assume Hypotheses \ref{hyp:E}, \ref{hyp:H}, \ref{hyp:hypeq}, \ref{Qexistshyp}, and \ref{hyp:transverse}. Let $Q_1(x)$ be the transversely constructed, symmetric primary pulse solution to \cref{genODE} from \cref{transverseint}. For any periodic parameterization $(m_0, \dots, m_{n-1}, \theta)$ with $\theta \notin \{-\pi + p^*, p^* \}$, there exists $r_* = r_*(m_0, \dots, m_{n-1}, \theta) > 0$ with the following property.
For any $r \in \mathcal{R}$ with $r \leq r_*$, there exists a periodic $n$-pulse solution $U(x) = U(x; m_0, \dots, m_{n-1}, \theta, r)$ to \cref{genODE}. The distances between consecutive copies of $Q_1(x)$ in $U(x)$ are given by $2X_i$, where the pulse distances $X_i$ are
\begin{equation}\label{Xi}
\begin{aligned}
	X_i(r; m_i, m_{n-1},\theta) &= \frac{1}{2 \alpha_0} |\log r| + \frac{1}{2\beta_0} t_i(r; m_i,m_{n-1}, \theta) + \tilde{L} &&\qquad i = 0, \dots, n-2 \\
	X_{n-1}(r; m_{n-1}, \theta) &= \frac{1}{2 \alpha_0} |\log r| + \frac{1}{2 \beta_0}( m_{n-1}\pi + \theta ) + \tilde{L}.
\end{aligned}
\end{equation}
The functions $t_i(r; m_i, m_{n-1}, \theta): \mathcal{R} \rightarrow \R$ are continuous in $r$ with 
\[
t_i(0; m_i, \theta) = m_i \pi + \theta^*(\theta; m_{n-1} - m_i) = m_i \pi + \mathcal{O}\left( e^{-\frac{\pi}{\rho}(m_{n-1} - m_i)} \right),
\]
and $\tilde{L}$ is a constant. Estimates for $U(x)$ in terms of the primary pulse $Q_1(x)$ are given below in \cref{solvewithjumps}.
\end{theorem}

\begin{remark}
It follows from the proof of \cref{2pulsebifurcation} that periodic single pulse solutions exist on the periodic domain $[-X, X]$ for all sufficiently large $X$ (see \cref{corr:1pexists}).
These are single-loop periodic orbits which lie close to the primary homoclinic orbit.
\end{remark}

\begin{figure}
\begin{center}
\includegraphics[width=9cm]{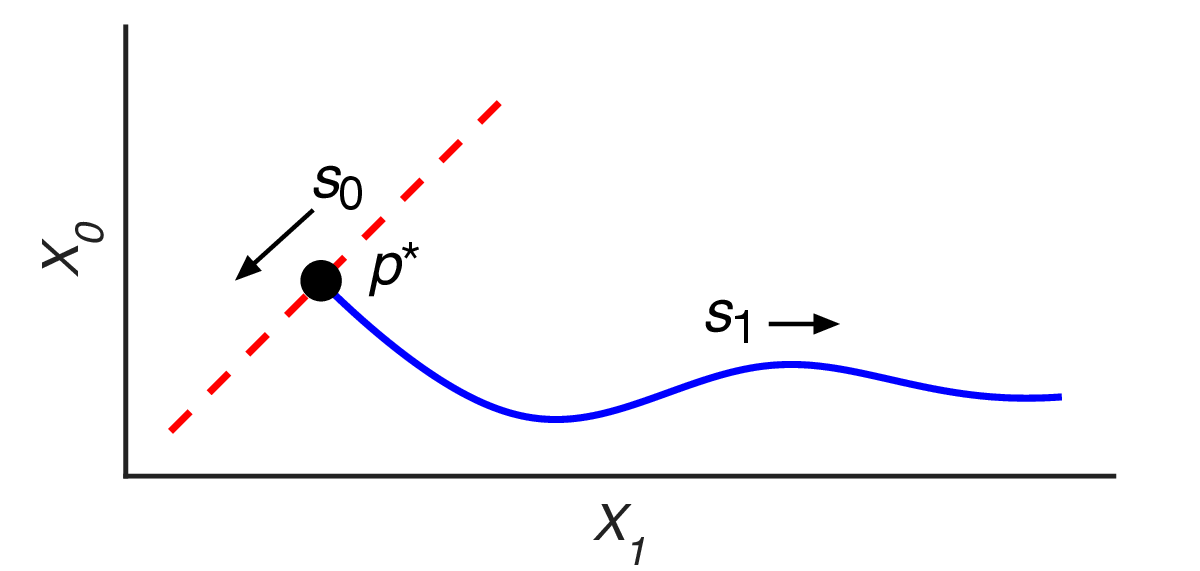}
\end{center}
\caption[Parameterization for periodic 2-pulses]{Parameterization for periodic 2-pulses. Red dashed line represents symmetric periodic 2-pulses, and blue solid line represents asymmetric periodic 2-pulses. Parameters $s_0$ and $s_1$ increase in the direction of the arrow.}
\label{fig:2pitch}
\end{figure}

The condition that $\theta \notin \{-\pi + p^*, p^* \}$ in \cref{th:perexist} is used to avoid bifurcation points which arise in the construction. For periodic 2-pulses, we can use the symmetry of the solutions and the reversibility of the system to give a complete bifurcation picture. In the next theorem, we show that for periodic 2-pulses, asymmetric solutions ($X_0 \neq X_1$) bifurcate from symmetric solutions ($X_0 = X_1$) in a series of pitchfork bifurcations (\cref{fig:KdV5periodicMP}, center panel). 
The symmetric 2-pulse solutions ($X_0 = X_1$) correspond to periodic single-pulse solutions with the period $2X_0$ repeated twice.
The parameterization in \cref{2pulsebifurcation}, which is shown in \cref{fig:2pitch}, is different from that in \cref{th:perexist}. The proof is given in \cref{sec:existproof}. 

\begin{theorem}\label{2pulsebifurcation}
Assume Hypotheses \ref{hyp:E}, \ref{hyp:H}, \ref{hyp:hypeq}, \ref{Qexistshyp}, and \ref{hyp:transverse}. Let $Q_1(x)$ be the transversely constructed, symmetric primary pulse solution to \cref{genODE} from \cref{Qexistshyp}. Then there exists $r_* > 0$ such that for all $r \in \mathcal{R}$ with $r \leq r_*$ and $m_0 \in \{0, 1\}$,
\begin{enumerate}[(i)]
	\item There exists a family of symmetric periodic 2-pulses $\tilde{Q}_2(x; m_0, s_0, r)$ parameterized by $s_0 \in [0, \pi)$. The pulse distances $\tilde{X}_i$ are given by
	\begin{equation}\label{2psymmdist}
		\tilde{X}_0(r, s_0) = \tilde{X}_1(r, s_0) = \frac{1}{2 \alpha_0} |\log r| + \frac{1}{2\beta_0} (m_0 \pi + s_0) + \tilde{L}.
	\end{equation}
	\item There exists a family of asymmetric periodic 2-pulses $Q_2(x; m_0, s_1, r)$ with pulse distances $X_1 > X_0$ parameterized by $s_1 \in [p^*, \infty)$. The pulse distances $X_i$ are given by
	\begin{equation}\label{2pasymmdist}
	\begin{aligned}
		X_0(r, m_0, s_1) &= \frac{1}{2 \alpha_0} |\log r| + \frac{1}{2\beta_0} t_0(r; m_0, s_1) + \tilde{L} \\
		X_1(r, s_1) &= \frac{1}{2 \alpha_0} |\log r| + \frac{1}{2\beta_0} s_1 + \tilde{L}, 
	\end{aligned}
	\end{equation}
	where $\tilde{L}$ is a constant, $t_0(r; m_0, s_1)$ is continuous in $r$ and $s_1$, $t_0(0; m_0, k \pi) = m_0 \pi$ for all nonnegative integers $k$, and 
	\begin{align}\label{t0est}
	t_0(0; m_0, s_1) = m_0 \pi + \mathcal{O}\left(e^{-\frac{1}{\rho} s_1 }\right).
	\end{align}

	\item The two families meet at a pitchfork bifurcation when $s_0 = p^*(m_0; r)$ and $s_1 = p^*$, where $p^*$ is defined in \cref{defrho}. The function $p^*(m_0; r)$ is continuous in $r$, and $p^*(m_0; r) \rightarrow p^*$ as $r \rightarrow 0$.
\end{enumerate}
\end{theorem}

\noi We note that in \cref{2pulsebifurcation}(ii), $X_1 > X_0$, which gives us the lower arms of the pitchforks in \cref{fig:2pitch}. For the upper arms, we swap $X_0$ and $X_1$, which we can do by symmetry.

\section{Spectrum of periodic multi-pulses}\label{sec:perstab}

We now locate the spectrum of the periodic multi-pulses which we constructed in the previous section. Let $Q_n(x) = (q_n(x), \partial_x q_n(x), \dots, \partial_x^{2m} q_n(x), 0)$ be any periodic $n$-pulse solution constructed according to \cref{th:perexist} on periodic domain $[-X, X]$. It is natural to pose the PDE eigenvalue problem \cref{genPDEeig2} on the space of periodic functions $H^{2m}_{\text{per}}[-X,X]$, where
\[
H^{2m}_{\text{per}}[-X,X] = \left\{ f \in H^{2m}(\R) : f^{(k)}(-X) = f^{(k)}(X) \text{ for } k = 0, 1, \dots, 2m \right\},
\]
although we note that in doing so, we are restricting ourselves to co-periodic perturbations. By \cref{Ekernel}, the linear operator $\partial_x \calL(q_n)$ has a kernel with algebraic multiplicity at least 2 and geometric multiplicity at least 1. In the next lemma, we show that $\partial_x \calL(q_n)$ has another kernel eigenfunction on $H^{2m}_{\text{per}}[-X,X]$. 

\begin{lemma}\label{qnkernel}
The linear operator $\partial_x \calL(q_n)$ posed on $H^{2m}_{\text{per}}[-X,X]$ has a kernel eigenfunction $v_n^c$, which is a solution to $\calL(q_n)v_n^c = 1$.
\begin{proof}
Since $1 \in H^{2m}_{\text{per}}[-X,X]$, $\calL(q_n) 1 = c \neq 0$. Since $\calL(q_n)$ is self-adjoint, for any $v \in \ker \calL(q_n)$,
\begin{align*}
\langle 1, v \rangle = \frac{1}{c} \langle \calL(q_n) 1, v \rangle
= \frac{1}{c} \langle  1, \calL(q_n)^* v \rangle = \frac{1}{c} \langle  1, \calL(q_n) v \rangle = 0,
\end{align*}
thus $1 \perp \ker \calL(q_n)^*$. By the Fredholm alternative, the equation $\calL(q_n) v_n^c = 1$ has a solution $v_n^c \in H^{2m}_{\text{per}}[-X,X]$. Differentiating with respect to $x$, $\partial_x \calL(q_n) v_n^c = 0$.
\end{proof}
\end{lemma}

\noi Using the same formulation as in \cref{sec:EVP}, the PDE eigenvalue problem \cref{genPDEeig2} on $H^{2m}_{\text{per}}[-X,X]$ is equivalent to the first order system with periodic boundary conditions
\begin{equation}\label{PDEeigsystemper3}
\begin{aligned}
V'(x) &= A(Q_n(x))V(x) + \lambda B V(x) \\
V(-X) &= V(X),
\end{aligned}
\end{equation}
where $V(x) \in C^0(\R,\C^{2m+1})$. In next lemma, we show that for small $\lambda$, the constant matrix $A(0) + \lambda B$ has a simple eigenvalue $\nu(\lambda)$ near 0.

\begin{lemma}\label{nulambdalemma}
There exists $\delta_1 > 0$ such that for $|\lambda| < \delta_1$, the matrix $A(0) + \lambda B$ has a simple eigenvalue $\nu(\lambda)$. Furthermore, $\nu(\lambda)$ is smooth in $\lambda$, $\nu(0) = 0$, $\nu'(0) = 1/c$, and for $|\lambda| < \delta_1$,
\begin{equation}\label{nulambda}
\nu(\lambda) = \frac{1}{c} \lambda + \mathcal{O}(|\lambda|^3).
\end{equation}
In addition, $\nu(-\lambda) = -\nu(\lambda)$ and $\nu(\overline{\lambda}) = \overline{\nu(\lambda)}$.
\begin{proof}
Let $p(\nu; \lambda)$ be the characteristic polynomial of $A(0) + \lambda B$. Since $p(0; 0) = 0$ and $\partial_\nu p(0; 0) = -c \neq 0$, by the implicit function theorem, there exists $\delta_1 > 0$ and a smooth function $\nu(\lambda)$ with $\nu(0) = 0$ such that for $|\lambda| < \delta_1$, $\nu(\lambda)$ is the unique solution to $p(\nu; \lambda) = 0$. The derivative $\nu'(0) = 1/c$ also follows from the implicit function theorem. By reversibility, $p(\nu; \lambda)$ only involves odd powers of $\nu$, thus $p(\nu; \lambda) = 0$ if and only if $p(-\nu; -\lambda) = 0$. Since the solution $\nu(\lambda)$ is unique, $\nu(-\lambda) = -\nu(\lambda)$. Conjugate symmetry follows similarly since $p(\nu; \lambda) = 0$ if and only if $p(\overline{\nu}; \overline{\lambda}) = 0$. Equation \cref{nulambda} follows from a Taylor expansion of $\nu(\lambda)$ about $\lambda = 0$. 
\end{proof}
\end{lemma}

We can now state the main theorem of this section, which provides a condition for \cref{PDEeigsystemper3} to have a solution. Since the spatial dynamics formulation \cref{PDEeigsystemper3} is equivalent to the PDE eigenvalue problem, this allows us to find the PDE eigenvalues near the origin. This theorem is analogous to \cite[Theorem 2]{Sandstede1998}, with the $n\times n$ matrix in that theorem replaced by a $2n\times 2n$ block matrix. The proof is given in \cref{sec:blockmatrixproof}.

\begin{theorem}\label{blockmatrixtheorem}
Assume Hypotheses \ref{hyp:E}, \ref{hyp:H}, \ref{hyp:hypeq}, \ref{Qexistshyp}, and \ref{hyp:transverse}, and \ref{hyp:dccpos}. Let $Q_1(x)$ be the transversely constructed, symmetric primary pulse solution to \cref{genODE} from \cref{Qexistshyp}, and let $Q(x) = (Q_1(x), 0)$. Let $\Psi(x)$ and $V^c(x)$ be defined as in \cref{varadjsolutions}. Choose any periodic parameterization $(m_0, \dots, m_{n-1}, \theta)$ with $\theta \notin \{-\pi + p^*, p^* \}$. Let $r_*$ be defined as in \cref{th:perexist}, and for $r \leq r_*$, let $Q_n(x; r)$ be the corresponding periodic $n$-pulse solution. Then there exists $r_1 \leq r_*$ and $\delta > 0$ with the following property. For $r \leq r_1$, there exists a bounded, nonzero solution $V(x)$ of \cref{PDEeigsystemper3} for $|\lambda| < \delta$ and $|\Re \lambda| < r^{1/4}$ if and only if
\begin{equation}\label{blockmatrixcond}
\det E(\lambda) = 0,
\end{equation}
where $E(\lambda)$ is the $2n \times 2n$ block matrix
\begin{equation}\label{blockeq}
E(\lambda) = 
\begin{pmatrix}
K(\lambda) - \frac{1}{2} \lambda \tilde{M} K^+(\lambda) & \lambda^2 M_c I \\
-\frac{1}{2} \lambda M_c K^+(\lambda) & S - \lambda^2 MI 
\end{pmatrix} +
\begin{pmatrix}C_1 & D_1 \\ C_2 & D_2
\end{pmatrix}.
\end{equation}
The individual terms in $E(\lambda)$ are as follows:
\begin{enumerate}[(i)]
\item $K(\lambda)$ is the periodic, bi-diagonal matrix
\begin{align*}
K(\lambda) &=  
\begin{pmatrix}
e^{-\nu(\lambda)X_1} & -e^{\nu(\lambda)X_0} \\
-e^{\nu(\lambda)X_1} & e^{-\nu(\lambda)X_0} \\
\end{pmatrix} && n = 2\\
K(\lambda) &=  
\begin{pmatrix}
e^{-\nu(\lambda)X_1} & & & & & -e^{\nu(\lambda)X_0} \\
-e^{\nu(\lambda)X_1} & e^{-\nu(\lambda)X_2} \\
& -e^{\nu(\lambda)X_2} & e^{-\nu(\lambda)X_3} \\
  & & \ddots & && \\
& & & & -e^{\nu(\lambda)X_{n-1}} & e^{-\nu(\lambda)X_0}
\end{pmatrix} && n > 2,
\end{align*}
where $\nu(\lambda)$ is defined in \cref{nulambdalemma}. $K^+(\lambda)$ is the same matrix with all terms positive.

\item $S$ is the symmetric banded matrix
\begin{align}\label{Asymm}
S &= \begin{pmatrix}
-a_1 - a_0 & a_1 + a_0 \\
a_1 + a_0 & - a_1 - a_0 
\end{pmatrix} && n = 2 \\
S &= \begin{pmatrix}
-a_{n-1} - a_0 & a_0 & & &  & a_{n-1}\\
a_0 & -a_0 - a_1 &  a_1 \\
& a_1 & -a_1 - a_2 &  a_2 \\
& \ddots & \ddots & \ddots \\
a_{n-1} & & & & a_{n-2} & -a_{n-2} - a_{n-1} \\
\end{pmatrix} && n > 2,
\end{align}
where $a_i = \langle \Psi(X_i), Q'(-X_i) \rangle$.

\item $M$, $M_c$, and $\tilde{M}$ are the Melnikov-type integrals
\begin{align*}
M &= \int_{-\infty}^\infty q(y) \partial_c q(y) dy, \quad
M_c = \int_{-\infty}^\infty \partial_c q(y) dy, \quad
\tilde{M} = \int_{-\infty}^{\infty} \left(v^c(y) - \frac{1}{c}\right) dy,
\end{align*}
where $v^c(y)$ is the first component of $V^c(y)$, and $q(y)$ is the first component of $Q(y)$.

\item The remainder matrices $C_i$ and $D_i$ are analytic in $\lambda$ and have uniform bounds
\begin{align*}
|C_1| &\leq C (|\lambda| + r^{1/2})^2, \quad
|D_1| \leq C |\lambda|(|\lambda| + r^{1/2})^2 \\
|C_2| &\leq C (|\lambda| + r^{1/2})^2, \quad
|D_2| \leq C (|\lambda| + r^{1/2})^3.
\end{align*}
\end{enumerate}
\end{theorem}

\noi The condition $|\Re \lambda| < r^{1/4}$ is used to simplify the analysis. We will see when we apply the theorem in the following sections that this condition is satisfied for sufficiently small $r$.

\subsection{Spectrum of periodic single pulse}\label{sec:persingle}

The simplest case is the periodic single pulse. There is a only single length parameter $X_0$, which is the same as the domain length $X$, and the block matrix $E(\lambda)$ is a $2\times 2$ matrix. The form of $E(\lambda)$ is given in the following lemma. Proofs of all results in this section are given in \cref{sec:singlepulse}. 

\begin{lemma}\label{lemma:1blockmatrix}
For a periodic single pulse, the block matrix $E(\lambda)$ from \cref{blockmatrixtheorem} is the $2 \times 2$ matrix
\begin{align}\label{1pblockmatrix}
E(\lambda) &= 
\begin{pmatrix}
-2 \sinh(\nu(\lambda) X) - \tilde{M}\lambda \cosh(\nu(\lambda) X) & M_c \lambda^2 \\
-M_c \lambda \cosh(\nu(\lambda)X) & - M \lambda^2
\end{pmatrix} +
\begin{pmatrix}
c_1 & d_1 \\ c_2 & d_2
\end{pmatrix},
\end{align}
where the remainder terms are scalars with bounds
\begin{align*}
|c_1|, |c_2| &\leq C |\lambda|(|\lambda| + r^{1/2}), \qquad |d_1|, |d_2| \leq C |\lambda|^2(|\lambda| + r^{1/2}).
\end{align*}
In addition, $\det E(-\lambda) = -\det E(\lambda)$, and 
\begin{equation}\label{1pblockmatrixdet}
\begin{aligned}
\det E(\lambda) &= \lambda^2 \Big( 2 M \sinh(\nu(\lambda)X)(1 + \mathcal{O}(|\lambda| + r^{1/2} )) + \lambda(M \tilde{M} + M_c^2 )\cosh(\nu(\lambda)X) \\
&\qquad+ \mathcal{O}(|\lambda|(|\lambda| + r^{1/2} )) \Big).
\end{aligned}
\end{equation}
\end{lemma}

\noi Using this lemma, we can compute the nonzero essential spectrum eigenvalues close to the origin for the periodic single pulse. We emphasize that this does not locate all of the essential spectrum eigenvalues, but only those near the origin, i.e. those of sufficiently small magnitude.

\begin{theorem}\label{theorem:1pess}
Assume Hypotheses \ref{hyp:E}, \ref{hyp:H}, \ref{hyp:hypeq}, \ref{Qexistshyp}, and \ref{hyp:transverse}, and \ref{hyp:dccpos}. Let $r_1$ and $\delta$ be as in \cref{blockmatrixtheorem}. Then there exists $r_2 \leq r_1$ such that for any $r \in \mathcal{R}$ with $r \leq r_2$, the following holds regarding the nonzero essential spectrum eigenvalues. Let $N$ be any positive integer such that $N c \pi/X < \delta$. Then the first $2N$ nonzero essential spectrum eigenvalues are given by $\lambda = \{ \pm \lambda_m^{\text{ess}} : m = 1, \dots, N \}$, where
\begin{align}\label{1pess}
\lambda_m^{\text{ess}}(r) = c \frac{m \pi i}{X} \left( \frac{1}{1 + c \frac{M\tilde{M} + M_c^2 }{2 M X}}\right) +  \mathcal{O}\left( \frac{m^3}{|\log r|^3} \right)
\end{align}
is on the imaginary axis.
\end{theorem}

\subsection{Spectrum of periodic double pulse}\label{sec:perdouble}

For the next application, we consider the periodic double pulse. In this case, the block matrix $E(\lambda)$ is a $4\times 4$ matrix, the form of which is given in the following lemma. Proofs of all results in this section are given in \cref{sec:doublepulse}. 

\begin{lemma}\label{lemma:2blockmatrix}
For a periodic 2-pulse, the block matrix $E(\lambda)$ from \cref{blockmatrixtheorem} is the $4 \times 4$ matrix
\begin{equation}\label{dpSmatrix}
\begin{aligned}
E(&\lambda) = \\
&\begin{pmatrix}
e^{-\nu(\lambda)X_1} -\frac{1}{2}\lambda \tilde{M} e^{-\nu(\lambda)X_1} & -e^{\nu(\lambda)X_0} -\frac{1}{2}\lambda \tilde{M} e^{\nu(\lambda)X_0} & M_c \lambda^2 & 0 \\
-e^{\nu(\lambda)X_1} -\frac{1}{2}\lambda \tilde{M} e^{\nu(\lambda)X_1} & e^{-\nu(\lambda)X_0} -\frac{1}{2}\lambda \tilde{M} e^{-\nu(\lambda)X_0} & 0 & M_c \lambda^2 \\
-\frac{1}{2}\lambda M_c e^{-\nu(\lambda)X_1} & -\frac{1}{2}\lambda M_c e^{\nu(\lambda)X_0} &-a-\lambda^2 M & a \\
-\frac{1}{2}\lambda M_c e^{\nu(\lambda)X_1} & -\frac{1}{2}\lambda M_c e^{-\nu(\lambda)X_0}  & a & -a-\lambda^2 M
\end{pmatrix} + R(\lambda),
\end{aligned}
\end{equation}
where
\begin{equation}\label{2pa}
a = \langle \Psi(X_0), Q'(-X_0) \rangle + \langle \Psi(X_1), Q'(-X_1) \rangle.
\end{equation}
The remainder matrix is a $4 \times 4$ matrix of the form
\begin{align}
R(\lambda) = 
\begin{pmatrix} 
c_1(\lambda) & \tilde{c}_1(\lambda) & \lambda d_1(\lambda) & \lambda \tilde{d}_1(\lambda) \\ 
-c_1(-\lambda) & -\tilde{c}_1(-\lambda) & -\lambda \tilde{d}_1(-\lambda) & -\lambda d_1(-\lambda) \\ 
c_2(\lambda) & \tilde{c}_2(\lambda) & d_0 + \lambda d_2(\lambda) & -d_0 + \lambda \tilde{d}_2(\lambda) \\ 
-c_2(-\lambda) & -\tilde{c}_2(-\lambda) & -d_0 - \lambda \tilde{d}_2(-\lambda) & d_0 - \lambda d_2(-\lambda)
\end{pmatrix},
\end{align}
where the individual entries are scalars with bounds
\begin{align*}
|d_0| &\leq C r^{3/2} \\
|c_i|, |\tilde{c}_i|, |d_i|, |\tilde{d}_i| &\leq C (|\lambda| + r^{1/2})^2 && i = 1, 2.
\end{align*}
In addition, $\det E(-\lambda) = -\det E(\lambda)$, and 
\begin{equation}\label{2pblockmatrixdet}
\begin{aligned}
\det E(&\lambda) = -2 \lambda^2 (2a + \lambda^2 M + R_1) \left( M \sinh(\nu(\lambda)X) + \lambda (M \tilde{M} + M_c^2 ) \cosh(\nu(\lambda)X) \right) \\
&+4 a \lambda^3 M_c^2 \sinh(\nu(\lambda)X_1)\sinh(\nu(\lambda)X_0) 
+ R_2 \lambda^2 \sinh(\nu(\lambda)(X_1 - X_0)) \\
&+ \lambda^2 R_3 \sinh(\nu(\lambda)X) + \lambda^3 R_4,
\end{aligned}
\end{equation}
where the $R_i$ are scalars with bounds 
\begin{equation*}
|R_1| \leq C r^{3/2}, \quad |R_2|, |R_3|, |R_4| \leq C(|\lambda| + r^{1/2})^4.
\end{equation*}
\end{lemma}

We will first consider the case where the interaction eigenvalues are ``out of the way'' of the essential spectrum eigenvalues. Since the interaction eigenvalues scale as $r^{1/2}$ and the essential spectrum eigenvalues scale as $1/|\log r|$, we can always choose $r$ sufficiently small so that this is the case. Provided we do this, the interaction eigenvalue pattern for asymmetric periodic 2-pulses is determined by the parameter $m_0$ used in the construction of the solution.

\begin{theorem}\label{theorem:2peigsassym}
Assume Hypotheses \ref{hyp:E}, \ref{hyp:H}, \ref{hyp:hypeq}, \ref{Qexistshyp}, and \ref{hyp:transverse}, and \ref{hyp:dccpos}. Let $r_1$ and $\delta$ be as in \cref{blockmatrixtheorem}. Then for every $m_0 \in \{0, 1\}$ and $s_1 > p^*$ there exists $r_2 = r_2(m_0, s_1) \leq r_1$ such that for any $r \in \mathcal{R}$ with $r \leq r_2$, the following hold regarding the spectrum associated with the asymmetric periodic 2-pulse $Q_2(x; m_0, s_1, r)$.

\begin{enumerate}[(i)]
\item Let $N$ be any positive integer such that $N c \pi/X < \delta$. Then the first $2N$ nonzero essential spectrum eigenvalues are given by $\lambda = \{ \pm \lambda_m^{\text{ess}} : m = 1, \dots, N \}$, where
\begin{align}\label{2pess}
\lambda_m^{\text{ess}}(r) = c \frac{m \pi i}{X}\left( \frac{1}{1 + c \frac{M\tilde{M} + M_c^2}{MX}} \right) +  \mathcal{O}\left( \frac{m^3}{|\log r|^3} \right)
\end{align}
is on the imaginary axis.

\item There is a pair of interaction eigenvalues located at $\lambda = \pm \lambda^{\text{int}}(r)$, where
\begin{align*}
\lambda^{\text{int}}(r) = \sqrt{-\frac{2a}{M}} + \mathcal{O}\left( r \right),
\end{align*}
$a$ is defined in \cref{2pa}, and $|\lambda^{\text{int}}(r)| < \frac{1}{2}|\lambda_1^{\text{ess}}(r)|$. These are real when $m_0 = 0$ and purely imaginary when $m_0 = 1$.
\item There is an eigenvalue at 0 with algebraic multiplicity 3. 
\end{enumerate}
\end{theorem}

\noi We note that since the interaction eigenvalues scale as $r^{1/2}$, and the essential spectrum eigenvalues are purely imaginary, the condition $|\Re \lambda| < r^{1/4}$ is satisfied for sufficiently small $r$.

\begin{remark}The essential spectrum eigenvalues are not identical for the periodic single pulse and the periodic double pulse. In particular, note the additional factor of 2 in the denominator of the term in parentheses in \cref{1pess}. To leading order, however, the nonzero essential spectrum eigenvalues are located at $c \frac{m \pi i}{X}$ for nonzero integer $m$ in both cases.
\end{remark}

Next, we consider the symmetric periodic 2-pulse. As long as we are away from the pitchfork bifurcation points (i.e. as long as $a \neq 0$) the results of \cref{theorem:2peigsassym} hold; the only difference is that the eigenvalue pattern is determined by the sign of $a$ rather than by $m_0$ (see \cref{lemma:chara} below). Thus we only need to consider what happens at the pitchfork bifurcation point, which is given by the following theorem.

\begin{theorem}\label{theorem:2peigssym}
Assume Hypotheses \ref{hyp:E}, \ref{hyp:H}, \ref{hyp:hypeq}, \ref{Qexistshyp}, and \ref{hyp:transverse}, and \ref{hyp:dccpos}, and let $r_1$ be defined as in \cref{blockmatrixtheorem}. Then there exists $r_2 \leq r_1$ such that for all $r \in \mathcal{R}$ with $r \leq r_2$ and for $m_0 \in \{0, 1\}$, there is eigenvalue at 0 with algebraic multiplicity 5 for the symmetric periodic 2-pulse $\tilde{Q}_2(x; m_0, p^*(m_0; r), r)$, where $p^*(m_0; r)$ is the pitchfork bifurcation point defined in \cref{2pulsebifurcation}.
\end{theorem}

Finally, we consider what happens when an essential spectrum eigenvalue collides with an interaction eigenvalue on the imaginary axis. For simplicity, we will only prove the result for the first collision. The existence of the first Krein bubble is given in the following theorem, which also provides numerically verifiable estimates for its size and location. 

\begin{theorem}\label{th:Kreinbubble}
Assume Hypotheses \ref{hyp:E}, \ref{hyp:H}, \ref{hyp:hypeq}, \ref{Qexistshyp}, and \ref{hyp:transverse}, and \ref{hyp:dccpos}. Choose $m_0 = 1$, and let $r_1$ be as in \cref{blockmatrixtheorem}. Let 
\begin{align}\label{deflambdastar}
\lambda_*(r) &= \sqrt{ \frac{-2a(r) - R_1 }{M} } = \sqrt{ \frac{-2a(r)}{M} } + \mathcal{O}(r),
\end{align}
where $R_1$ is defined in \cref{lemma:2blockmatrix}, and let $Q_2^*(x; r)$ be the periodic 2-pulse solution from \cref{2pulsebifurcation} with domain size $X = X_*(r)$, where
\begin{align}\label{defXstar}
X_*(r) &= c \left( \frac{\pi i}{\lambda_*(r)} - \frac{M \tilde{M} + M_c^2 }{M}\right).
\end{align}
Define $T_1(r) > 0$ by 
\begin{align}\label{defT0}
T_1(r) = \frac{M_c^2 }{2 \pi M^2 c } |\lambda_*|^5 X_0^2.
\end{align}
For $s \in \left[-2\sqrt{T_1(r)}, 2 \sqrt{T_1(r)}\right]$, let $Q_2(x; s, r)$ be the periodic 2-pulse solution with domain size $X = X(s, r)$, where
\begin{equation}
X(s,r) = X_*(r) + \frac{2 c \pi s}{|\lambda_*(r)|^2} + \mathcal{O}(s^2).
\end{equation}
Then there exists $r_2 \leq r_1$ such that for all $r \in \mathcal{R}$ with $r \leq r_2$, the following holds for the linearization of the PDE about $Q_2(x; s, r)$.
\begin{enumerate}[(i)]
	\item There is a pair of eigenvalues located at
	\begin{align}\label{Kreineigs}
	\lambda = \lambda_*(r) - s i \pm \sqrt{ T_1(r) -  s^2} + \mathcal{O}\left(r^{5/4}|\log r|^{1/2} \right).
	\end{align}

	\item For
	\[
	s = s_\pm(r) = \pm \sqrt{T_1(r)}\left( 1 +  \mathcal{O}\left(\frac{1}{|\log r|} \right) \right),
	\]
	there is a double eigenvalue on the imaginary axis at 
    \[
    \lambda = \lambda_*(r) + s_\pm(r) i + \mathcal{O}\left(\frac{1}{|\log r|} \right),
    \]
    which occurs when 
	\begin{equation}\label{KreinDeltaX}
	X(s, r) = X_*(r) \pm \Delta X(r) + \mathcal{O}(s^2), \qquad \Delta X(r) = \frac{2 c \pi}{|\lambda_*(r)|^2}\sqrt{T_1(r)}.
	\end{equation}
\item For $s \in (s_-, s+)$, equation \cref{Kreineigs} describes, to leading order, a circle of radius $\sqrt{T_1(r)}$ in the complex plane, which is the Krein bubble. The pair of eigenvalues is symmetric across the imaginary axis.
\item For $s \in \Big[-2\sqrt{T_1(r)}, s_-\Big) \cup \Big(s_+, 2 \sqrt{T_1(r)}\Big]$, the eigenvalues \cref{Kreineigs} are on the imaginary axis.
\end{enumerate}
\end{theorem}

\noi We note that maximum real part of the Krein bubble is order $r^{5/4}|\log r|$, thus the condition $|\Re \lambda| < r^{1/4}$ is satisfied for sufficiently small $r$.

\begin{remark}\label{remark:kreinbubbles}
It is straightforward to adapt \cref{th:Kreinbubble} to locate subsequent Krein bubbles. For any positive integer $N$, there exists $r_2 = r_2(N)$ with $r_2 \leq r_1$ such that for $r \leq r_2$ and $m = 1, \dots, N$, a Krein bubble occurs when the $m$-th essential spectrum eigenvalue collides with the interaction eigenvalue on the imaginary axis. The radius of $m$-th Krein bubble in the complex plane is approximately $\sqrt{T_1(r) / m}$, and the Krein collisions occur at approximately $X = X^m_*(r) \pm \Delta X^m(r)$, where 
\begin{equation}\label{eq:Kreinradiusm}
X^m_*(r) = c \left( \frac{m \pi i}{\lambda_*(r)} - \frac{M \tilde{M} + M_c^2 }{M}\right), \qquad \Delta X^m(r) = \sqrt{m} \Delta X(r).
\end{equation}
Note that this requires $N$ to be chosen first, and $r_2$ depends on $N$. See \cref{sec:conclusions} for a discussion on what occurs with subsequent Krein bubbles when $r$ is fixed.
\end{remark}

\section{Numerical Results}\label{sec:numerics}

In this section, we present numerical results for the existence and spectrum of periodic multi-pulse solutions to KdV5. We start with the construction of the primary pulse solution. For $p = -1$ and $c = 36/169$, the exact solution to \cref{KdV5eq4} is known \cite{Pelinovsky2007}
\begin{equation}\label{KdV5exactsol}
q(x) = \frac{105}{338}\sech^4\left(\frac{x}{2\sqrt{13}} \right).
\end{equation}
We use AUTO \cite{AUTO} for parameter continuation in $c$ and $p$ until $-2 \sqrt{c} < p < 2 \sqrt{c}$, so that \cref{hyp:hypeq} is satisfied. Following the AUTO demo \texttt{kdv}, we formulate the problem using equation \cref{KdV5ham2}, and use a small parameter $\epsilon$ to break the Hamiltonian structure. We impose periodic boundary conditions and rescale the domain from $[-X, X]$ to $[0, 1]$, using the domain size $X$ as a parameter. 

To construct a periodic double pulse $q_2(x)$, we discretize equation \cref{KdV5eq4} using Fourier spectral differentiation matrices to enforce periodic boundary conditions. As an initial ansatz, we take two copies of the primary pulse joined together at the distances predicted by \cref{2pulsebifurcation}. We then solve for the periodic double pulse using Matlab's \texttt{fsolve} function. This same procedure can also be used to construct arbitrary periodic multi-pulses. We can also vary the domain size $X$ by parameter continuation in AUTO. Using this, we verify that asymmetric periodic 2-pulses bifurcate from symmetric periodic 2-pulses in a series of pitchfork bifurcations (\cref{fig:KdV5periodicMP}, center panel).  

\begin{figure}
\begin{center}
\begin{tabular}{cc}
\includegraphics[width=8cm]{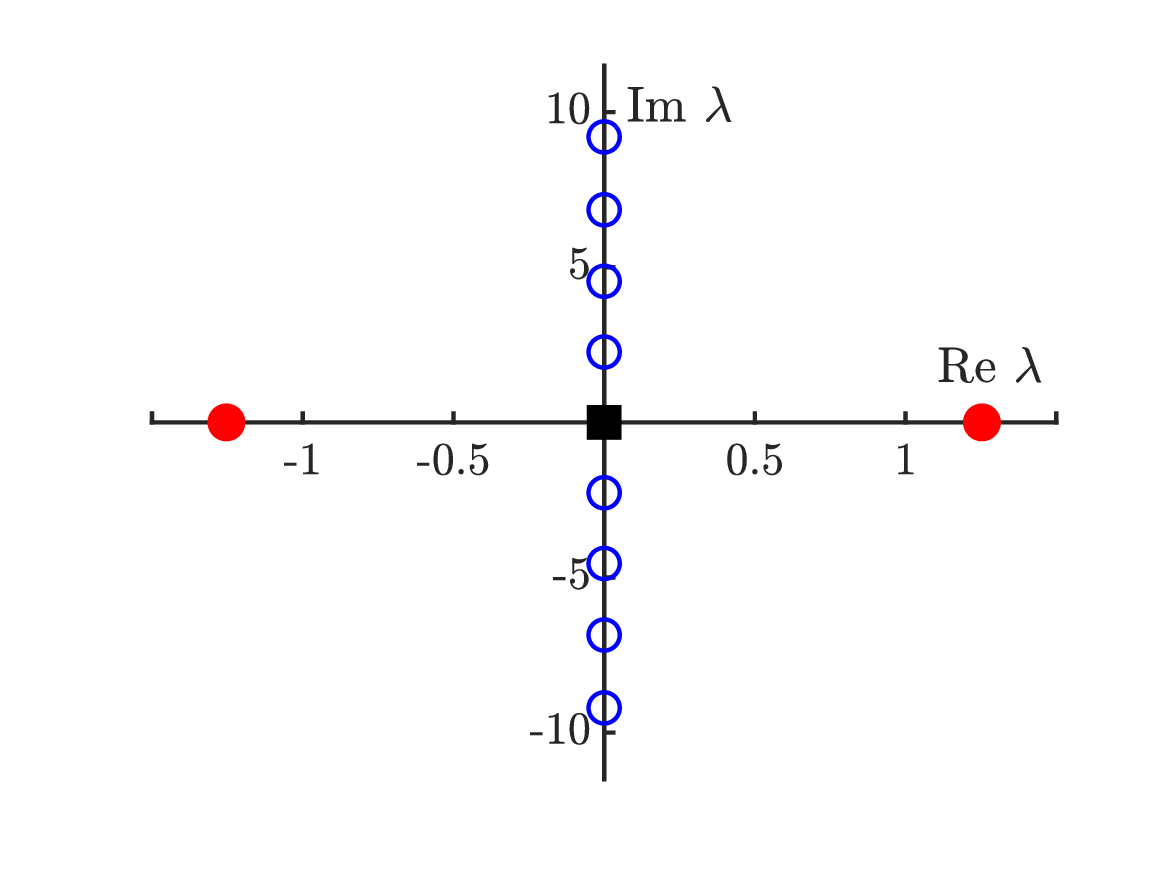} &
\includegraphics[width=8cm]{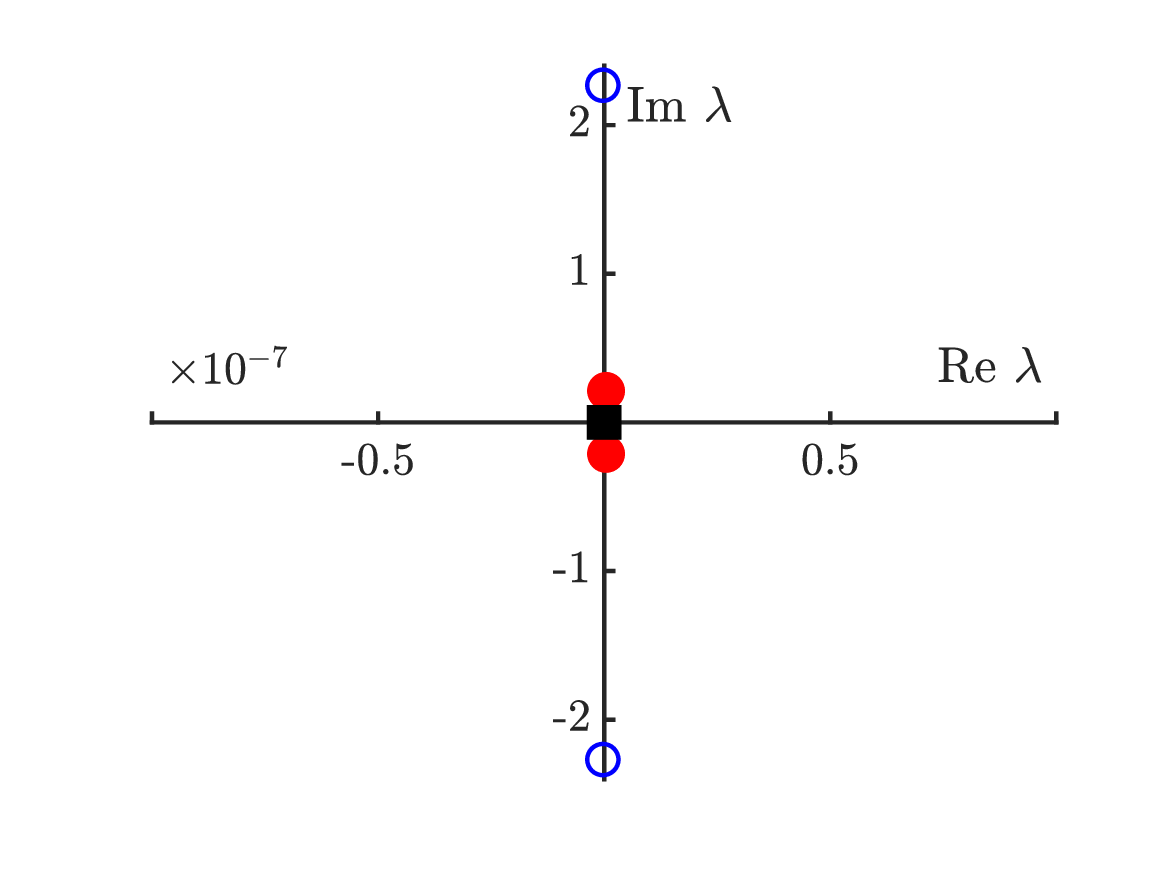}
\end{tabular}
\end{center}
\caption[Spectrum of periodic 2-pulse solutions]{Spectrum of periodic 2-pulse solutions for KdV5, showing the interaction eigenvalues (red dots), essential spectrum eigenvalues (blue open circles), and double eigenvalue at origin (black square). Left panel ($m_0 = 0$) has real interaction eigenvalues, right panel ($m_0 = 1$) has imaginary eigenvalues. Fourier spectral methods with $N = 1024$ grid points, $p = -1$, $c = 20$, $X = 30$.}
\label{fig:KdV5eigs1}
\end{figure}

\begin{figure}
\begin{center}
\includegraphics[width=8cm]{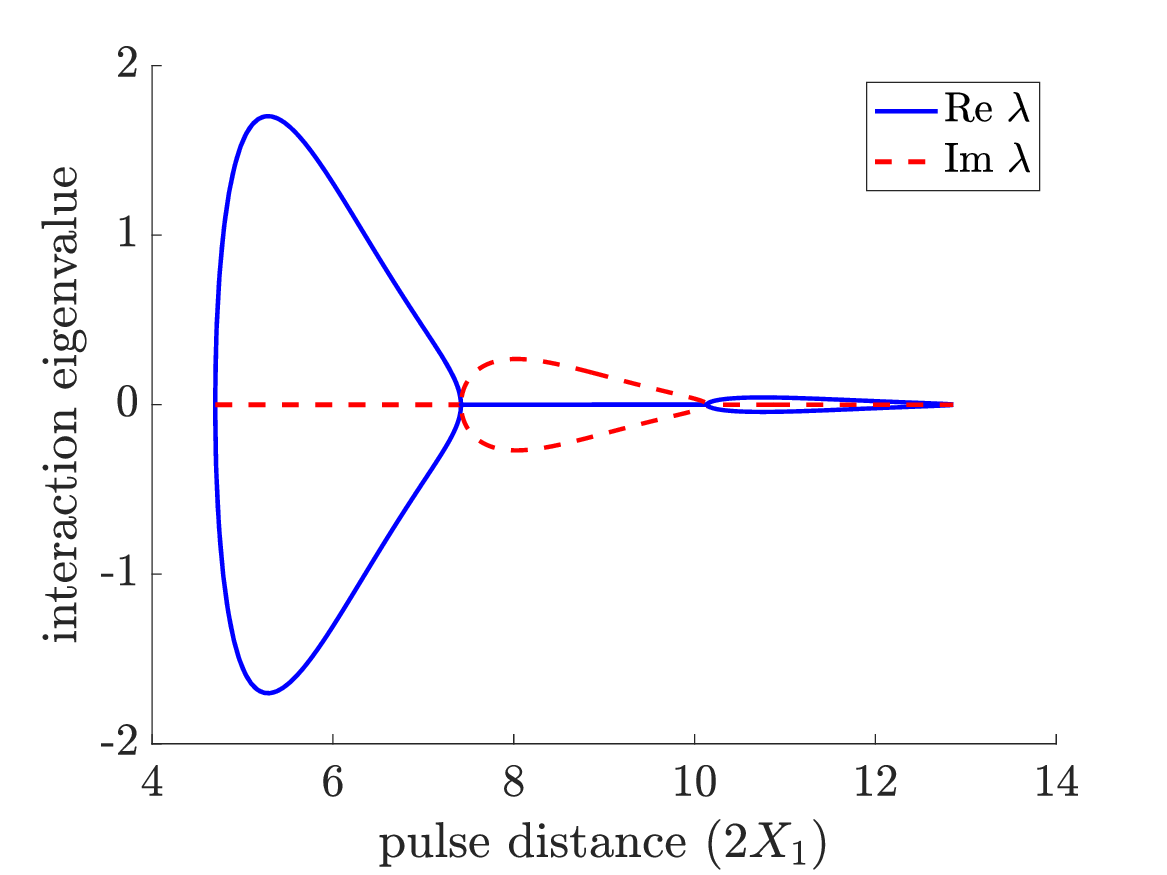}
\includegraphics[width=8cm]{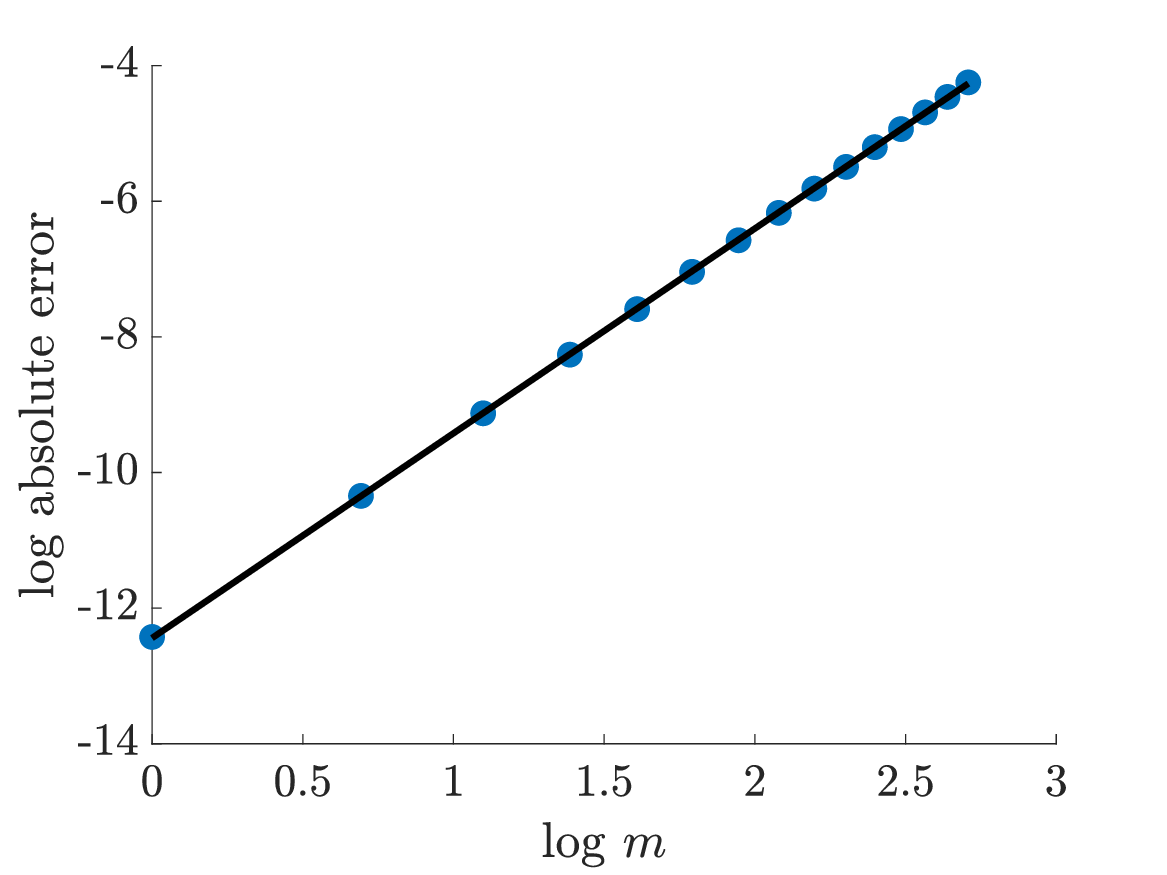}
\end{center}
\caption{Left panel: real part (blue solid line) and imaginary part (red dashed line) of interaction eigenvalues vs. pulse distance for symmetric periodic 2-pulses $(X_0 = X_1)$. Parameters are $p = -1$ and $c = 20$. Right panel: log of absolute error vs. $\log m$ for the first 15 essential spectrum eigenvalues with least squares linear regression line for periodic single pulse. Slope of regression line is 3.05. Parameters $p = -1$, $c = 10$, and $X = 200$.}
\label{fig:periodiceigbif}
\end{figure}

Next, we compute the spectrum of $\partial_x \calE''(q_2)$ by discretizing the linear operator using Fourier spectral differentiation matrices and using Matlab's \texttt{eig} function (\cref{fig:KdV5eigs1}). For asymmetric periodic 2-pulses $(X_0 \neq X_1)$, the interaction eigenvalue pattern depends only on the integer $m_0$ from the periodic parameterization. For $m_0 = 0$, the interaction eigenvalues are real, and for $m_0 = 1$, the interaction eigenvalues are purely imaginary (the real part of the eigenvalues computed with \texttt{eig} is less than $10^{-9}$). We can also compute the interaction eigenvalues for symmetric periodic 2-pulses (\cref{fig:periodiceigbif}, left panel). At each pitchfork bifurcation point, an eigenvalue bifurcation occurs, where a pair of interaction eigenvalues collides at 0 and switches from real to purely imaginary (or vice versa). The full interaction eigenvalue pattern (up to the first Krein bubble) is shown in the center panel of \cref{fig:KdV5periodicMP}.

We then compute the essential spectrum eigenvalues for periodic single pulses using \texttt{eig}, and compare the results to the leading order formula from \cref{theorem:1pess}. Plotting the log of absolute value of the error versus $\log m$ and constructing a least squares linear regression line (\cref{fig:periodiceigbif}, right panel), the absolute error is proportional to $m^3$, with a relative error in the exponent of less than 0.02, as predicted by \cref{theorem:1pess}. The results are similar for periodic double pulses using the leading order formula from \cref{theorem:2peigsassym}. 

\begin{figure}
\begin{center}
\includegraphics[width=9cm]{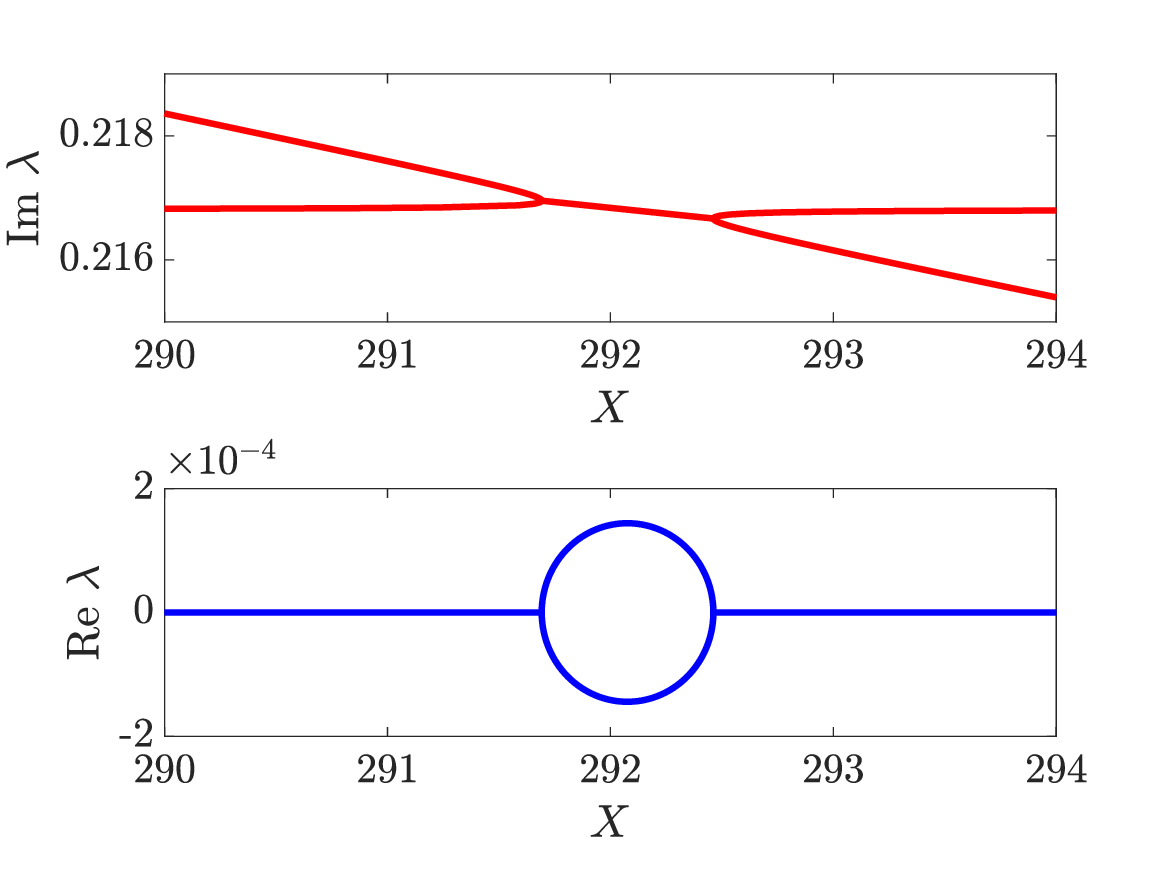}
\hspace{-0.75cm}
\includegraphics[width=8cm]{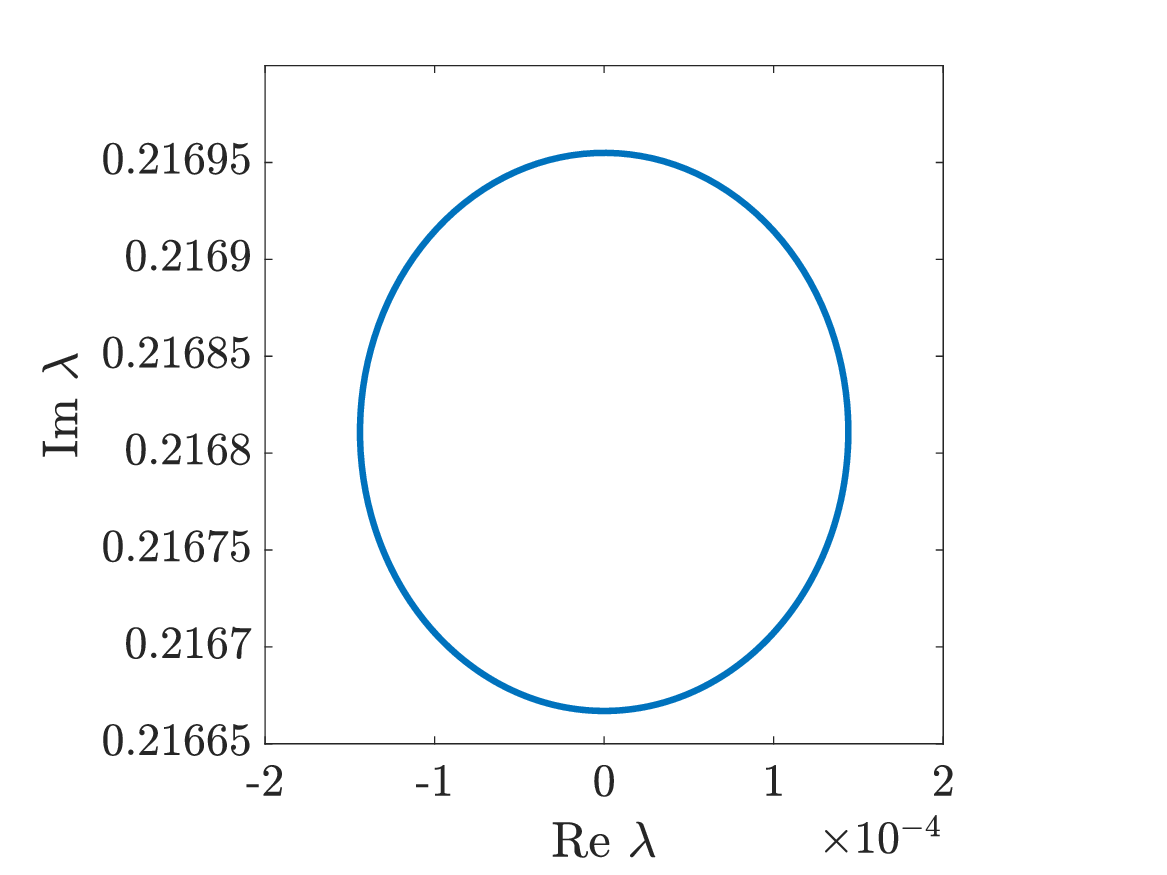}
\end{center}
\caption{Left panel shows collision of first essential spectrum eigenvalue with purely imaginary interaction eigenvalue as $X$ is increased. Imaginary part of eigenvalues in top panel (red), real part of eigenvalues in bottom panel (blue). Right panel plots the imaginary part vs. the real part of the eigenvalues in the Krein bubble as $X$ varies. Parameter continuation with AUTO in periodic domain length $X$, $p = -1$, $c = 20$.}
\label{fig:kreinbubble1}
\end{figure}

\begin{figure}
\begin{center}
\includegraphics[width=8cm]{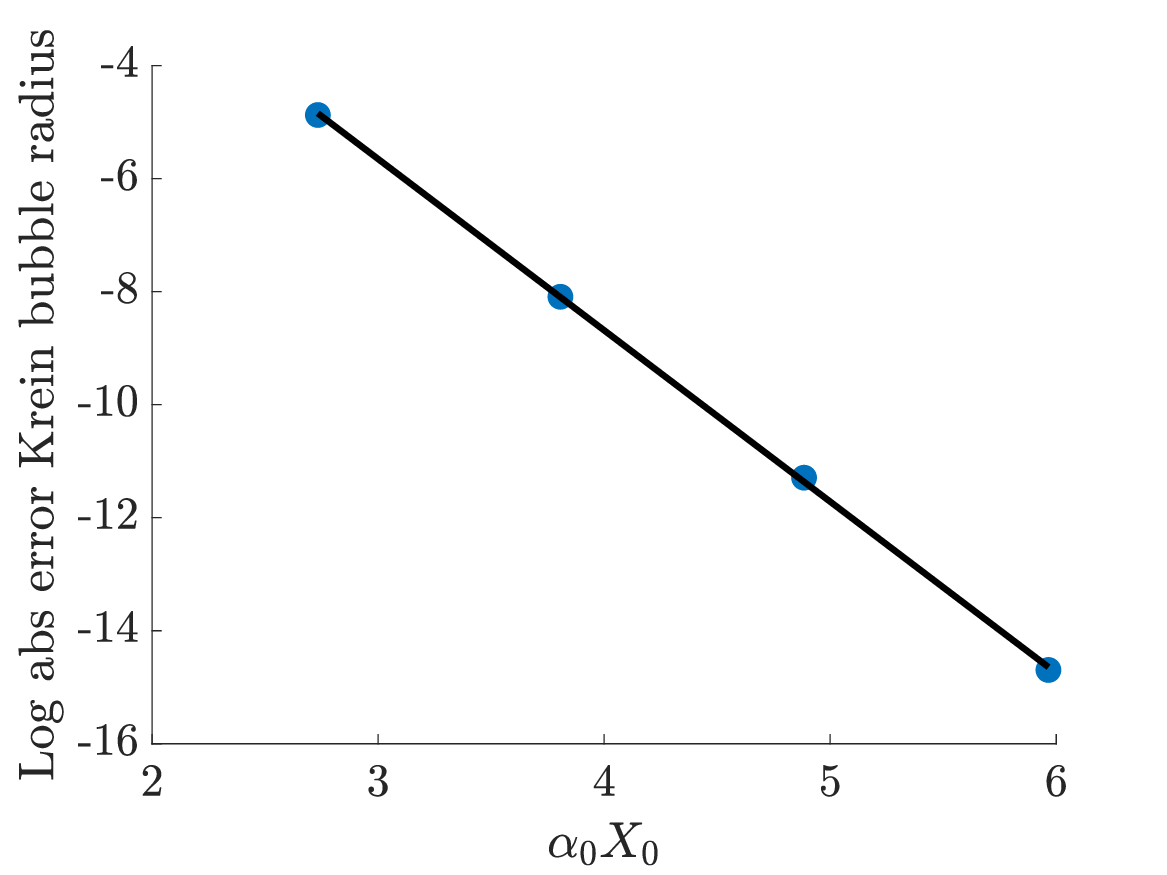}
\includegraphics[width=8cm]{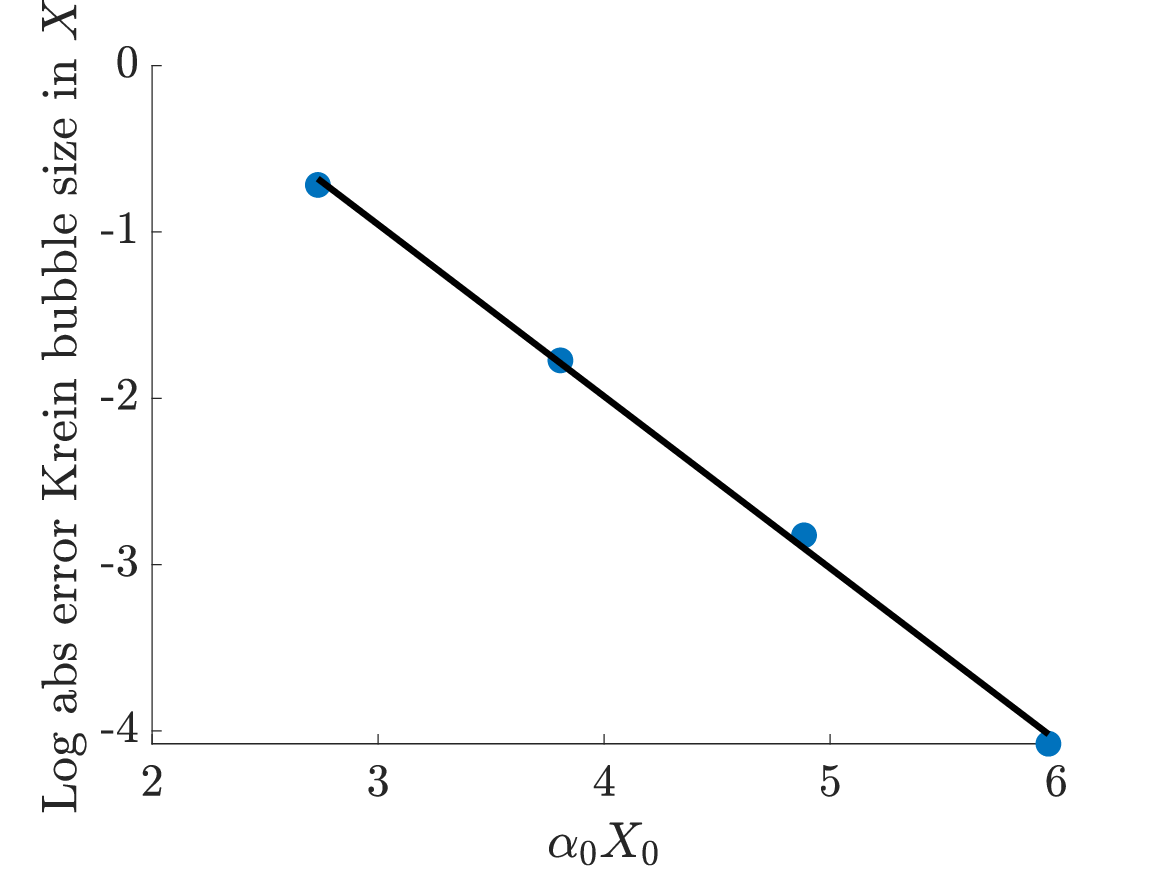}
\end{center}
\caption{Plot of log of absolute error for Krein bubble radius in the complex plane (left panel) and Krein bubble size in $X$ (right panel) vs. $\alpha_0 X_0$, with least squares linear regression lines for first four periodic double pulses with $m_0 = 1$. Slopes of regression lines -3.0319 (left) and -1.0332 (right). Parameters $a = 0.25$ and $b = 3$ in \cref{KdVpabc}, which corresponds to $p = 5.5$and $c = 10.5625$.
}
\label{fig:kreinerrors}
\end{figure}

Next, we look at what happens when we increase the periodic domain parameter $X$ using parameter continuation with AUTO. As predicted by \cref{th:Kreinbubble}, there is a brief instability bubble when the first essential spectrum eigenvalue collides with the imaginary interaction eigenvalue (\cref{fig:kreinbubble1}).
If $p$ and $c$ in \cref{KdV5eq4} are related by
\begin{align}\label{KdVpabc}
p &= -2(a-b), \quad c = (a+b)^2 && a, b > 0,
\end{align}
then the eigenvalues of $DF(0)$ are the quartet $\pm \sqrt{a} \pm \sqrt{b} i$. Choosing $a = 0.25$ and $b = 3$, so that the tail oscillations of the primary pulse are sufficiently rapid but do not decay too quickly, we can construct the first four periodic double pulses with $m_0 = 1$, together with the eigenfunctions corresponding to the imaginary interaction eigenvalue, to a sufficient degree of accuracy so that AUTO converges for both the existence and eigenvalue problems. 
\cref{fig:kreinerrors} plots the log of the absolute error of the Krein bubble radius in the complex plane ($\sqrt{T_1}$ from \cref{defT0}) and the log of the absolute error of the Krein bubble size in $X$ ($\Delta X$ from \cref{KreinDeltaX}) versus $\alpha_0 X_0$. The slopes of the least squares linear regression lines suggest that the Krein bubble radius in the complex plane is given by $\sqrt{T_1} + \mathcal{O}(e^{-3 \alpha_0 X_0})$ and that the Krein bubble radius in $X$ is given by $\Delta X + \mathcal{O}(e^{-\alpha_0 X_0})$, with relative errors in the exponent less than 0.01 and 0.03 (respectively). The leading order terms agree with \cref{th:Kreinbubble}, while the error term is higher order than predicted. The results for subsequent Krein bubbles discussed in \cref{remark:kreinbubbles} can similarly be verified numerically.

Finally, we present results of numerical timestepping experiments to illustrate the effects of the Krein bubble on the PDE dynamics of perturbations of periodic double pulses. Let $u_1(x)$ be the periodic single pulse solution to \cref{KdV5eq4}. The initial condition for the timestepping is the sum of two well-separated copies of the periodic single pulse,
\begin{equation}\label{eq:timestepIC}
u(x) = u_1(x-L) + u_1(x+L).
\end{equation}
The two pulses are separated by a distance $2L$, which is chosen to be close to the pulse separation distance for a periodic double pulse. Timestepping was performed using a pseudo-spectral method for spatial discretization and a fourth-order Runge-Kutta method for time evolution, as in \cite{Pelinovsky2007}. High frequency oscillations resulting from large essential spectrum eigenmodes were damped using a lowpass filter.

\begin{figure}
\begin{center}
\begin{tabular}{cc}
\includegraphics[width=8cm]{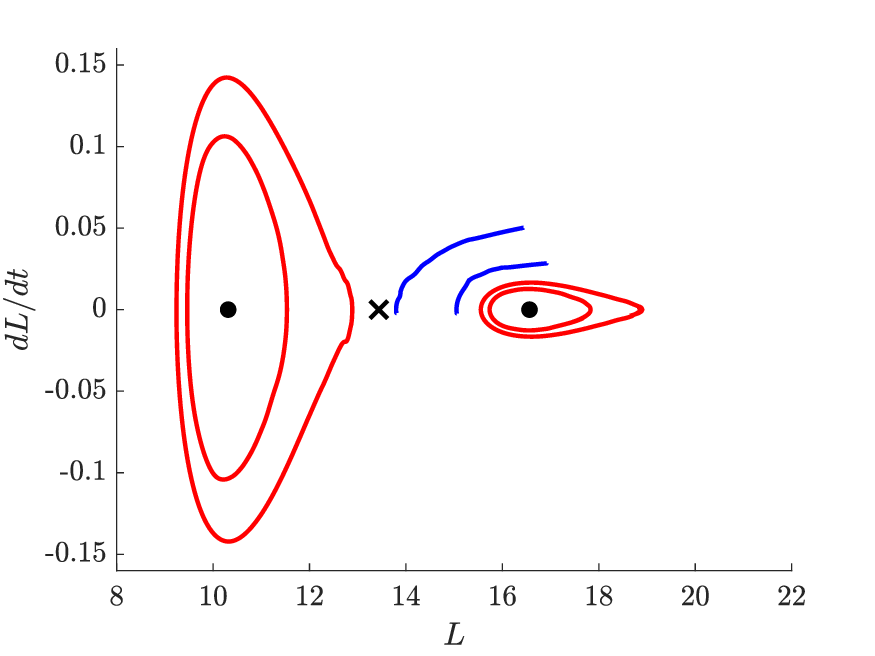} &
\includegraphics[width=8cm]{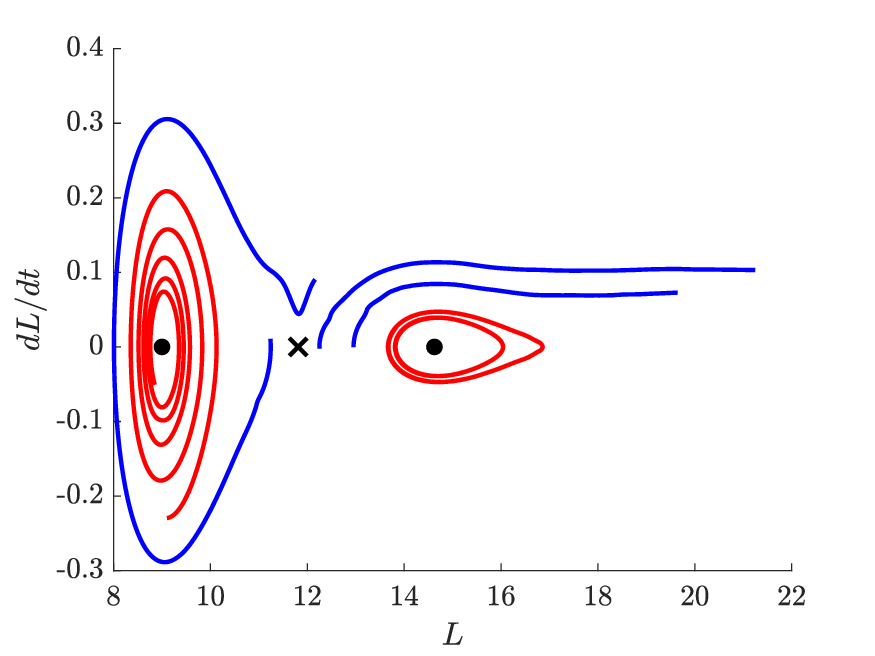} \\
\includegraphics[width=8cm]{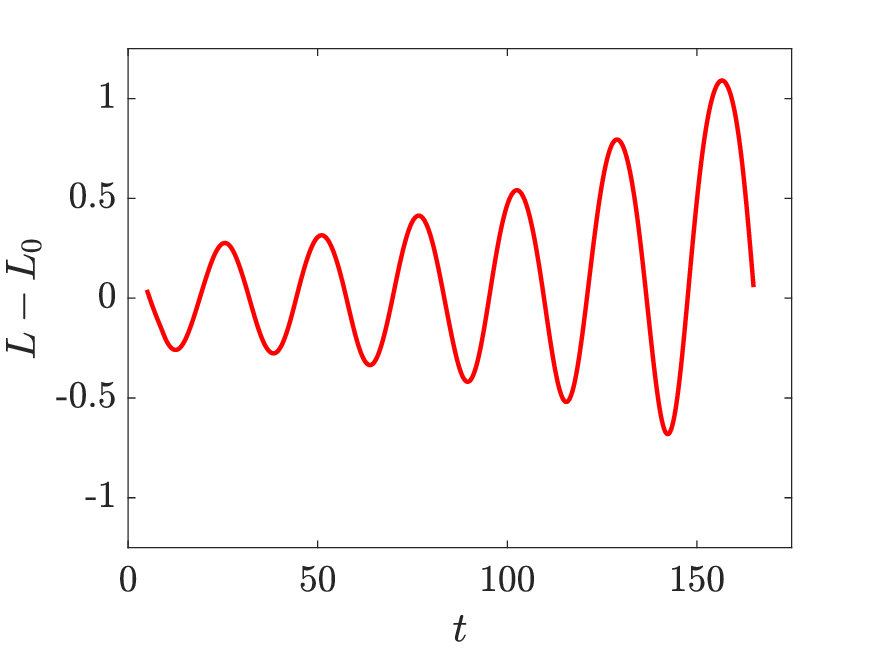} &
\includegraphics[width=8cm]{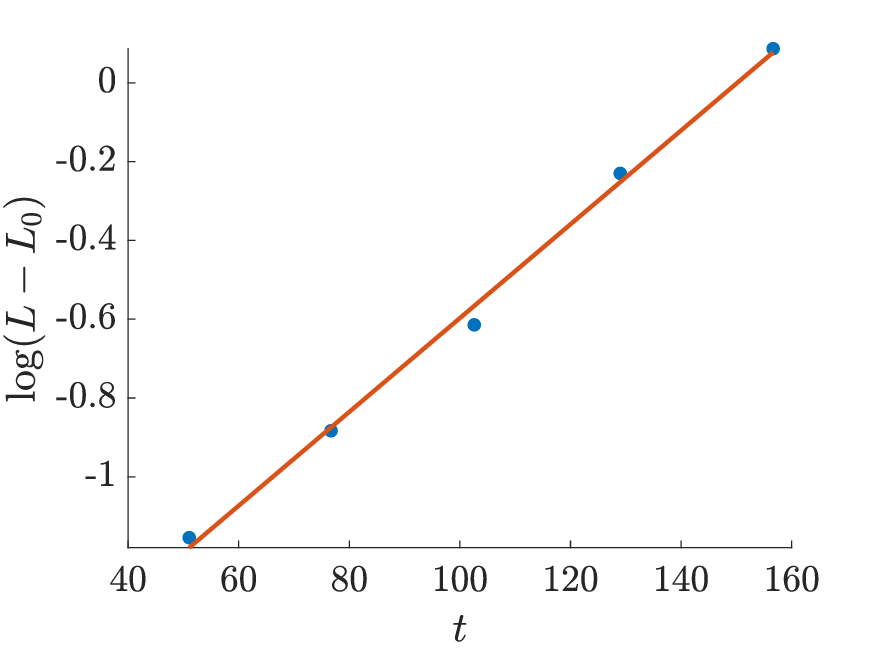}
\end{tabular}
\end{center}
\caption{Phase plane showing $dL/dt$ vs. $L$ for timestepping simulations, where the two peaks are separated by a distance $2L$. Top left (outside Krein bubble): parameters $p=1.1$ and $c=2.1025$, domain size $X=40$, and interaction eigenvalues $0.1120i$, $0.0384$, and $0.0134i$ (from left to right). Top right (inside Krein bubble), parameters $p=1.7$ and $c=2.7225$, domain size $X=36$, and interaction eigenvalues $0.0124 + 0.2498i$, $0.0968$, and $0.0398i$ (from left to right). Bottom left: $L-L_0$ vs. $t$ for solution starting near leftmost equilibrium point in top right phase portrait, where the two peaks in the unperturbed periodic double pulse are separated by a distance $2 L_0$. Bottom right: $\log(L-L_0)$ vs. $t$ at maxima of this solution.}
\label{fig:timestepbubble}
\end{figure}

The top panel of \cref{fig:timestepbubble} plots phase portraits $(L, dL/dt)$ for two different parameter configurations (see \cite[Figure 10]{Pelinovsky2007} for similar phase portraits). Neutrally stable periodic double pulses ($m_0$ odd) are marked with a black dot, and unstable periodic double pulses ($m_0$ even) are marked with black X. The top left panel of \cref{fig:timestepbubble} is the phase portrait corresponding to a parameter configuration outside the Krein bubble. The interaction eigenvalues for both periodic doubles pulses with $m_0$ odd are purely imaginary. These correspond to neutrally stable centers in the phase portrait, and the frequency of oscillation about these equilibria is within 5\% of the imaginary part of the corresponding interaction eigenvalue. The interaction eigenvalue for the periodic double pulse with $m_0$ even is real, which corresponds to an unstable saddle equilibrium in the phase portrait. The top right panel of \cref{fig:timestepbubble} is the phase portrait corresponding to a parameter configuration inside the Krein bubble. Since the interaction eigenvalue for the first periodic double pulse (leftmost equilibrium point) has a small, positive real part, trajectories starting near this unstable equilibrium slowly spiral outward (\cref{fig:timestepbubble}, bottom left). The average frequency of these oscillations is 0.2380, which is within 5\% of the imaginary part of the Krein bubble eigenvalue, and the exponential growth rate of the maxima of this solution (\cref{fig:timestepbubble}, bottom right) is 0.0119, which is within 5\% of the real part of the Krein bubble eigenvalue.

\section{Proof of Theorem \ref{transverseint}}\label{sec:transverseintproof}

The proof is similar to that of \cite[Lemma 6.2 and Lemma 6.4]{Kapitula2020}. By \cref{Qexistshyp}, $Q(0; c_0) \neq 0$, $H(Q(0; c_0); c_0) = 0$, and $\nabla_U H(Q(0; c_0); c_0) \neq 0$. By the implicit function theorem, there exists $\delta_0 > 0$ such that for $c \in (c_0 - \delta_0, c_0 + \delta_0)$, the 0-level set $H^{-1}(0; c)$ contains a smooth $(2m-1)$-dimensional manifold $K(c)$, with $K(c_0)$ containing $Q(0; c_0)$. The existence result and the smoothness of the map $c \mapsto Q(x; c)$ for $c \in (c_0 - \delta_0, c_0 + \delta_0)$ follow from the transverse intersection of $\tilde{W}^s(0; c_0)$ and $\tilde{W}^u(0; c_0)$ in $K(c_0) \subset H^{-1}(0; c_0)$, the implicit function theorem, and the smoothness of $F$. Symmetry with respect to the reversor $R$ follows from the symmetry of $Q(0; c_0)$ and $H$.

Fix $c \in (c_0 - \delta_0, c_0 + \delta_0)$. Since $Q(x; c)$ solves \cref{genODE}, $Q(x; c) \in C^1(\R, \R^{2m})$. By the stable manifold theorem, $Q(x; c)$ is exponentially localized, i.e. for every $\epsilon > 0$ there exists a constant $C$ such that for all $x \in \R$, $|Q(x; c)| \leq C e^{-(\alpha_0 - \epsilon)|x|}$. Substituting $Q(x; c)$ into \cref{genODE} and differentiating with respect to $c$, $\partial_c Q(x; c)$ satisfies
\begin{equation}\label{Qcprime}
[\partial_c Q(x; c)]' = D_U F(Q(x;c); c) \partial_c Q(x; c) - B Q(x;c),
\end{equation}
where $B$ is defined in \cref{defAB}. Define the linear operator $\calL: C^1(\R, \R^{2m}) \rightarrow C^0(\R, \R^{2m})$ by
\begin{equation}\label{defLtv}
\calL = \frac{d}{dx} - D_U F(Q(x;c); c).
\end{equation}
By equation \cref{Qcprime}, $B Q(x;c) \in \ran \calL$ and is exponentially localized. Since $DF(0; c)$ is hyperbolic, it follows from \cite[Lemma~4.2]{Palmer1984} and the roughness theorem for exponential dichotomies \cite{Coppel1978} that $\calL$ is Fredholm with index 0. By \cref{Qexistshyp}, $\ker \calL = \R Q'(x; c)$, thus the set of all bounded solutions to \cref{Qcprime} is given by $\{\partial_c Q(x; c) + \R Q'(x; c)\}$.

To show that $\partial_c Q(x; c)$ is exponentially localized, we reformulate equation \cref{Qcprime} in an exponentially weighted space. Choose $\epsilon \in (0,\alpha_0)$ and let $\eta(x)$ be a standard mollifier function \cite[Section~C.5]{Evans2010}. Let
\begin{equation}\label{defQZ}
Q(x; c) = Z(x; c) e^{-(\alpha_0 - \epsilon)r(x)},
\end{equation}
where $r(x) = \eta(x) * |x|$ is smooth, and for $|x| > 1$, $r(x) = |x|$ and $r'(x) = \pm 1$. Substituting \cref{defQZ} into \cref{Qcprime} and simplifying, we obtain the weighted equation
\begin{equation}\label{Zcprime}
[\partial_c Z(x; c)]' = [D_U F(Q(x;c); c) + (\alpha_0 - \epsilon) r'(x) ] Z(x; c) - e^{(\alpha_0 - \epsilon)r(x)} B Q(x;c),
\end{equation}
where the last term on the RHS is bounded. Define the weighted linear operator $\calL_{\alpha_0 - \epsilon}: C^1(\R, \R^{2m}) \rightarrow C^0(\R, \R^{2m})$ by
\begin{equation}\label{defLtvalpha}
\calL_{\alpha_0 - \epsilon} = \frac{d}{dx} - D_U F(Q(x;c); c) - (\alpha_0 - \epsilon) r'(x) \calI.
\end{equation}
Since $B Q(x;c) \in \ran \calL$, $e^{(\alpha_0 - \epsilon)r(x)} B Q(x;c) \in \ran \calL_{\alpha_0 - \epsilon}$. Since $D_U F(0; c)-(\alpha_0 - \epsilon)\calI$ is still hyperbolic with the same unstable dimension as $D_U F(0; c)$, it follows again from \cite[Lemma 4.2]{Palmer1984} that $\calL_{\alpha_0 - \epsilon}$ is Fredholm with index 0. Since $Q'(x; c)$ is exponentially localized by the stable manifold theorem, $e^{(\alpha_0 - \epsilon)r(x)} Q'(x; c)$ is bounded, thus since $Q'(x; c) \in \ker \calL$, $e^{(\alpha_0 - \epsilon)r(x)} Q'(x; c) \in \ker \calL_{\alpha_0 - \epsilon}$. Since any element in $\ker \calL_{\alpha_0 - \epsilon}$ gives an element of $\ker \calL$ via \cref{defQZ}, $\ker \calL_{\alpha_0 - \epsilon} = \R e^{(\alpha_0 - \epsilon)r(x)} Q'(x; c)$. Since $e^{(\alpha_0 - \epsilon)r(x)} B Q(x;c) \in \ran \calL_{\alpha_0 - \epsilon}$, the set of all bounded solutions to \cref{Zcprime} is given by $\{ \partial_c Z(x; c) + \R e^{(\alpha_0 - \epsilon)r(x)} Q'(x;c) \}$, which implies that $\partial_c Q(x; c) = \partial_c Z (x; c) e^{-(\alpha_0 - \epsilon)r(x)}$ is exponentially localized.

\section{Proof of existence results}\label{sec:existproof}

We will construct a periodic $n$-pulse $U(x)$ using Lin's method. For convenience of notation, we will denote the primary pulse by $Q(x)$ instead of $Q_1(x)$. Rather than taking $U(x)$ to be a piecewise perturbation of $Q(x)$, we adapt the technique in \cite{Sandstede1997} and instead take a piecewise ansatz of the form
\[
U_i^\pm(x) = Q^\pm(x; \beta_i^\pm) + \tilde{Q}_i^\pm(x),
\]
where $\beta_i^\pm$ parameterize the stable and unstable manifolds $\tilde{W}^s(0)$ and $\tilde{W}^u(0)$ near $Q(0)$. The functions $Q^\pm(x; \beta_i^\pm)$ lie in these manifolds, and the $\tilde{Q}_i^\pm$ are small remainder functions. In essence, we use the parameters $\beta_i^\pm$ to break the homoclinic orbit $Q(x)$, and the remainder functions $\tilde{Q}_i^\pm$ to glue the pieces back together. We will show that we can find a unique piecewise solution which generically has $n$ jumps in a specified direction. A periodic multi-pulse solution exists if and only if these $n$ jumps are all 0.

\subsection{Setup}

Using \cref{nondegencond}, we decompose the tangent spaces of the stable and unstable manifolds at $Q(0)$ as 
\begin{equation*}
\begin{aligned}
T_{Q(0)}\tilde{W}^u(0) &= \R Q'(0) \oplus Y^- \\
T_{Q(0)}\tilde{W}^s(0) &= \R Q'(0) \oplus Y^+.
\end{aligned}
\end{equation*}
It follows from \cref{nondegencond} that $Q'(x)$ is the unique bounded solution to the variational equation
\begin{equation}\label{vareq1}
V' = DF(Q(x)) V,
\end{equation}
and that there exists a unique bounded solution $\Psi(x)$ to the adjoint variational equation
\begin{equation}\label{adjvareq1}
W' = -DF(Q(x))^* W.
\end{equation}
(In both cases, uniqueness is up to scalar multiple.) Since we have a conserved quantity $H$, the following lemma gives the exact form of $\Psi(x)$.

\begin{lemma}\label{psiform}
$\Psi(x) = \nabla H(Q(x))$, where $H$ is the conserved quantity from \cref{hyp:H}. In addition, $\Psi(-x) = R \Psi(x)$, where $R$ is the standard reversor operator, and the last component of $\Psi(x)$ is $q'(x)$.
\begin{proof}
Differentiating $\langle F(U; c), \nabla_U H(U; c) \rangle = 0$, 
\[
0 = D F(Q(x))^* \nabla H(Q(x)) + D^2 H(Q(x))^* F(Q(x)).
\]
Using standard vector calculus identities, equation \cref{eqODE}, and the fact that the Hessian is self-adjoint,
\begin{align*}
-D F(Q(x))^* \nabla H(Q(x)) &= D^2 H(Q(x)) Q'(x) = \frac{d}{dx} \nabla H(Q(x)),
\end{align*}
thus $\nabla H(Q(x))$ is a solution to \cref{adjvareq1}. Since $\nabla H$ is continuous and $Q(x)$ is exponentially localized, $\nabla H(Q(x))$ is bounded, thus by uniqueness we can take $\Psi(x) = \nabla H(Q(x))$. Using \cref{genODErev} and the symmetry relation $Q(-x) = R Q(x)$, 
\begin{align*}
[R \Psi(-x)]' = -R DF(Q(-x)) \Psi(-x) 
= -R [-RDF(Q(x))R] \Psi(-x) = DF(Q(x))[ R \Psi(-x) ],
\end{align*}
thus $\Psi(-x) = R \Psi(x)$ by uniqueness. By the definition of $F$, if $\Psi(x) = (\psi_1(x), \dots, \psi_{2m}(x))$ is a solution to \cref{adjvareq1}, then $\psi_{2m}(x)$ solves $\calL(q)^* w = [\calE''(q) + c]^* w = 0$. Since $\calL(q)$ is self-adjoint and $\calL(q)q' = 0$, $\psi_{2m}(x) = q'(x)$.
\end{proof}
\end{lemma}

In the next lemma, we collect a few important results about solutions to \cref{vareq1} and \cref{adjvareq1}.

\begin{lemma}\label{eigadjoint}
Consider the linear ODE $V' = A(x)V$ and the corresponding adjoint equation $W' = -A(x)^* W$, where $A(x)$ is a smooth $n \times n$ matrix. Then
\begin{enumerate}[(i)]
\item $\dfrac{d}{dx}\langle V(x), W(x) \rangle = 0$, thus the inner product is constant in $x$.
\item If $W(x)$ is bounded, and $V(x) \rightarrow 0$ as $x \rightarrow \infty$ or $x \rightarrow -\infty$, then $\langle V(x), W(x) \rangle = 0$ for all $x \in \R$. The same holds if we reverse the roles of $W$ and $V$.
\item If $\Phi(y, x)$ is the evolution operator for $V'(x) = A(x)V(x)$, then $\Phi(x, y)^*$ is the evolution operator for the adjoint equation $W'(y) = -A(y)^* W(y)$.
\end{enumerate}
\begin{proof}
For part (i), 
\begin{align*}
\dfrac{d}{dx}\langle V(x), W(x) \rangle &= 
\langle V'(x), W(x) \rangle + \langle V(x), W'(x) \rangle \\
&= \langle A(x)V(x), W(x) \rangle + \langle V(x), -A(x)^* W(x) \rangle = 0.
\end{align*}
Part (ii) follows from part (i), the Cauchy-Schwartz inequality and the continuity of the inner product. For part (iii), take the derivative of $\Phi(y, x)\Phi(x, y) = I$ with respect to $y$ to get
\begin{align*}
0 &= \left(\frac{d}{dy}\Phi(y, x)\right) \Phi(x, y) +
\Phi(y, x)\left(\frac{d}{dy}\Phi(x, y)\right) 
= A(y) + \Phi(y, x)\left(\frac{d}{dy}\Phi(x, y)\right).
\end{align*}
Rearranging and taking the transpose of both sides yields
\begin{align*}
\frac{d}{dy}\Phi(x, y)^* &= -A(y)^* \Phi(x, y)^*.
\end{align*}
\end{proof}
\end{lemma}

\noi By \cref{eigadjoint}(ii), $\Psi(0) \perp \R Q'(0) \oplus Y^+ \oplus Y^-$, thus we can decompose $\R^{2m}$ as
\begin{equation}\label{R2mdecomp}
\R^{2m} = \R Q'(0) \oplus Y^+ \oplus Y^- \oplus \R \Psi(0).
\end{equation}

\subsection{Piecewise ansatz}

First, we write the unstable and stable manifolds as graphs over their tangent spaces. Following \cite{Sandstede1997}, we can parameterize $\tilde{W}^u(0)$ and $\tilde{W}^s(0)$ near $Q(0)$ by the smooth functions $Q^-(\gamma, \beta^-)$ and $Q^+(\gamma, \beta^+)$, where $\gamma \in \R$ and $\beta^\pm \in Y^\pm$. These functions are chosen so that $Q^+(\gamma, 0) - Q^-(\gamma, 0) \in \R \Psi(0)$, and $Q^+(0, 0) = Q^-(0, 0) = Q(0)$. We will always take $\gamma = 0$. Let $Q^\pm(x; \beta^\pm)$ be the unique solutions to \cref{genODE} on $\R^\pm$ with initial conditions $Q^\pm(0, \beta^\pm)$ at $x = 0$. $Q^-(x; \beta^-)$ lies in the unstable manifold $\tilde{W}^u(0)$ and $Q^+(x; \beta^+)$ lies in the stable manifold $\tilde{W}^s(0)$. 

We will look for a $n-$periodic solution $U(n)$ to \cref{genODE} which is piecewise of the form
\begin{equation}\label{Upiecewise}
\begin{aligned}
U_i^-(x) &= Q^-(x; \beta_i^-) + \tilde{Q}_i^-(x) \qquad\qquad x \in [-X_{i-1}, 0] \\
U_i^+(x) &= Q^+(x; \beta_i^+) + \tilde{Q}_i^+(x) \qquad\qquad x \in [0, X_i]
\end{aligned}
\end{equation}
for $i = 0, \dots, n-1$, where $U_i^-: [-X_{i-1}, 0] \rightarrow \R$ and $U_i^+: [0, X_i] \rightarrow \R$ are continuous. The subscripts $i$ are taken $\Mod n$ since we are on a periodic domain, and the pieces are glued together end-to-end as in \cite{Sandstede1998}, with one additional join needed to ``close the loop''. Since $Q^\pm(0; \beta_i^\pm) \in \R Q'(0) \oplus Y^\pm$, we are free to choose $\tilde{Q}_i^\pm(x)$ so that
\begin{align*}
\tilde{Q}_i^-(0) &\in \R \Psi(0) \oplus Y^- \\
\tilde{Q}_i^+(0) &\in \R \Psi(0) \oplus Y^+.
\end{align*}
To construct a periodic $n$-pulse, we will solve the following system of equations
\begin{align}
(U_i^\pm(x))' - F(U_i^\pm(x)) &= 0 \label{exsystem1} \\
U_i^+(X_i) - U_{i+1}^-(-X_i) &= 0 \label{exsystem2} \\
U_i^+(0) - U_i^-(0) &= 0 \label{exsystem3}
\end{align}
for $i = 0, \dots, n-1$. Equation \cref{exsystem2} is a matching condition at the pulse tails, and equation \cref{exsystem3} is a matching condition at the pulse centers.

\subsection{Exponential Dichotomy}\label{sec:existdichot}

Let $\Phi_\pm(x, y; \beta^\pm)$ be the family of evolution operators for
\begin{align}\label{qpmODEs}
[V^\pm(x)]' &= D F\left(Q^\pm(x, \beta^\pm)\right) V^\pm(x) && x \in \R^\pm.
\end{align}
Choose any $\alpha$ slightly less than $\alpha_0$. In the next lemma, we decompose these evolution operators in exponential dichotomies on $\R^+$ and $\R^-$. 

\begin{lemma}\label{dichotomy1}
There exist projections
\begin{align*}
&P_+^s(y; \beta^+), \quad P_+^u(y; \beta^+) = I - P_+^s(y; \beta^+) && y \geq 0 \\
&P_-^u(y; \beta^-), \quad P_-^s(y; \beta^-) = I - P_-^u(y; \beta^-) && y \leq 0 
\end{align*}
on $\R^\pm$ such that the evolution operators $\Phi_\pm(x, y; \beta^\pm)$ can be decomposed as $\Phi_\pm(x, y; \beta^\pm) = \Phi^s_\pm(x, y; \beta^\pm) + \Phi^u_\pm(x, y; \beta^\pm)$, where
\begin{align*}
\Phi^s_\pm(x, y; \beta^\pm) &= \Phi_\pm(x, y; \beta^\pm) P^s_\pm(y; \beta^\pm) \\
\Phi^u_\pm(x, y; \beta^\pm) &= \Phi_\pm(x, y; \beta^\pm) P^u_\pm(y; \beta^\pm).
\end{align*}
We have the estimates
\begin{align*}
|\Phi^s_+(x, y, \beta^+)| &\leq C e^{-\alpha(x - y)} && 0 \leq y \leq x \\
|\Phi^u_+(x, y, \beta^+)| &\leq C e^{-\alpha(y - x)} && 0 \leq x \leq y \\
|\Phi^u_-(x, y, \beta^-)| &\leq C e^{-\alpha(y - x)} && 0 \geq y \geq x \\
|\Phi^s_-(x, y, \beta^-)| &\leq C e^{-\alpha(x - y)} && 0 \geq x \geq y,
\end{align*}
which also hold for derivatives with respect to the initial conditions $\beta^\pm$. In addition, the projections satisfy the commuting relations
\begin{align*}
\Phi_\pm(x, y; \beta^\pm) P^{s/u}_\pm(y; \beta^\pm) 
= P^{s/u}_\pm(x; \beta^\pm) \Phi_\pm(x, y; \beta^\pm).
\end{align*}
The projections can be chosen so that at $y = 0$ we have, independent of $\beta^+$ and $\beta^-$,
\begin{align*}
\ker P^s_+(0; \beta^+) &= \R \Psi(0) \oplus Y^-, \qquad
\ker P^u_-(0; \beta^-) = \R \Psi(0) \oplus Y^+ \\
\ran P^u_+(0; \beta^+) &= \R \Psi(0) \oplus Y^-, \qquad
\ran P^s_-(0; \beta^-) = \R \Psi(0) \oplus Y^+.
\end{align*}
Let $E_0^s$ and $E_0^u$ be the stable and unstable eigenspaces of $DF(0)$, and let $P_0^s$ and $P_0^u$ be the corresponding eigenprojections. For any $\alpha$ with $0 < \alpha < \alpha_0$, we have the following estimates, which are independent of $\beta_i^\pm$.
\begin{equation}\label{projdiffest}
\begin{aligned}
&|P^u_+(x; \beta^+) - P_0^u| \leq C e^{-\alpha x},
&&|P^s_+(x; \beta^+) - P_0^s| \leq C e^{-\alpha x} \\
&|P^u_-(x; \beta^-) - P_0^u| \leq C e^{\alpha x},
&&|P^s_-(x; \beta^-) - P_0^s| \leq C e^{\alpha x} .
\end{aligned}
\end{equation}

\begin{proof}
Since $DF(0)$ is hyperbolic by \cref{hyp:hypeq}, and $|\Re \nu| \geq \alpha_0$ for all eigenvalues $\nu$ of $DF(0)$, the exponential dichotomy results follow from \cite[Lemma 5.1]{Sandstede1997}, which follows from \cite[Lemma 1.1]{Sandstede1993}. The estimates \cref{projdiffest} follow from \cite[Lemma 1.1]{Sandstede1993} and \cite[Lemma 2.1]{Sandstede1993}.
\end{proof}
\end{lemma}

\subsection{Fixed Point Formulation}

Next, we formulate equation \cref{exsystem1} as a fixed point problem. Substituting the piecewise ansatz \cref{Upiecewise} into \cref{exsystem1}, and using the fact that $Q^\pm(x; \beta_i^\pm)$ solves \cref{exsystem1} on $\R^\pm$,
\begin{align*}
(\tilde{Q}_i^\pm(x))' &= F(Q^\pm(x; \beta_i^\pm) + \tilde{Q}_i^\pm(x)) - F(Q^\pm(x; \beta_i^\pm)) && i = 0, \dots, n-1.
\end{align*}
Expanding the RHS in a Taylor series about $Q^\pm(x; \beta_i^\pm)$, this becomes
\begin{align}\label{Vpiecewise1}
(\tilde{Q}_i^\pm(x))' &= DF(Q^\pm(x; \beta_i^\pm)) \tilde{Q}_i^\pm(x) + G_i^\pm(x; \beta_i^\pm) && i = 0, \dots, n-1,
\end{align}
where $G_i^\pm(x; \beta_i^\pm) = \mathcal{O}(|\tilde{Q}_i^\pm(x)|^2)$. As in \cite{Sandstede1997}, derivatives of $G_i^\pm$ with respect to the parameters $\beta_i^\pm$ are also quadratic in $\tilde{Q}_i^\pm$. Using the exponential dichotomy, we rewrite \cref{Vpiecewise1} in integrated form as the fixed point equations
\begin{equation}\label{FPequations}
\begin{aligned}
\tilde{Q}_i^+(x) &= \Phi^u_+(x, X_i; \beta_i^+) a_i^+ \\
&\qquad+ \int_{X_i}^x \Phi_+^u(x, y; \beta_i^+) G_i^+(y; \beta_i^+)dy 
+ \int_0^x \Phi_+^s(x, y; \beta_i^+) G_i^+(y; \beta_i^+)dy \\ 
\tilde{Q}_i^-(x) &= \Phi^s_-(x, -X_{i-1}; \beta_i^-) a_{i-1}^- \\ 
&\qquad+\int_{-X_{i-1}}^x \Phi_-^s(x, y; \beta_i^-) G_i^-(y; \beta_i^-)dy 
+ \int_0^x \Phi_-^u(x, y; \beta_i^-) G_i^-(y; \beta_i^-)dy,
\end{aligned}
\end{equation}
for $i = 0, \dots, n-1$, where $a_i^+ \in E_0^u$ and $a_i^- \in E_0^s$. Define the exponentially weighted norms
\begin{equation}\label{expwtnorm}
\begin{aligned}
\|V\|_{X, +} &= \sup_{x \in [0, X]} e^{\alpha(X - x)}|V(x)| \\
\|V\|_{X, -} &= \sup_{x \in [-X, 0]} e^{\alpha(X + x)}|V(x)|,
\end{aligned}
\end{equation}
and let $K_{X, \pm}$ be the Banach spaces of continuous functions on $[0, X]$ and $[-X, 0]$ equipped with these norms. Let $B_{X, \pm}(\rho)$ be the ball of radius $\rho$ about $0$ in $K_{X, \pm}$.

\subsection{Inversion}

As in \cite{Sandstede1997}, we will solve for the remainder functions $\tilde{Q}_i^\pm$ and parameters $\beta_i^\pm$ in a series of lemmas. First, we will solve equation \cref{FPequations} for $\tilde{Q}_i^\pm$.

\begin{lemma}\label{solveforV}
There exist $\delta, \rho > 0$ such that for $X_i > 1/\delta$ and $|a_i^\pm|, |\beta_i^\pm| < \delta$, where $i = 0, \dots, n-1$, there exist unique solutions $\tilde{Q}_i^- \in B_{X_{i-1}, -}(\rho)$ and $\tilde{Q}_i^+ \in B_{X_i, +}(\rho)$ to \cref{FPequations}. These depend smoothly on $(a_{i-1}^-, \beta_i^-)$ and $(a_i^+, \beta_i^+)$, respectively, and we have the estimates
\begin{equation}\label{Vest}
\begin{aligned}
\|\tilde{Q}_i^-\|_{X_{i-1}, -} &\leq C |a_{i-1}^-| \\
\|\tilde{Q}_i^+\|_{X_i, +} &\leq C |a_i^+|,
\end{aligned}
\end{equation}
where the constant $C$ depends only on $\delta$. The estimates hold for derivatives of $\tilde{Q}_i^\pm$ with respect to $\beta_i^\pm$.
\begin{proof}
The proof follows \cite[Lemma 5.2]{Sandstede1997}. For the first term on the RHS of \cref{FPequations},
\begin{align*}
e^{\alpha(X_i - x)} | \Phi^u_+(x, X_i; \beta_i^+) a_i^+ | 
&\leq C e^{\alpha(X_i - x)} e^{-\alpha(X_i - x)} |a_i^+| = C |a_i^+|,
\end{align*}
and for the second term, since $G_i^+$ is quadratic in $\tilde{Q}_i^+$ and $\tilde{Q}_i^+ \in K_{X_i, +}$, 
\begin{align*}
e^{\alpha(X_i - x)} &\left| \int_{X_i}^x \Phi_+^u(x, y; \beta_i^+) G_i^+(y; \beta_i^+)dy  \right| 
\leq C e^{\alpha(X_i - x)} \int_x^{X_i} e^{-\alpha(y - x)}|\tilde{Q}_i^+(y)|^2 dy \\
&\leq C e^{\alpha(X_i - x)} \int_x^{X_i} 
e^{-\alpha(y - x)}(e^{-\alpha(X_i - y)})^2|e^{\alpha(X_i - y)} \tilde{Q}_i^+(y)|^2 dy \leq C.
\end{align*}
The third term is similarly bounded. Thus the RHS of the fixed point equation \cref{FPequations} is a smooth map $K_{X_i, +} \mapsto K_{X_i, +}$. Define $H_i^+: K_{X_i, +} \times E_0^s \times Y^+ \rightarrow K_{X_i, +}$ by
\begin{align*}
H_i^+(\tilde{Q}_i^+(x), &a_i^+, \beta_i^+) = \tilde{Q}_i^+(x) - \Phi^u_+(x, X_i; \beta_i^+) a_i^+ - \int_{X_i}^x \Phi_+^u(x, y; \beta_i^+) G_i^+(y; \beta_i^+)dy \\
&- \int_0^x \Phi_+^s(x, y; \beta_i^+) G_i^+(y; \beta_i^+)dy.
\end{align*}
Since $Q(x)$ satisfies \cref{genODE}, $H_i^+(0, 0, 0) = 0$, and since $G_i^+$ is quadratic in $\tilde{Q}_i^+(x)$, the Fr\'echet derivative of $H_i^+$ with respect to $\tilde{Q}_i^+(x)$ at $(\tilde{Q}_i^+(x), a_i^+, \beta_i^+) = (0, 0, 0)$ is the identity. Using the implicit function theorem, we can solve for $\tilde{Q}_i^+(x)$ in terms of $(a_i^+, \beta_i^+)$ for sufficiently small $|a_i^+|$ and $|\beta_i^+|$. Since the map $H_i^+$ is smooth, this dependence is smooth. The estimate on $\tilde{Q}_i^+$ comes from the first term on the RHS of \cref{FPequations}, since the remaining terms are quadratic in $\tilde{Q}_i^+$. Since the exponential dichotomy estimates from \cref{dichotomy1} hold for derivatives with respect to $\beta_i^+$, these estimates do as well. We can similarly solve for $\tilde{Q}_i^-$ in terms of $(a_{i-1}^-, \beta_i^-)$.
\end{proof}
\end{lemma}

Next, we will solve equation \cref{exsystem2} to match the pieces \cref{Upiecewise} at $\pm X_i$ and obtain the initial conditions $a_i^\pm$.

\begin{lemma}\label{solvefora}
For $X_i$ chosen as in \cref{solveforV}, and for $i = 0, \dots, n-1$, there is a unique pair of initial conditions $(a_i^+, a_i^-) \in E_0^s \times E_0^u$ such that $U_i^+(X_i) - U_{i+1}^-(-X_i) = 0$. The pair $(a_i^+, a_i^-)$ depends smoothly on $(\beta_i^+, \beta_{i+1}^-)$, and we have the estimate
\begin{equation}\label{aest}
|a_i^\pm| \leq C e^{-\alpha X_i},
\end{equation}
which holds as well for derivatives with respect to $\beta_i^\pm$. In addition,
\begin{equation}\label{aiformula}
\begin{aligned}
a_i^+ &= -P^u_0 \left( Q^+(X_i; \beta_i^+) - Q^-(-X_i; \beta_{i+1}^-) \right) + \mathcal{O}( e^{-2 \alpha X_i} ) \\
a_i^- &= P^s_0 \left( Q^+(X_i; \beta_i^+) - Q^-(-X_i; \beta_{i+1}^-) \right) + \mathcal{O}( e^{-2 \alpha X_i}).
\end{aligned}
\end{equation}

\begin{proof}
Evaluating the fixed point equations \cref{FPequations} at $\pm X_i$ and substituting them into \cref{Upiecewise}, the matching condition \cref{exsystem2} can be written as $H_i(a_i^+, a_i^-, \beta_i^+, \beta_{i+1}^-) = 0$, where $H_i: E_0^s \times E_0^u \times Y^+ \times Y^- \rightarrow \R^{2m}$ is defined by
\begin{align*}
H_i(a_i^+, &a_i^-, \beta_i^+, \beta_{i+1}^-) 
= a_i^+ - a_i^- + (P^u_+(X_i; \beta_i^+) -  P^u_0)a_i^+ - (P^s_-(-X_i; \beta_{i+1}^-) - P^s_0) a_i^-  \\
&+ Q^+(X_i; \beta_i^+) - Q^-(-X_i; \beta_{i+1}^-)\\
&+ \int_0^{X_i} \Phi_+^s(X_i, y; \beta_i^+) G_i^+(y,\beta_i^+)dy
- \int_0^{-X_i} \Phi_-^u(-X_i, y; \beta_{i+1}^-) G_{i+1}^-(y,\beta_{i+1}^-)dy,
\end{align*}
where we have substituted $\tilde{Q}_i^\pm(x)$ from \cref{solveforV} into $G_i^\pm$. Since $Q(x)$ satisfies \cref{genODE}, $H_i(0, 0, 0, 0) = 0$. Since $G_i^\pm$ is quadratic in $\tilde{Q}_i^\pm$, thus quadratic in $a_i^\pm$ by \cref{solveforV},
\[
\frac{\partial}{\partial a_i^\pm} H_i(0, 0, 0, 0) = \pm 1 + \mathcal{O} (e^{-\alpha X_i}),
\]
where we also used the estimate \cref{projdiffest}. For sufficiently large $X_i$, $D_{a_i^\pm} H_i(0, 0, 0, 0)$ is invertible in a neighborhood of $(0, 0, 0, 0)$, thus we can use the implicit function theorem to solve for $a_i^\pm$ in terms of $\beta_i^\pm$. The estimate \cref{aest} then comes from the stable manifold theorem, since $Q^\pm(\pm X_i; \beta_i^\pm) = \mathcal{O}(e^{-\alpha X_i})$. To obtain the expressions \cref{aiformula}, we apply the eigenprojections $P^u_0$ and $P^s_0$ (respectively) to $H_i(a_i^+, a_i^-, \beta_i^+, \beta_{i+1}^-) = 0$. The bound on the remainder term comes from the bound \cref{aest}, together with the estimates from \cref{solveforV} and equation \cref{projdiffest}. 
\end{proof}
\end{lemma}

It remains to solve equation \cref{exsystem3}, which is the matching condition at $x = 0$. Before doing that, we will use the flow-box method to make a smooth change of coordinates which will ``straighten out'' the stable and unstable manifolds near $Q(0)$ so that their non-intersecting directions are $Y^+$ and $Y^-$.

\begin{lemma}\label{straightenW}
There exists a differentiable map $S: \R \times Y^- \times Y^+ \times \R \Psi(0) \rightarrow \R^{2m}$ such that $S(0, 0, 0, 0) = Q(0)$, $S$ is invertible in a neighborhood of $Q(0)$, and for sufficiently small $\beta^\pm$,
\begin{align*}
S^{-1}(Q^-(0; \beta^-)) &= \beta^-  \\
S^{-1}(Q^+(0; \beta^+)) &= \beta^+.
\end{align*}
\begin{proof}
Let $\Theta_x(U_0)$ be the solution operator which maps $U_0 \in \R^{2m}$ to the point $U(x)$, where $U(\cdot)$ is the unique solution to \cref{genODE} with $U(0) = U_0$. Define the map $S: \R \times Y^- \times Y^+ \times \R \Psi(0) \rightarrow \R^{2m}$ by 
$S(x; \beta^-, \beta^+, \gamma) = \Theta_x\left(Q(0) + Q^-(0; \beta^-) + Q^-(0; \beta^+) + \gamma \Psi(0)\right)$. 
For small $x$ and $\beta^\pm$, the stable and unstable manifolds are the surfaces $\tilde{W}^s = S(x; 0, \beta^+, 0)$ and $\tilde{W}^u = S(x; \beta^-, 0, 0)$. Their one-dimensional intersection is the homoclinic orbit $Q(x) = S(x; 0, 0, 0)$, and $S(0, 0, 0, 0) = Q(0) \neq 0$. The partial derivatives of $S$ are
\begin{align*}
S_x(0, 0, 0, 0) &= F(Q(0)) = Q'(0) \\
S_{\beta^-}(0, 0, 0, 0) &= (Q^-)_{\beta^-}(0; 0) = Y^- \\
S_{\beta^+}(0, 0, 0, 0) &= (Q^+)_{\beta^+}(0; 0) = Y^+ \\
S_{\gamma}(0, 0, 0, 0) &= \Psi(0),
\end{align*}
which span $\R^{2m}$ by \cref{R2mdecomp}. Since the Jacobian of $S$ is invertible at the origin, $S$ is invertible near $Q(0)$ by the inverse function theorem.
\end{proof}
\end{lemma}

After applying this coordinate change near $Q(0)$, the matching condition \cref{exsystem3} is equivalent to projecting $U_i^+(0) - U_i^-(0) = 0$ onto $\R \Q'(0)$, $Y^+$, $Y^-$, and $\R \Psi(0)$ and solving separately on each subspace. Since $P_{\R Q'(0)}(Q^\pm(0; \beta^\pm)) = 0$ and $\tilde{Q}_i^\pm(0) \in \R \Psi(0) \oplus Y^+ \oplus Y^-$, the equation $P_{\R Q'(0)}(U_i^+(0) - U_i^-(0)) = 0$ is automatically satisfied. Since $P_{Y^\pm}(Q^\pm(0; \beta^\pm)) = \beta^\pm$ and $P_{\R \Psi(0)}(Q^\pm(0; \beta^\pm)) = 0$ due to the change of coordinates, it remains to solve the equations
\begin{align}
P_{Y^+}(\tilde{Q}_i^+(0) - \tilde{Q}_i^-(0)) + \beta_i^+  &= 0 \label{PY+match} \\
P_{Y^-}(\tilde{Q}_i^+(0) - \tilde{Q}_i^-(0)) - \beta_i^- &= 0 \label{PY-match} \\
P_{\R \Psi(0)}(\tilde{Q}_i^+(0) - \tilde{Q}_i^-(0)) &= 0. \label{PZmatch}
\end{align}
In the next lemma we solve \cref{PY+match} and \cref{PY-match} to obtain the parameters $\beta_i^\pm$.

\begin{lemma}\label{solveforbeta}
For $X_i$ chosen as in \cref{solveforV}, and for $i = 0, \dots, n-1$, there exist $(\beta_i^+, \beta_i^-) \in Y^+ \times Y^-$ such that $P_{Y^+ \oplus Y^-}(U_i^+(0) - U_i^-(0)) = 0$. In addition,
\begin{equation}\label{betaest}
\begin{aligned}
|\beta_i^+| &\leq C e^{-2 \alpha X_{i-1}} \\
|\beta_i^-| &\leq C e^{-2 \alpha X_i}.
\end{aligned}
\end{equation}
\begin{proof}
Evaluating the fixed point equations \cref{FPequations} at $0$ and substituting them into \cref{Upiecewise}, equations \cref{PY+match} and \cref{PY-match} can be written as $H_i(\beta_i^+, \beta_i^-) = 0$, where $H_i: Y^+ \oplus Y^- \rightarrow Y^+ \oplus Y^-$ is defined by
\begin{equation}\label{defHPY}
H_i(\beta_i^+, \beta_i^-) = 
\begin{pmatrix}
\beta_i^+ - P_{Y^+}\left(\Phi^s_-(0, -X_{i-1}, \beta_i^-) a_{i-1}^- 
- \int_{-X_{i-1}}^0 \Phi_-^s(0, y, \beta_i^-) G_i^-(y; \beta_i^-) dy\right) \\
\beta_i^- + P_{Y^-}\left( \Phi^u_+(0, X_i; \beta_i^+) a_i^+ 
+ \int_{X_i}^0 \Phi_+^u(0, y; \beta_i^+) G_i^+(y; \beta_i^+)dy \right)
\end{pmatrix},
\end{equation}
where we have substituted our expressions for $\tilde{Q}_i^\pm$ and $a_i^\pm$ from \cref{solveforV} and \cref{solvefora}. Using the estimates from these lemmas together with \cref{dichotomy1},
\begin{equation}\label{DHexp}
D H_i(\beta_i^+, \beta_i^-) = 
\begin{pmatrix}
1 & \mathcal{O}(e^{-2 \alpha X_{i-1}} ) \\
\mathcal{O}(e^{-2 \alpha X_i}) &  1 
\end{pmatrix},
\end{equation}
which is independent of $\beta_i^\pm$, thus $D H_i(\beta_i^+, \beta_i^-)$ is invertible for sufficiently large $X_i$. By the inverse function theorem, $(\beta_i^+, \beta_i^-) = H_i^{-1}(0, 0)$. The estimates \cref{betaest} follow from \cref{defHPY} and Lemmas \ref{dichotomy1}, \ref{solveforV}, and \ref{solvefora}.
\end{proof}
\end{lemma}

We have found a unique solution to \cref{exsystem1} and \cref{exsystem2} such that \cref{exsystem3} is satisfied except for $n$ jumps in the direction of $\Psi(0)$. We summarize what we have obtained so far in the following lemma.

\begin{lemma}\label{solvewithjumps}
There exists $\delta > 0$ such that for $|X_i| \geq 1/\delta$, where $i = 0, \dots, n-1$, there is a unique solution $U(x)$ to equations \cref{exsystem1}, \cref{exsystem2}, and \cref{exsystem3} which is continuous except for $n$ jumps in the direction of $\Psi(0)$. $U(x)$ can be written piecewise in the form 
\begin{equation}\label{Upiecewise2}
\begin{aligned}
U_i^-(x) &= Q^-(x; \beta_i^-) + \tilde{Q}_i^-(x) \qquad\qquad x \in [-X_{i-1}, 0] \\
U_i^+(x) &= Q^+(x; \beta_i^+) + \tilde{Q}_i^+(x) \qquad\qquad x \in [0, X_i],
\end{aligned}
\end{equation}
where the pieces are glued together end-to-end in a loop, and we have the estimates
\begin{enumerate}[(i)]
\item
\begin{equation}\label{tildeQbounds}
\begin{aligned}
|\tilde{Q}_i^-(x)| &\leq C e^{-\alpha(X_{i-1} + x)}e^{-\alpha X_{i-1}} \\
|\tilde{Q}_i^+(x)| &\leq C e^{-\alpha(X_i - x)}e^{-\alpha X_i}.
\end{aligned}
\end{equation}
\item 
\begin{equation}
\begin{aligned}\label{Qpmbounds}
|Q^-(x; \beta_i^-) - Q(x)| &\leq C e^{-2 \alpha X_i} e^{\alpha x} \\
|Q^+(x; \beta_i^+) - Q(x)| &\leq C e^{-2 \alpha X_{i-1}} e^{-\alpha x}.
\end{aligned}
\end{equation}
\item
\begin{equation}\label{VQpm}
\begin{aligned}
\tilde{Q}_i^+(X_i) &= Q^-(-X_i; \beta_{i+1}^-) + \mathcal{O}(e^{-2 \alpha X_i}) \\
\tilde{Q}_{i+1}^-(-X_i) &= Q^+(X_i; \beta_i^+) + \mathcal{O}(e^{-2 \alpha X_i}).
\end{aligned}
\end{equation}
\end{enumerate}
These estimates hold in addition for derivatives with respect to $x$.
\begin{proof}
Part (i) follows from the estimates \cref{Vest} and \cref{aest} together with the definition of the exponentially weighted norm \cref{expwtnorm}. Part (ii) follows from the estimate \cref{betaest}, smooth dependence on initial conditions, and the stable manifold theorem. For part (iii), we solved the matching condition $Q^+(X_i; \beta_i^+) + \tilde{Q}_i^+(X_i) = Q^-(-X_i; \beta_{i+1}^-) + \tilde{Q}_i^-(-X_i)$ in \cref{solvefora}. Applying the projections $P^u_-(-X_i, \beta_{i+1}^-)$ and $P^s_+(X_i, \beta_i^+)$ in turn to this and using \cref{projdiffest}, \cref{FPequations}, and the estimates from the previous lemmas in this section, we obtain the estimates \cref{VQpm}.
\end{proof}
\end{lemma}

\subsection{Jump conditions}

Equation \cref{PZmatch} gives us $n$ jump conditions in the direction of $\Psi(0)$. As in \cite{SandstedeStrut,Sandstede1998} these will only be satisfied for certain values of the pulse distances $X_i$. Since we have a conservative system, we only need to satisfy $n-1$ of these jump conditions, which are given in the following lemma.

\begin{lemma}\label{jumplemma1}
A periodic $n$-pulse solution exists if and only if for $i = 0, \dots, n-2$,
\begin{align}\label{jumpcond2}
\xi_i = \langle \Psi(-X_i), Q(X_i) \rangle - \langle \Psi(-X_{i-1}), Q(X_{i-1}) \rangle + R_i = 0,
\end{align}
where the remainder term has bound
\begin{align*}
|R_i| \leq C ( e^{-3 \alpha X_i} +  e^{-3 \alpha X_{i-1}}).
\end{align*}
\begin{proof}
Evaluating the fixed point equations \cref{FPequations} at $x=0$ and substituting them into \cref{PZmatch},
\begin{equation}\label{PZmatch2}
\begin{aligned}
\langle \Psi(0), &\tilde{Q}_i^+(0) - \tilde{Q}_i^-(0) \rangle = \langle \Psi(0), \Phi^u_+(0, X_i; \beta_i^+) a_i^+ \rangle
- \langle \Psi(0), \Phi^s_-(0, -X_{i-1}, \beta_i^-) a_{i-1}^- \rangle \\
&+ \int_{X_i}^0 \langle \Psi(0), \Phi_+^u(0, y; \beta_i^+) G_i^+(y; \beta_i^+) \rangle dy - \int_{-X_{i-1}}^0 \langle \Psi(0), \Phi_-^s(0, y, \beta_i^-) G_i^-(y; \beta_i^-) \rangle dy,
\end{aligned}
\end{equation}
where we have substituted our expressions for $\tilde{Q}_i^\pm$, $a_i^\pm$, and $\beta_i^\pm$ from \cref{solveforV}, \cref{solvefora}, and \cref{solveforbeta}. Using \cref{dichotomy1}, smooth dependence on initial conditions, and the bound \cref{betaest},
\begin{equation}\label{phibetaest}
\begin{aligned}
|\Phi_+^u(0, x; \beta_i^+) - \Phi_+^u(0, x; 0)| &\leq C e^{-2 \alpha X_{i-1}} e^{-\alpha x} && x \geq 0\\
|\Phi_-^s(0, x; \beta_i^-) - \Phi_-^s(0, x; 0)| &\leq C e^{-2 \alpha X_i}e^{\alpha x} &&  x \leq 0.
\end{aligned}
\end{equation}
For the term involving $a_i^+$ in \cref{PZmatch2}, we substitute \cref{aiformula} from \cref{solvefora} and use \cref{phibetaest}, the estimate \cref{Qpmbounds}, and \cref{eigadjoint}(iii) to obtain
\begin{align*}
\langle \Psi&(0), \Phi^u_+(0, X_i; \beta_i^+) a_i^+ \rangle \\
&= -\langle \Psi(0), \Phi^u_+(0, X_i; \beta_i^+)( P^u_0 ( Q^+(X_i; \beta_i^+) - Q^-(-X_i; \beta_{i+1}^-)) + \mathcal{O}( e^{-2 \alpha X_i} ) ) \rangle \\
&= -\langle \Psi(0), \Phi^u_+(0, X_i; 0) P^u_0 ( Q^+(X_i; \beta_i^+) - Q^-(-X_i; \beta_{i+1}^-)) \rangle + \mathcal{O}( e^{-3 \alpha X_i} + e^{-2\alpha X_{i-1}}e^{-\alpha X_i} ) \\
&= -\langle \Psi(X_i), P^u_+(X_i; 0) P^u_0 ( Q(X_i) - Q(-X_i) ) \rangle + \mathcal{O}( e^{-3 \alpha X_i} + e^{-3\alpha X_{i-1}}).
\end{align*}
By \cref{projdiffest}, $P^u_+(X_i; 0) P^u_0 ( Q(X_i) - Q(-X_i)) = Q(-X_i) + \mathcal{O}(e^{-2\alpha X_i})$, thus it follows from \cref{eigadjoint} that
\begin{align*}
\langle \Psi(0), \Phi^u_+(0, X_i; \beta_i^+) a_i^+ \rangle = 
\langle \Psi(X_i), Q(-X_i) \rangle + \mathcal{O}( e^{-3 \alpha X_i} + e^{-3\alpha X_{i-1}}).
\end{align*}
Similarly, we have
\begin{align*}
\langle \Psi(0), \Phi^s_-(0, -X_{i-1}, \beta_i^-) a_{i-1}^- \rangle = 
\langle \Psi(-X_{i-1}), Q(-X_{i-1}) \rangle + \mathcal{O}( e^{-3 \alpha X_i} + e^{-3\alpha X_{i-1}}).
\end{align*}
For the integral terms, we use the estimates from \cref{solvewithjumps} and the fact that $G_i^\pm$ is quadratic in $\tilde{Q}_i^\pm$ to obtain the estimate
\begin{align*}
\left| \int_{X_i}^0 \langle \Psi(0), \Phi_+^u(0, y; \beta_i^+) G_i^+(y; \beta_i^+) \rangle dy \right| \leq C e^{-3\alpha X_i}.
\end{align*}
The other integral term has a similar bound. By reversibility,
\begin{align*}
\langle \Psi(X_i), Q(-X_i) \rangle 
&= \langle R\Psi(-X_i), R Q(X_i) \rangle 
= \langle \Psi(-X_i), Q(X_i) \rangle.
\end{align*}
Combining everything above, we obtain equation \cref{jumpcond2} and the remainder bound for $R_i$. As in \cite[p. 2093]{SandstedeStrut}, since \cref{eqODE} is a conservative system, if $n-1$ of the jump conditions are satisfied, the final jump condition must automatically be satisfied. Since it does not matter which condition we eliminate, we choose to eliminate the last one. 
\end{proof}
\end{lemma}

As a corollary, periodic single pulse solutions exist for sufficiently large $X_0$, since in that case there are no jump conditions. These are unimodal periodic orbits which are close to the primary homoclinic orbit $Q(x)$.

\begin{corollary}\label{corr:1pexists}
Periodic single pulse solutions exist for sufficiently large $X_0$.
\end{corollary}

\subsection{Rescaling and parameterization}

Following \cite[Section 6]{Sandstede1998}, we will introduce a change of variables with a built-in scaling parameter to facilitate the analysis. Define the set
\begin{align}
\mathcal{R} &= \left\{ \exp\left(-\frac{2 \pi m}{\rho}\right) : m \in \N_0 \right\} \cup \{ 0 \},
\end{align}
where $\rho = \beta_0 / \alpha_0$. Since $\mathcal{R}$ is closed and bounded, it is compact, thus complete. For $r \in \mathcal{R}$, define $X^* = X^*(r)$ by 
\begin{equation}\label{Xstar}
X^* = -\frac{1}{2\alpha_0}\log r - \frac{\phi}{2\beta_0},
\end{equation}
so that
\begin{equation}\label{defr}
r = e^{-\alpha_0(2X^* + \phi/\beta_0)}.
\end{equation}
The constant $\phi$ comes from \cite[Lemma 6.1]{Sandstede1998} (see \cref{jumplemma3} below for details). We will use $r$ as a scaling parameter for the system. For $i = 0, \dots n-1$, define
\begin{equation}\label{bjscale}
b_i = e^{-2 \alpha_0 (X_i - X^*)},
\end{equation}
where we have chosen $X_i \geq X^*$ for $i = 0, \dots, n-1$. The quantities $b_i$ are length parameters for the system. In terms of $r$ and $b_i$,
\begin{equation}\label{Xiscale}
X_i = -\frac{1}{2\alpha_0}\log(b_i r) - \frac{\phi}{2 \beta_0}.
\end{equation}
In the next lemma, we rewrite the system \cref{jumpcond2} using this rescaling.

\begin{lemma}\label{jumplemma3}
A periodic multi-pulse solution exists if and only if for $i = 0, \dots, n-2$,
\begin{equation}\label{Geq}
G_i(b_1, \dots, b_{n-1}, r) = b_i \sin \left( -\rho \log b_i \right) - b_{n-1} \sin \left( -\rho \log b_{n-1} \right) + \mathcal{O}(r^{\gamma / 2 \alpha}) = 0,
\end{equation}
where $r \in \mathcal{R}$ and $0 < \gamma \leq 1$. All derivatives of the remainder term with respect to $b_i$ are also $\mathcal{O}(r^{\gamma / 2 \alpha})$. 
\begin{proof}
Using \cite[Lemma 6.1(i)]{Sandstede1998}, for $x > 0$ sufficiently large,
\begin{equation}\label{IPalphabeta}
\langle \Psi(-x), Q(x) \rangle
= p_0 e^{-2 \alpha_0 x} \sin(2 \beta_0 x + \phi) + \mathcal{O}\left(e^{-(2 \alpha_0 + \gamma) x}\right),
\end{equation}
where $0 < \gamma \leq 1$, $p_0 > 0$, and $\phi$ are constants which come from \cite[Lemma 6.1]{Sandstede1998} (note that we use $p_0$ in place of $s$ in that lemma, and that there is no dependence on a parameter $\mu$).
Substituting \cref{IPalphabeta} into \cref{jumpcond2} and rescaling using \cref{bjscale} and \cref{Xiscale},
\begin{equation}\label{diff2}
\begin{aligned}
p_0 e^{\alpha_0 \phi / \beta_0 } b_i r \sin \left( - \rho \log (b_i r) \right) 
- p_0 e^{\alpha_0 \phi / \beta_0 } b_{i-1} r \sin \left( -\rho \log (b_{i-1} r) \right) + \mathcal{O}(r^{1 + \gamma / 2 \alpha}) = 0.
\end{aligned}
\end{equation}
Dividing both sides by $r p_0 e^{\alpha_0 \phi / \beta_0 } > 0$ and simplifying, we obtain the jump conditions
\begin{align}\label{diff3}
\xi_i = b_i \sin \left( -\rho \log b_i \right) - b_{i-1} \sin \left( -\rho \log b_{i-1} \right) + \mathcal{O}(r^{\gamma / 2 \alpha}) &= 0,
\end{align} 
since $\sin \left( -\rho \log (b_i r) \right) = \sin \left( -\rho \log b_i \right)$ for $r \in \mathcal{R}$. For $i = 0, \dots, n-2$, let
\begin{equation}\label{Gidef}
G_i(b_1, \dots, b_{n-1}, r) = \sum_{k = 0}^i \xi_k.
\end{equation}
After canceling terms, we obtain the equations \cref{Geq}, which are equivalent to \cref{jumpcond2} via an invertible linear transformation.
\end{proof}
\end{lemma}

\begin{remark}
In \cref{jumplemma3}, we rewrote \cref{jumpcond2} so that equation $i$ involves $b_i$ and a common parameter $b_{n-1}$. Since we are on a periodic domain, that choice was arbitrary;  the final length parameter $b_{n-1}$ was chosen for notational convenience.
\end{remark}

When $r = 0$, the equations \cref{Geq} all have the same form. Let
\begin{equation}\label{defH}
H(b_0, b_1) = b_0 \sin \left( -\rho \log b_0 \right) - b_1 \sin \left( -\rho \log b_1 \right).
\end{equation}
In the next lemma, we will show that pitchfork bifurcations occur on the diagonal in the zero set of $H(b_0, b_1)$.

\begin{lemma}\label{pitchforkH}
A discrete family of pitchfork bifurcations occurs along the diagonal in the zero set of $H(b_0, b_1)$ at $(b_0, b_1) = (b_k^*, b_k^*)$ for $k \in \Z$, where $b^*_k = e^{-\frac{1}{\rho} (p^* + k \pi) }$ and $p^* = \arctan \rho$. Locally, the arms of the pitchfork bifurcations open upwards along the diagonal.
\begin{proof}
First, we note that the partial derivative
$H_{b_0}(b_0, b_1) = \sin \left( - \rho \log b_0 \right) - \rho \cos \left( - \rho \log b_0 \right)$
is zero if and only if $b_0 = b_k^*$ for integer $k$, which gives the locations of the bifurcation points. Next, we change coordinates so that the pitchfork bifurcation will occur along the horizontal axis. Let $b_0 = y-x$ and $b_1 = y+x$, which is a rotation by $-\pi/4$. Making this substitution, we obtain
\begin{equation}\label{Hxy}
H(x, y) = 
(y - x) \sin \left( -\rho \log(y - x) \right) - (y + x) \sin \left( - \rho \log (y + x) \right).
\end{equation}
For all $y$, $H(-x, y) = -H(x, y)$, which is the required odd symmetry for a pitchfork bifurcation. Let $(x_0, y_0) = \left(0, b^*_k \right)$. Evaluating the relevant partial derivatives of $H$ at $(x_0, y_0)$,
\begin{align*}
H_x(x_0, y_0) &= 0, \quad H_y(x_0, y_0) = 0, \quad H_{xx}(x_0, y_0) = 0, \quad H_{yy}(x_0, y_0) = 0 \\
H_{xy}(x_0, y_0) &= (-1)^k 2 \rho \sqrt{1 + \rho^2} \: \exp{\left(\frac{1}{\rho} (\arctan \rho - k \pi) \right)} \neq 0 \\
H_{xxx}(x_0, y_0)
&= -(-1)^k 2 \rho \sqrt{1 + \rho^2} \: \exp{\left(\frac{2}{\rho} (\arctan \rho - k \pi) \right)} \neq 0,
\end{align*}
thus a pitchfork bifurcation occurs at $(0, b^*_k)$ for all $k \in \Z$. To leading order, near the bifurcation points $(0, b_k^*)$, the arms of the pitchforks are upwards-opening parabolas of the form $y = b_k^* + c_k x^2$, where $c_k = \frac{1}{6}\exp{\left(\frac{1}{\rho} (\arctan \rho - k \pi) \right)} > 0
$. The result follows upon reverting to the $(b_0, b_1)$ coordinate system.
\end{proof}
\end{lemma}

Now that we have located the pitchfork bifurcations, we will construct a natural parameterization for the zero set of $H(b_0, b_1)$. We only need to consider $b_0 \geq b_1$ since the zero set is symmetric across the diagonal. For any nonnegative integers $m_1 \geq m_0$, the point
\[
(b_0, b_1) = \left( \exp\left( -\frac{1}{\rho}m_0 \pi \right), \exp\left( -\frac{1}{\rho}m_1 \pi \right) \right)
\]
is in the zero set of $H$. We will use these points to anchor our parameterization, and will use a phase parameter $\theta$ to connect these anchor points. 

\begin{lemma}\label{thetaparamlemma}
For any nonnegative integers $m_0$ and $m_1$ with $m_1 \geq m_0$, there is a smooth family of solutions
\begin{align*}
\left( b_0( m_0, m_1, \theta), b_1( m_0, m_1, \theta) \right) && \theta \in [-\pi + p^*, p^*]
\end{align*}
to $H(b_0, b_1) = 0$. This parameterization is given explicitly by
\begin{equation}\label{thetaparam}
\begin{aligned}
b_0( m_0, m_1, \theta) &= \exp\left( -\frac{1}{\rho}\left(m_0 \pi + \theta^*(\theta; m_1 - m_0) \right) \right) \\
b_1( m_0, m_1, \theta) &= \exp\left( -\frac{1}{\rho}\left(m_1 \pi + \theta \right) \right),
\end{aligned}
\end{equation}
where the functions $\theta^*(\theta; m) : [-\pi + p^*, p^*] \rightarrow \R$ are smooth in $\theta$ for all nonnegative integers $m$ and have the following properties:
\begin{enumerate}[(i)]
\item $\theta^*(\theta; 0) = \theta$.
\item $\theta^*(0; m) = 0 \text{ for all } m$.
\item $|\theta^*(\theta; m)| \leq |\theta|$.
\item $\theta^*(p^*; m) = \theta^*(-\pi+p^*; m+1)$.
\item $|\theta^*(\theta; m)| \leq C \exp\left(-\frac{m \pi}{\rho} \right)$.
\end{enumerate}
In particular, $\theta^*(p^*; 0) = \theta^*(\pi - p^*; 1) = p^*$.
\begin{proof}
We substitute \cref{thetaparam} into $H(b_0, b_1) = 0$ and solve for $\theta^*$ in terms of $\theta$.
\begin{enumerate}
\item First, we show that $\theta^*$ only depends on the difference $m_1 - m_0$. Substituting \cref{thetaparam} into $H(b_0, b_1) = 0$ and simplifying, we obtain the equation
\begin{equation}\label{thetaeq1}
e^{-\frac{1}{\rho}\theta^*}\sin \theta^* = e^{-\frac{1}{\rho}(m_1 - m_0) \pi} (-1)^{m_1 - m_0} e^{-\frac{1}{\rho}\theta}\sin \theta = 0,
\end{equation}
which only depends on $m_1 - m_0$. Letting $m = m_1 - m_0$, it suffices to solve
\begin{equation}\label{thetaeq2}
g(\theta^*) = t(m) g(\theta)
\end{equation}
for all nonnegative integers $m$, where $g(\theta) = e^{-\frac{1}{\rho}\theta}\sin \theta$ and $t(m) = (-1)^m e^{-\frac{1}{\rho}m \pi}$. 

\item For $m = 0$, $t(0) = 1$, thus $\theta^*(\theta; 0) = \theta$.

\item Let $I = \left[-\pi + p^*, p^*\right]$. For $m \geq 1$, we first show that $g(\theta)$ is invertible on $I$. Since
\begin{equation}\label{gprime}
g'(\theta) = e^{ -\frac{1}{\rho} \theta } \left( \cos \theta - \frac{1}{\rho} \sin \theta \right),
\end{equation}
$g'(p^*) = 0$, $g'(-\pi+p^*) = 0$, and the only critical point of $g'(\theta)$ on $I$ is a local maximum at $\theta = - \pi + 2 p^*$, $g(\theta)$ is strictly increasing, thus invertible, on $I$. Let
\begin{equation*}
g(I) = \left[-e^{\frac{1}{\rho}\pi} T, T\right], \quad T = e^{-\frac{1}{\rho}p^*} \sin p^* = \frac{\rho}{\sqrt{1+\rho^2}}e^{-\frac{1}{\rho}\arctan \rho}.
\end{equation*}
Then $g: I \rightarrow g(I)$ is a bijection, and $g^{-1}: g(I) \rightarrow I$ is also strictly increasing. Since the only zero of $g(\theta)$ on $I$ occurs at $\theta = 0$, $g^{-1}(0) = 0$. We can now solve \cref{thetaeq2} for $\theta^*$. For all $\theta \in I$,
\begin{equation}\label{RHSbounds}
	t(m)g(\theta) \in
	\begin{cases}
	[-e^{-\frac{1}{\rho}(m-1) \pi} T, e^{-\frac{1}{\rho}m \pi} T] & m \text{ even }\\
	[-e^{-\frac{1}{\rho}m \pi} T, e^{-\frac{1}{\rho}(m-1) \pi} T] & m \text{ odd.}
	\end{cases}
\end{equation}
Since $t(m)g(\theta) \subset g(I)$ for all $\theta \in I$, define
\begin{align}\label{solvethetastar}
\theta^*(\theta; m) &= g^{-1}\left( t(m) g(\theta) \right) && \theta \in I, m \geq 1.
\end{align}
Since $g(0) = 0$, $\theta^*(0, m) = 0$ for all $m$.

\item Next, we show that $|\theta^*(\theta; m)| \leq |\theta|$. For $m = 0$, we have equality. For $m = 1$ and $\theta = -\pi + p^*$, $t(-\pi + p^*)g(-\pi + p^*) = T$, thus $\theta^*(-\pi + p^*; 1) = -\pi + p^*$. For any other $m$ and $\theta$, it follows from \cref{RHSbounds} that $t(m)g(\theta) \in [e^{-\frac{1}{\rho}\pi}T, T)$, which is strictly contained in $g(I)$. Since $g^{-1}$ is strictly increasing and $g(0) = 0$, $|\theta^*(\theta; m) \leq |\theta|$.

\item We now show that consecutive parameterizations match up at the endpoints. For $m \geq 0$, 
\begin{align*}
\theta^*(p^*; m) &= g^{-1}\left( e^{-\frac{1}{\rho}m \pi} (-1)^m e^{-\frac{1}{\rho}p^*} \sin p^* \right)
\end{align*}
and
\begin{align*}
\theta^*(-\pi + p^*; m+1) &= g^{-1}\left( e^{-\frac{1}{\rho}(m+1) \pi} (-1)^{m+1} e^{-\frac{1}{\rho}(-\pi + p^*)} \sin (-\pi + p^*) \right) \\
&=g^{-1}\left( e^{-\frac{1}{\rho}m \pi} (-1)^m e^{-\frac{1}{\rho}p^*} \sin p^* \right),
\end{align*}
which are equal. In particular, $\theta^*(p^*; 0) = \theta^*(\pi - p^*; 1) = p^*$.

\item Finally, we obtain a bound on $\theta^*(\theta; m)$. For $m \geq 2$ and $\theta \in I$, it follows from \cref{RHSbounds} that $t(m)g(\theta)$ is contained in an interval $\tilde{I}$, which is strictly contained in $I$. Since $g'(0) = 0$ only at the endpoints of $I$ and is positive in the interior of $I$, $g'(\theta)$ is bounded below on $\tilde{I}$, thus $[g^{-1}]'(\theta)$ is bounded for $\theta \in g(\tilde{I})$. Since $|t(m)g(\theta)| \leq e^{-\frac{1}{\rho}(m - 1)\pi}T$ by \cref{RHSbounds} and $g(0) = 0$, we obtain the bound (v), which is independent of $\theta$.
\end{enumerate}
\end{proof}
\end{lemma}

\subsection{Proof of Theorem \ref{th:perexist}}

By \cref{jumplemma3}, a periodic $n$-pulse exists if and only if $G_i(b_0, \dots, b_{n-1}, r) = 0$ for $i = 0, \dots, n-2$. When $r = 0$, $G_i(b_0, \dots, b_{n-1}, 0) = H(b_i, b_{n-1})$, and we can use the parameterization from \cref{thetaparamlemma}. Choose any periodic parameterization $(m_0, \dots, m_{n-1}, \theta)$ and define
\begin{equation}\label{thetaparammulti}
\begin{aligned}
b_i( m_i, m_{n-1}, \theta) &= \exp\left( -\frac{1}{\rho}\left(m_i \pi + \theta^*(\theta; m_{n-1} - m_i) \right) \right) \qquad\qquad i = 0, \dots, n-2 \\
b_{n-1}( m_{n-1}, \theta) &= \exp\left( -\frac{1}{\rho}\left(m_{n-1} \pi + \theta \right) \right),
\end{aligned}
\end{equation}
so that by \cref{thetaparamlemma},
\begin{align*}
H(b_i( m_i, m_{n-1}, \theta), b_{n-1}( m_{n-1}, \theta) ) &= 0 && i = 0, \dots, n-2.
\end{align*}
Substituting $b_{n-1}(m_{n-1}, \theta)$ for $b_{n-1}$ in \cref{Geq}, define 
$b = (b_0, \dots, b_{n-1}) \in \R^{n-1}$ and 
$G: \R^{n-1} \times \mathcal{R} \rightarrow \R^{n-1}$ by $G = (G_0, \dots, G_{n-2})^T$, where 
\begin{equation*}
G_i(b, r) = b_i \sin \left( -\rho \log b_i \right) - b_{n-1}(m_{n-1}, \theta) \sin \left( -\rho \log b_{n-1}(m_{n-1}, \theta) \right) + \mathcal{O}(r^{\gamma / 2 \alpha}).
\end{equation*}
Let $b^* = (b_0( m_0, m_{n-1}, \theta), \dots, b_{n-2}( m_{n-2}, m_{n-1}, \theta) )$. Then $G(b^*, 0) = 0$, and the Jacobian matrix $D_b G(b^*,0)$ is diagonal, with
\begin{align*}
\partial_{b_i} G_i(b^*, 0)
&= \sin \left( -\rho \log b_i(m_i, \theta) \right) - \rho \cos \left( -\rho \log b_i(m_i, \theta) \right).
\end{align*}
From \cref{pitchforkH}, $\partial_{b_i} G_i(b, 0) = 0$ if and only if $b_i$ is one of the pitchfork bifurcation points $b_k^*$. By \cref{thetaparamlemma}, if we exclude $\theta \in \{ -\pi + p^*, p^* \}$, the Jacobian $D_b G(b^*,0)$ is invertible, thus we can use the implicit function theorem to solve for $b$ in terms of $r$ near $b^*$. Specifically, there exists $r_* > 0$ and a unique continuous function $b: \mathcal{R} \cap [0, r_*] \rightarrow \R^{n-1}$ given by $b(r) = \left( b_0(r), \dots, b_{n-2}(r) \right)$ such that $b(0) = b^*$ and $G(b,r) = 0$ if and only if $b = b(r)$. For $i = 0, \dots, n-2$, define $t_i(r; m_i, \theta): \mathcal{R} \rightarrow \R$ by
\begin{equation}\label{defti}
t_i(r; m_i, m_{n-1}, \theta) = -\rho \log b_i(r),
\end{equation}
which is continuous in $r$, and $t_i(0; m_i, m_{n-1}, \theta) = m_i \pi + \theta^*(\theta; m_{n-1} - m_i)$ by \cref{thetaparammulti}. Equations \cref{Xi} are obtained by substituting $b_i(r)$, $i = 0, \dots, n-2$, and $b_{n-1}(m_{n-1}, \theta)$ for $b_i$ in \cref{Xiscale} and using \cref{defti}.

\subsection{Periodic 2-pulse}

For the periodic 2-pulse, we have a single jump condition
\begin{equation}\label{2pulsedefG}
G(b_0, b_1, r) = b_0 \sin \left( -\rho \log b_0 \right) - b_1 \sin \left( -\rho \log b_1 \right) + \mathcal{O}(r^{\gamma / 2 \alpha}) = 0.
\end{equation}
First, we will show that the pitchfork bifurcations along the diagonal persist for small $r$. Recall that by \cref{def:perparam}, $m_0 \in \{0, 1\}$ for a periodic 2-pulse.

\begin{lemma}\label{pitchpersist}
There exists $r_2 > 0$ such that for $m_0 \in \{0, 1\}$ and $r \leq r_2$, there is a non-degenerate pitchfork bifurcation in the zero set of $G(b_0, b_1, r)$ at $(b_{m_0}^*(r),b_{m_0}^*(r))$, and 
\begin{equation*}
b_{m_0}^*(r) \rightarrow b_{m_0}^* \text{ as } r \rightarrow 0.
\end{equation*}
\begin{proof}
Take $m_0 = 0$. The proof is identical for $m_0 = 1$. First, we show the required odd symmetry relation. By \cref{solvewithjumps}, for an ordered pair of pulse distances $(X_0, X_1)$ with $X_i$ sufficiently large, there exists a unique piecewise solution $(U_0^-(x), U_0^+(x), U_1^-(x), U_1^+(x))$ which is continuous except for two jumps 
\begin{equation}\label{xijumps}
\xi_0(X_0, X_1) = \langle \Psi(0), U_0^+(0) - U_0^-(0) \rangle, \quad
\xi_1(X_0, X_1) = \langle \Psi(0), U_1^+(0) - U_1^-(0) \rangle 
\end{equation}
in the direction of $\Psi(0)$. By symmetry, $(U_1^+(-x), U_1^-(-x), U_0^+(-x), U_0^-(-x))$ is also a solution for $(X_0, X_1)$, thus it must be the same solution by uniqueness. In particular, $U_1^+(0) = U_0^-(0)$ and $U_1^-(0) = U_0^+(0)$, thus $\xi_0(X_0, X_1) = -\xi_1(X_0, X_1) = -\xi_0(X_1, X_0)$. Since swapping $X_0$ and $X_1$ is equivalent to swapping $b_0$ and $b_1$ in \cref{2pulsedefG}, $G(b_0, b_1, r) = -G(b_1, b_0, r)$ for sufficiently small $r$. 

Making the change of coordinates $(b_0, b_1) \mapsto (x, y)$ as in \cref{pitchforkH}, $G(-x, y, r) = -G(x, y, r)$ for sufficiently small $r$. The persistence of the pitchfork bifurcation follows from a Lyapunov-Schmidt reduction.  By \cref{pitchforkH}, $G_x(0, b_0^*, 0) = 0$ and $G_{xy}(0, b_0^*, 0) \neq 0$, thus by the implicit function theorem there exists $r_2 > 0$, an open interval $(-a, a)$, and a unique smooth function $y = y^*(x, r)$ such that $y^*(0, 0) = b_0^*$ and $G_x(x, y^*(x, r), r) = 0$ for all $x \in (-a, a)$ and $r < r_2$. It follows that a pitchfork bifurcation occurs at $(x, y, r) = (0, y^*(0, r), r)$ by evaluating the appropriate partial derivatives of $G$ as in \cref{pitchforkH}. Letting $b_0^*(r) = y^*(0, r)$, the result follows upon reverting to the original $(b_0, b_1)$ coordinates.
\end{proof} 
\end{lemma}

Next, we show that the arms of the pitchfork persist for sufficiently small $r$. By symmetry, it suffices to show this for the lower arm.

\begin{lemma}\label{armpersists}
Choose any $\delta > 0$. Then there exists $r_3 > 0$ such that for $m_0 \in \{0, 1\}$ and $r \leq r_3$, the portion of the zero set of $G(b_0, b_1, r)$ corresponding to lower arm of the pitchfork at $(b_{m_0}^*(r), b_{m_0}^*(r))$ is parameterized by
\begin{align*}
(b_0, b_1) = (b_0(s; m_0, r), b_1(s)) && s \in [p^* + \delta, \infty),
\end{align*}
where
\begin{align}\label{defb1}
b_1(s) &= e^{-\frac{1}{\rho}s}.
\end{align}

\begin{proof}
As in the previous lemma, take $m_0 = 0$. The proof is identical for $m_0 = 1$. For every positive integer $m_1$ with $m_1 \geq 1$, let
\begin{equation}\label{thetaparam2}
\begin{aligned}
\tilde{b}_0( m_1, \theta) = \exp\left( -\frac{1}{\rho}\theta^*(\theta; m_1) \right), \quad
\tilde{b}_1( m_1, \theta) = \exp\left( -\frac{1}{\rho}(m_1 \pi + \theta) \right),
\end{aligned}
\end{equation}
where $\theta^*$ is defined in \cref{thetaparamlemma}. Since these families connect at their endpoints, define $b_1(s)$ by \cref{defb1}. For $s > p^*$, let $m(s) = \left\lceil \frac{s - p^*}{\pi} \right\rceil$ and $\theta(s) = s - m(s) \pi$, and define $\tilde{b}_0(s)$ by
\begin{equation}\label{defb0}
\tilde{b}_0(s) = b_0(m(s), \theta(s)),
\end{equation}
so that the continuous curve $(\tilde{b}_0(s), b_1(s))$ for $s > p^*$ parameterizes the lower arm of the pitchfork when $r = 0$. Define the Banach space $X = C_b([p^* + \delta, \infty), \R)$ of bounded continuous functions equipped with the uniform norm. For $b \in X$, define $\tilde{G}: X \times \mathcal{R} \rightarrow X$ by
\begin{align*}
\left[\tilde{G}(b, r)\right](s) &= G(b(s), b_1(s), r) = b(s) \sin(-\log b(s) ) - e^{-\frac{1}{\rho}s} \sin s + \mathcal{O}(r^{\gamma/2\alpha}).
\end{align*}
Then $\tilde{G}(\tilde{b}_0(s), r) = 0$, and $D_{b_0} \tilde{G}(b, 0) = G_{b_0}(b(s), b_1(s), 0)$.

By \cref{pitchforkH}, $G_{b_0}(b_0, b_1, 0) = 0$ if and only if $(b_0, b_1)$ is one of the pitchfork bifurcation points; by \cref{thetaparamlemma}, these occur on the curve $(\tilde{b}_0(s), b_1(s))$ only when $s = p^*$. From the proof of \cref{thetaparamlemma}, $|G_{b_0}(\tilde{b}_0(s), b_1(s), 0)|$ is bounded below for $s \geq p^* + \delta$, thus $D_{b_0} \tilde{G}(\tilde{b}_0, 0)$ is invertible with bounded inverse. Using the implicit function theorem for Banach spaces, there exists $r_3 > 0$ and a unique smooth function $b: \mathcal{R} \rightarrow X$ with $b(0) = \tilde{b}_0$ such that for all $r \leq r_3$, $\tilde{G}(b(r), r) = 0$. It follows from the definition of $\tilde{G}$ that $G([b(r)](s), b_1(s), r) = 0$ for all $r \leq r_3$ and $s \in [p^* + \delta, \infty)$. The result follows by taking $b_0(s; m_0, r) = [b(r)](s)$.
\end{proof}
\end{lemma} 

Finally, we show that for sufficiently small $r$, the lower arm of the pitchfork connects to the pitchfork bifurcation point, which extends the parameterization in \cref{armpersists} to $s \in [p^*, \infty)$.

\begin{lemma}\label{pitchforkconnects}
There exists $r_4 > 0$ such that for $m_0 \in \{0, 1\}$ and $r \leq r_4$, the parameterization in \cref{armpersists} can be extended to $s \in [p^*, \infty)$, and 
\begin{align*}
\left(b_0(p^*; m_0, r), b_1(p^*)\right) = \left(b_{m_0}^*(r), b_{m_0}^*(r)\right),
\end{align*}
which is the pitchfork bifurcation point from \cref{pitchpersist}.

\begin{proof}
For simplicity, take $m_0 = 0$. The proof is identical for $m_0 = 1$. Change variables $(b_0, b_1) \mapsto (x, y)$ as in \cref{pitchforkH}, so that the pitchfork bifurcation takes place on the horizontal axis. Let $r_2$ be as in \cref{pitchpersist}. Then there is a nondegenerate pitchfork bifurcation at $(x, y) = (b_0^*(r), 0)$, and there exists $y_1 \geq 0$ such that for $r \leq r_2$, the arm of the pitchfork is uniquely parameterized by $(x, y) = (x_0(y, r), y)$ for $y \in [0, y_1]$, where $x_0(0, r) = b_0^*(r)$. Take $\delta = y_1/2$, and let $r_3$ be as in \cref{armpersists}. In the $(x, y)$ coordinate system, the lower arm of the pitchfork is uniquely parameterized by 
$(x, y) = (x_1(y, r), y)$ for  $y \in [y_1/2, \infty)$ for $r < r_3$. Let $r_4 = \min\{ r_2, r_3 \}$. Then by the uniqueness of the two parameterizations, $(x_1(y, r), y) = (x_0(y, r), y)$ for $y \in [y_1/2, y_1]$ and $r \leq r_4$. Since the two parameterizations overlap on an interval, the lower arm of the pitchfork connects to the pitchfork bifurcation point. Returning to the $(b_0, b_1)$ coordinate system, we can extend the parameterization in \cref{armpersists} to the pitchfork bifurcation point, which occurs when $s = p^*$.
\end{proof}
\end{lemma}

\subsection{Proof of Theorem \ref{2pulsebifurcation}}

Let $r_* = r_4$, where $r_4$ is defined in \cref{pitchforkconnects}. From the proof of \cref{pitchpersist}, $G(b_0, b_0, r) = -G(b_0, b_0, r)$, thus symmetric solutions with $b_0 = b_1$ exist for sufficiently small $r$. To parameterize these, for $m_0 \in \{0, 1\}$ and $s_0 \in [0, \pi)$, let
\[
b_0(m_0, s_0) = b_1(m_0, s_0) = 
\exp\left( -\frac{1}{\rho}(m_0 \pi + s_0) \right).
\]
The pulse distances \cref{2psymmdist} are obtained by substituting this into \cref{Xiscale}. Let $p^*(m_0; r) = -\rho \log(b_{m_0}^*(r))$, where $b_{m_0}^*(r)$ is the pitchfork bifurcation point defined in \cref{pitchpersist}. Then the pitchfork bifurcation occurs when $s_0 = p^*(m_0; r)$, and $p^*(m_0; r) \rightarrow p^*$ as $r \rightarrow 0$.

For asymmetric periodic 2-pulses, taking $s_1 = s$ in \cref{pitchforkconnects}, the lower arms of the pitchforks are parameterized by $(b_0, b_1) = (b_0(s_1; m_0, r), b_1(s_1))$ for $s_1 \in [p^*, \infty)$. The formula for $X_1(r, s_1)$ in \cref{2pasymmdist} follows by substituting \cref{defb1} into \cref{Xiscale}. Let $t_0(r; m_0, s_1) = -\rho \log\left( b_0(s_1; m_0, r) \right)$, which is continuous in $r$ and $s_1$. Using this together with \cref{Xiscale} we obtain the formula for $X_0(r; m_0, s_1)$ in \cref{2pasymmdist}. From \cref{armpersists}, $t_0(0; m_0, s_1) = m_0 \pi + \theta^*(\theta(s_1); m(s_1) - m_0)$, where $m(s_1) = \lceil \frac{s_1 - p^*}{\pi} \rceil$ and $\theta(s_1) = s_1 - m(s_1) \pi$. Using the estimate for $\theta^*(\theta; m)$ from \cref{thetaparamlemma},
\[
\theta^*(\theta(s_1); m(s_1) - m_0) \leq C \exp\left(-\frac{1}{\rho} ( m(s_1) - m_0 )\pi \right) \leq C \exp\left(-\frac{1}{\rho} s_1 \right),
\]
from which the estimate \cref{t0est} follows. By \cref{pitchforkconnects}, the pitchfork bifurcation point is reached when $s_1 = p^*$.

\section{Proof of Theorem \ref{blockmatrixtheorem}}\label{sec:blockmatrixproof}

We will use Lin's method as in \cite{Sandstede1998} to construct eigenfunctions which are solutions to \cref{PDEeigsystemper3}. To do this, we will take a piecewise linear combination of the kernel eigenfunctions $\partial_x Q_n(x)$ and $\partial_c Q_n(x)$ and a ``center'' eigenfunction as our ansatz, and we will join these together using small remainder functions. As long as the individual pulses in $Q_n(x)$ are well-separated, Lin's method will yield a unique solution which solves \cref{PDEeigsystemper3} but which has $n$ discontinuities. In contrast to \cite{Sandstede1998}, these $n$ jumps line in the two-dimensional subspace spanned by $\Psi(0)$ and $W_0$, which gives us $2n$ jump conditions. Finding the eigenvalues near 0 amounts to solving these jump conditions, which will give us both the interaction eigenvalues and the essential spectrum eigenvalues.

\subsection{Preliminaries}

For convenience, we define
\begin{equation*}
\begin{aligned}
A(Q(x); \lambda) &= A(Q(x)) + \lambda B,
\end{aligned}
\end{equation*}
and note that $A(0; \lambda) = A(0) + \lambda B$ is a constant matrix.
It follows from \cref{defAB} and the symmetry relations \cref{genODErev} that 
\begin{equation}\label{AQsymmetry}
A(Q(x); \lambda) = -R A(Q(-x); -\lambda)R,
\end{equation}
where $R$ is standard reversor operator. Let $\alpha_0$ and $\beta_0$ be defined as in \cref{hyp:hypeq}. Choose any $\eta > 0$ with $2 \eta < \alpha_0$, and let $\alpha = \alpha_0 - \eta$. Let $\delta_1$ be as in \cref{nulambdalemma}, and choose $\delta_2 \leq \delta_1$ sufficiently small so that for all $|\lambda| < \delta_2$, $|\nu(\lambda)| < \eta$, where $\nu(\lambda)$ is the simple eigenvalue of $A(0; \lambda)$ close to 0 which is defined in \cref{nulambdalemma}, and $|\Re \nu| > \alpha$ for any other eigenvalue $\nu$ of $A(0; \lambda)$. To greatly simplify our analysis, we place the additional assumption on the real part of $\lambda$
\begin{equation}\label{relambdabound}
|\Re \lambda| \leq r^{1/4} = C e^{-\frac{1}{2}\alpha X^*},
\end{equation}
where the scaling parameter $r$ is defined in \cref{defr} and $X^*$ is defined in \cref{Xstar}. We will verify that this assumption is satisfied for sufficiently small $r$ when we consider applications of the theorem. It then follows from \cref{relambdabound} and \cref{nulambdalemma} that 
\begin{equation}
|e^{\nu(\lambda)X^*}| \leq e^{|\Re \nu(\lambda)| X^*} \leq \exp \left( C X^* e^{-\frac{1}{2} \alpha_0 X^* } \right) \leq C.
\end{equation}
Since the periodic parameterization $(m_0, \dots, m_{n-1}, \theta)$ is fixed,
\begin{align}\label{expnubound}
|e^{\nu(\lambda)X_i}| \leq C  && i = 0, \dots, n-1.
\end{align}

\subsection{Conjugation lemma}

The conjugation lemma allows us to make a smooth change of coordinates to convert certain linear ODEs of the form $W'(x) = A(x) W(x)$ into a constant coefficient system. The statement of the lemma is identical to that in \cite{Zumbrun2009}, except that the parameter vector $\Lambda$ here lives in an arbitrary Banach space. The proofs of \cref{conjlemma} and \cref{corr:adjconj} are straightforward modifications of the proof of \cite[Corollary 2.3]{Zumbrun2009}.

\begin{lemma}[Conjugation lemma]\label{conjlemma}
Let $V \in \C^N$, and consider the family of ODEs on $\R$
\begin{equation}\label{EVPconj}
V(x)' = A(x; \Lambda) V(x) + F(x),
\end{equation}
where $\Lambda \in \Omega$ is a parameter vector and $\Omega$ is a Banach space. Assume that
\begin{enumerate}[(i)]
	\item The map $\Lambda \mapsto A(\cdot; \Lambda)$ is analytic in $\Lambda$.
	\item $A(x; \Lambda) \rightarrow A^\pm(\lambda)$ (independent of $\Lambda$) as $x \rightarrow \pm \infty$, and there exists $\delta > 0$ such that for $|\Lambda| < \delta$, we have the uniform exponential decay estimates 
	\begin{align}
	\left| \frac{\partial^k}{\partial x^k} A(x; \Lambda) - A^\pm(\Lambda) \right| 
	&\leq C e^{-\theta |x|} && 0 \leq k \leq K,
	\end{align}
	where $\alpha > 0$, $C > 0$, and $K$ is a nonnegative integer.
\end{enumerate}
Then in a neighborhood of any $\Lambda_0 \in \Omega$ there exist invertible linear transformations
\begin{equation}\label{conjlemmaP}
\begin{aligned}
P^+(x, \Lambda) &= I + \Theta^+(x, \Lambda) \\
P^-(x, \Lambda) &= I + \Theta^-(x, \Lambda) 
\end{aligned}
\end{equation}
defined on $\R^+$ and $\R^-$, respectively, such that
\begin{enumerate}[(i)]
\item The change of coordinates $V = P^\pm Z$ reduces \cref{EVPconj} to the equations on $\R^\pm$
\begin{align}\label{conjZ}
Z'(x) = A^\pm(\Lambda) Z(x) + P^\pm(x, \Lambda)^{-1} F(x).
\end{align}

\item For any fixed $0 < \tilde{\theta} < \theta$, $0 \leq k \leq K+1$, and $j \geq 0$ we have the decay rates
\begin{align}\label{conjthetadecay}
\left| \partial_\Lambda^j \partial_x^k \Theta^\pm \right| \leq C(j, k)e^{-\tilde{\theta}|x|}.
\end{align}
\end{enumerate}
\end{lemma}

\begin{corollary}\label{corr:adjconj}
Take the same hypotheses as in \cref{conjlemma}, and let $P^\pm(x; \Lambda)$ be the conjugation operators for $V(x)' = A(x; \Lambda) V(x)$ on $\R^\pm$. Then the change of coordinates $W = [(P^\pm)^{-1}]^* Z$ on $\R^\pm$ reduces the adjoint equation $W'(x) = -A(x; \Lambda)^* W(x)$ to the equation $Z'(x) = -A^\pm(\Lambda)^* Z(x)$.
\end{corollary}

\subsection{Solutions in center subspace}\label{sec:conjvareq}

First, we apply the conjugation lemma to
\begin{equation}\label{vareqlambda}
V'(x) = A(Q(x); \lambda) V'(x).
\end{equation}
For all $\lambda$, $A(Q(x); \lambda)$ decays exponentially to the constant-coefficient matrix $A(0; \lambda)$. Since $DF(0)$ is hyperbolic, $|A(Q(x); \lambda) - A(0; \lambda)| \leq C e^{-(\alpha_0 + \epsilon) |x|}$ for small $\epsilon$; the price to pay is a larger constant $C$. Using the conjugation lemma on $\R^+$ with $\Lambda = \lambda$ and $\Lambda_0 = 0$, there exists $\delta_3 \leq \delta_2$ and an invertible linear transformation 
\begin{equation}\label{defPplus}
P^+(x; \lambda) = I + \Theta^+(x; \Lambda),
\end{equation}
such that for all $|\lambda| < \delta_3$, the change of coordinates $V(x) = P^+(x; \lambda) Z^+(x)$ conjugates \cref{vareqlambda} into the constant-coefficient equation $(Z^+)'(x) = A(0; \lambda) Z^+(x)$. The function $\Theta^+(x; \lambda)$ has the uniform decay rate
\begin{equation}\label{Thetadecay}
|\Theta^+(x; \lambda)| \leq C e^{-\alpha_0 |x|},
\end{equation}
which holds for derivatives with respect to $x$ and $\lambda$. 

For $x \in \R^-$ and $|\lambda| < \delta_3$, define $P^-(x; \lambda)$ by
\begin{equation}\label{defPminus}
P^-(x; \lambda) = RP^+(-x; -\lambda)R.
\end{equation}
By a straightforward adaptation of the proof of the conjugation lemma in \cite{Zumbrun2009} and the symmetry relation \cref{AQsymmetry}, the change of coordinates $V(x) = P^-(x; \lambda) Z^-(x)$ conjugates \cref{vareqlambda} into the constant-coefficient equation $(Z^-)'(x) = A(0; \lambda) Z^-(x)$ on $\R^-$.

Let $E^{u/s/c}(0)$ be the stable, unstable, and center eigenspaces of $A(0)$, and $P^{u/s/c}(0)$ be their respective eigenprojections. Let $E^{u/s/c}(\lambda)$ and $P^{u/s/c}(\lambda)$ be the corresponding eigenspaces and eigenprojections for $A(0; \lambda)$, which are smooth in $\lambda$. $E^s(\lambda)$ and $E^u(\lambda)$ are $m$-dimensional, and $E^c(\lambda)$ is 1-dimensional. In the next lemma, we collect some useful results about $A(0; \lambda)$.

\begin{lemma}\label{lemma:Afacts}
We have the following results concerning $A(0; \lambda)$.
\begin{enumerate}[(i)]
	\item $A(0; -\lambda) = -R A(0; \lambda)R$, where $R$ is the standard reversor operator. In particular, $A(0) = -R A(0)R$.
	\item If $V$ is an eigenvector of $A(0)$ corresponding to eigenvalue $\mu$, then $RV$ is an eigenvector of $A(0)$ corresponding to eigenvalue $-\mu$, and $\overline{V}$ is an eigenvector of $A(0)$ corresponding to eigenvalue $\overline{\mu}$. 
	\item Let $P^{c,*}(0)$ be the center eigenprojection of $-A(0)^*$. Then $P^{c,*}(0) = [P^c(0)]^*$. 
	\item Let $P^{s,*}(0)$ and $P^{u,*}(0)$ be the stable and unstable eigenprojections of $-A(0)^*$. Then $P^{s,*}(0) = [P^u(0)]^*$ and $P^{u,*}(0) = [P^s(0)]^*$
\end{enumerate}	
\begin{proof}
Part (i) can be verified by multiplying out $R A(0; \lambda)R$. For part (ii), $A(0)V = \mu V$ implies $A(0)\overline{V} = \overline{\mu}\overline{V}$ since $A(0)$ is real. Using part (i), $-R A(0) R V = \mu V$, which rearranges to $A(0) (RV) = -\mu(RV)$. For part (iii), $\ker [P^c(0)]^* = (\ran P^c(0))^\perp = \spn\{V_0\}^\perp$  and $\ran [P^c(0)]^* = (\ker P^c(0))^\perp$. For any eigenvector $W$ of $-A(0)^*$ with nonzero eigenvalue $\mu$,
\begin{align*}
\langle W, V_0 \rangle = \frac{1}{\overline{\mu}} \langle \mu W, V_0 \rangle
= \frac{1}{\overline{\mu}}\langle -A(0)^* W, V_0 \rangle = -\frac{1}{\overline{\mu}}\langle W, A(0) V_0 \rangle = 0,
\end{align*}
thus $\ker [P^c(0)]^*$ is the direct sum of the stable and unstable eigenspaces of $-A(0)^*$. Similarly, for any eigenvector $V$ of $A(0)$ with nonzero eigenvalue $\mu$, $\langle W_0, V \rangle = 0$, thus $\ran [P^c(0)]^* = \spn\{W_0\}$, which proves (iii).

For part (iv), let $W$ be an eigenvector of $-A(0)^*$ with nonzero eigenvalue $\mu$. Then for all eigenvectors of $A(0)$ with eigenvalue $\tilde{\mu}$, we have 
\begin{align*}
\langle W, V \rangle = \frac{1}{\overline{\mu}} \langle \mu W, V \rangle
= \frac{1}{\overline{\mu}}\langle -A(0)^* W, V \rangle = -\frac{1}{\overline{\mu}}\langle W, A(0) V \rangle = -\frac{1}{\overline{\mu}}\langle W, \tilde{\mu} V_0 \rangle = -\frac{\tilde{\mu}}{\overline{\mu}}\langle W, V \rangle,
\end{align*}
thus $\langle W, V \rangle = 0$ unless $\tilde{\mu} = -\overline{\mu}$. For the stable eigenprojection, $\ker [P^s(0)]^* = (\ran P^s(0))^\perp$ and $\ran [P^s(0)]^* = (\ker P^s(0))^\perp = (E^u(0) \oplus E^c(0))^\perp$. Let $W$ be in the unstable eigenspace of $-A(0)^*$, and let $\mu$ be the corresponding eigenvalue with $\Re \mu > 0$. Then $W$ is perpendicular to all eigenvectors in $E^c(0) \oplus E^u(0)$, thus $W \in \ran [P^s(0)]^*$. Similarly, if $W$ is in the unstable or center eigenspace of $-A(0)^*$, then $W \in \ker [P^s(0)]^*$. This proves (iv) for $P^{u,*}(0)$. The result for $P^{s,*}(0)$ is similar.
\end{proof}
\end{lemma}

Let $\Phi(x, y; \lambda) = e^{A(0; \lambda)(x-y)}$ be the evolution operator of the constant-coefficient equation $Z'(x) = A(0; \lambda) Z(x)$. Then the evolution operator of the unconjugated equation \cref{vareqlambda} on $\R^\pm$ is given by
\begin{align}\label{deftildephi}
\tilde{\Phi}^\pm(x, y; \lambda) &= P^\pm(x; \lambda) \Phi(x, y; \lambda) P^+(y; \lambda)^{-1} && x, y \in \R^\pm.
\end{align}
We can decompose \cref{deftildephi} in exponential trichotomies on $\R^\pm$ via the operators
\begin{align}\label{tildephitrich}
\tilde{\Phi}^{s/u/c, \pm}(x, y; \lambda) &= P^\pm(x; \lambda) \Phi(x, y; \lambda) P^{s/u/c}(\lambda) P^\pm(x; \lambda)^{-1} && x, y \in \R^\pm,
\end{align}
where we have estimates
\begin{equation}\label{stdtrichbounds}
\begin{aligned}
|\tilde{\Phi}^{s,\pm}(x, y; \lambda)| &\leq C e^{-\alpha(x - y)} \qquad\qquad y \leq x \\
|\tilde{\Phi}^{u,\pm}(x, y; \lambda)| &\leq C e^{-\alpha(y - x)} \qquad\qquad x \leq y \\
|\tilde{\Phi}^{c,\pm}(x, y; \lambda)| &\leq C e^{\eta|x - y|},
\end{aligned}
\end{equation}
which are uniform in $\lambda$. 

The equation $Z'(x) = A(0; \lambda) Z(x)$ has a solution $Z(x) = V_0(\lambda)e^{\nu(\lambda)x}$ which lies in the one-dimensional eigenspace $E^c(\lambda)$ spanned by $V_0(\lambda)$. In the next lemma, we show that equation \cref{vareqlambda} has solutions $V^\pm(x; \lambda)$ on $\R^\pm$ which approach $V_0(\lambda)e^{\nu(\lambda)x}$ as $x \rightarrow \pm \infty$. 

\begin{lemma}\label{lemma:Vpm}
For sufficiently small $|\lambda|$, equation \cref{vareqlambda} has solutions
\begin{align}\label{Vpmlambda}
V^\pm(x; \lambda) &= e^{\nu(\lambda)x}(V_0(\lambda) + V_1^\pm(x; \lambda)) && x \in \R^\pm,
\end{align}
where $|V_1^\pm(x; \lambda)| \leq C e^{-\alpha_0 |x|}$ and $V^-(x; \lambda) = R V^+(-x; -\lambda)$.
\begin{proof}
Let $Z(x) = e^{\nu(\lambda)x}V_0(\lambda)$, and define
\begin{equation}\label{defVplus}
V^+(x; \lambda) = P^+(x; \lambda) Z(x) = e^{\nu(\lambda)x}P^+(x; \lambda)V_0(\lambda).
\end{equation}
By \cref{conjlemmaP}, $V^+(x; \lambda) = e^{\nu(\lambda)x}( V_0(\lambda) + V_1^+(x; \lambda))$, where we define $V_1^+(x; \lambda) = \Theta^+(x; \lambda) V_0(\lambda)$. Similarly, define 
\begin{align*}
V^-(x; \lambda) &= P^-(x; \lambda) Z(x) = RP^+(-x; -\lambda)R e^{\nu(\lambda)x} V_0(\lambda) \\
&= e^{\nu(\lambda)x} R(I + \Theta^+(-x; -\lambda))R V_0(\lambda) = e^{\nu(\lambda)x}( V_0(\lambda) + R\Theta^+(-x; -\lambda) V_0(-\lambda) ),
\end{align*}
and let $V_1^-(x; \lambda) = R\Theta^+(-x; -\lambda) V_0(-\lambda)$. The decay rate for $V_1^\pm(x; \lambda)$ comes from the conjugation lemma.
\end{proof}
\end{lemma}

We use this result to prove the existence of $V^c(x)$ in part (i) of \cref{varadjsolutions}. Using a dimension counting argument, $\dim W^{cs}(0) = m + 1$ and $\dim W^{cu}(0) = m + 1$. Since $\Psi(0) \perp T_{Q(0)}W^{cs}(0) + T_{Q(0)}W^{cu}(0)$, $\dim T_{Q(0)}W^{cs}(0) + T_{Q(0)}W^{cu}(0) \leq 2m$, which implies that $\dim T_{Q(0)}W^{cs}(0) \cap T_{Q(0)}W^{cu}(0) = 2$. Since $\dim T_{Q(0)}W^s(0) \cap T_{Q(0)}W^s(0) = 1$ by \cref{nondegenlemma}, there exists $Y^0 \in T_{Q(0)}W^{cs}(0) \cap T_{Q(0)}W^{cu}(0)$ which is linearly independent from $Q'(0)$ with $Y^0 \notin T_{Q(0)}W^s(0) \cap T_{Q(0)}W^s(0)$. Define $V^\pm(x; 0)$ as in \cref{lemma:Vpm}. Since $V^-(0; 0) = R V^+(0; 0)$ and $V^\pm(0; 0) \in \Span \{Y^0 \}$, $V^+(0; 0) = V^-(0; 0)$, thus we define 
\[
V^c(x) = \begin{cases}
V^+(x; 0) & x \geq 0 \\
V^-(-x; 0) & x \leq 0,
\end{cases}
\]
from which it follows that $V^c(-x) = R V^c(x)$. In addition, by reversibility, $Y^- = R Y^+$, thus since $Y^+$ and $Y^-$ only have trivial intersection, $V^c(0)$ can contain no component in $Y^+ \oplus Y^-$. 

The next lemma evaluates two important inner products involving $V^\pm(0; \lambda)$.

\begin{lemma}\label{lemma:VpmPsiIP}
We have the following inner products
\begin{equation}\label{VpmIPs}
\begin{aligned}
\langle W_0, V^\pm(0; \lambda) \rangle &= 1 \mp \frac{1}{2} \tilde{M}\lambda + \mathcal{O}(|\lambda|^2) \\
\langle \Psi(0), V^\pm(0; \lambda) \rangle &= \mp \frac{1}{2} M_c \lambda + \mathcal{O}(|\lambda|^2),
\end{aligned}
\end{equation}
where
\begin{align*}
\tilde{M} &= \int_{-\infty}^{\infty} \left(v^c(y) - \frac{1}{c}\right) dy < \infty, \quad
M_c = \int_{-\infty}^\infty \partial_c q(y) dy,
\end{align*}
$q(y)$ is the first component of $Q(y)$, and $v^c(y)$ is the first component of $V^c(y)$.
\begin{proof}
Define $\tilde{V}^+(x; \lambda)$ by
\begin{equation}\label{deftildeV}
V^+(x; \lambda) = e^{\nu(\lambda)x}\tilde{V}^+(x; \lambda),
\end{equation}
so that $\tilde{V}^+(x, \lambda) \rightarrow V_0(\lambda)$ as $x \rightarrow \infty$. Differentiating with respect to $\lambda$ at $\lambda = 0$, 
\begin{equation}\label{VtildeVderiv}
\partial_\lambda V^+(x; 0)
= \partial_\lambda \tilde{V}^+(x; 0) + \frac{1}{c} x V^c(x),
\end{equation}
since $\tilde{V}^+(x; 0) = V^c(x)$ and $\nu'(0) = 1/c$. Substituting \cref{deftildeV} into \cref{vareqlambda}, differentiating with respect to $\lambda$ at $\lambda = 0$, and simplifying, 
\begin{align}\label{tildeVsolves}
[\partial_\lambda \tilde{V}^+(x; 0)]' &= A(Q(x))\partial_\lambda \tilde{V}^+(x; 0) + \left( B - \frac{1}{c}I\right) V^c(x).
\end{align}

For convenience, let $Y(x) = \partial_\lambda \tilde{V}^+(x; 0)$. Using the exponential trichotomy \cref{tildephitrich} for $\lambda = 0$ and noting that $\tilde{\Phi}^{c,+}(x, y; 0) = \langle W_0, \cdot \rangle V^c(x)$, we can formally write $Y(x)$ in integrated form as
\begin{equation}\label{Yintegform}
\begin{aligned}
Y(x) &= \tilde{\Phi}^{s,+}(x,0; 0)Y_0^s 
+ \int_0^x \tilde{\Phi}^{s,+}(x,y; 0)\left( B - \frac{1}{c}I \right) V^c(y) dy \\
&+ \int_{\infty}^x \tilde{\Phi}^{u,+}(x,y; 0)\left( B - \frac{1}{c}I \right) V^c(y) dy + V^c(x) \int_{\infty}^x \left\langle W_0, \left( B - \frac{1}{c}I \right) V^c(y) \right\rangle dy.
\end{aligned}
\end{equation}
Using the estimates \cref{stdtrichbounds}, the first and second integrals are finite for all $x$ since $V^c(x)$ is bounded. To prove that \cref{Yintegform} is a valid expression for $Y(x)$, it remains to show that the third integral is finite for all $x$. Using the expression for $V^c(x)$ from \cref{varadjsolutions}, 
\[
\left\langle W_0, \left( B - \frac{1}{c}I \right) V^c(y) \right\rangle = v^c(y) - \frac{1}{c},
\]
where $v^c(x)$ is the first component of $V^c(y)$. Since $|V^c(y) - V_0|\leq C e^{-\alpha y}$, and the first component of $V_0$ is $1/c$, 
\begin{align*}
\left| \int_{\infty}^x \left\langle W_0, \left( B - \frac{1}{c}I \right) V^c(y) \right\rangle dy \right| 
\leq C \int_x^{\infty} e^{-\alpha y} dy = C \frac{e^{-\alpha x}}{\alpha},
\end{align*}
which is finite, thus the formal expression \cref{Yintegform} is valid. 

Evaluating \cref{Yintegform} at $x = 0$ and taking the inner product with $W_0$, 
\begin{equation}\label{W0IPwithtildeV}
\langle W_0, \partial_\lambda \tilde{V}^+(x; 0) \rangle
= \langle W_0, V^c(0) \rangle \int_{\infty}^0 \left(v^c(y) - \frac{1}{c}\right) dy = -\frac{1}{2} \tilde{M},
\end{equation}
since $v^c(y)$ is an even function and $\langle W_0, V^c(0) \rangle = 1$. Expanding $V^+(0; \lambda)$ in a Taylor series about $\lambda = 0$ and using \cref{VtildeVderiv}, \cref{W0IPwithtildeV}, and $\langle W_0, V^c(0)\rangle = 1$, we obtain the first equation in \cref{VpmIPs}.

Evaluating \cref{Yintegform} at $x = 0$ and taking the inner product with $\Psi(0)$, 
\begin{align*}
\langle \Psi(0), \partial_\lambda \tilde{V}^+(x; 0)  \rangle
&= \int_{\infty}^0 \langle \Psi(0), \tilde{\Phi}^{u,+}(0,y; 0) B V^c(y) \rangle dy \\
&= \int_{\infty}^0 \langle \Psi(0), \tilde{\Phi}^{u,+}(0,0; 0) \tilde{\Phi}(0,y; 0) B V^c(y) \rangle dy,
\end{align*}
since $\tilde{\Phi}^{u,+}(0,y; 0) V^c(y) = 0$. For the projection $\tilde{\Phi}^{u,+}(0,0; 0)$, since $\ker \tilde{\Phi}^{u,+}(0, 0; 0)^* \perp \ran \tilde{\Phi}^{u,+}(0, 0; 0)$, $\ran \tilde{\Phi}^{u,+}(0, 0; 0)^* \perp \ker \tilde{\Phi}^{u,+}(0, 0; 0)$, and $\ker \tilde{\Phi}^{u,+}(0, 0; 0) = T_{Q(0)} W^s(0) \oplus T_{Q(0)} W^c(0)$, $\tilde{\Phi}^{u,+}(0, 0; 0)^*$ acts as the identity on $(T_{Q(0)} W^s(0) \oplus T_{Q(0)} W^c(0))^\perp$. Since $\Psi(0) \perp T_{Q(0)} W^s(0)$, $\Psi(0) \in \ran \tilde{\Phi}^{u,+}(0, 0; 0)^*$, thus $\tilde{\Phi}^{u,+}(0, 0; 0)^* \Psi(0) = \Psi(0)$. It follows that
\begin{align*}
\langle \Psi(0), \partial_\lambda \tilde{V}^+(x; 0)  \rangle
&= \int_{\infty}^0 \langle \tilde{\Phi}(0,y; 0)^* \tilde{\Phi}^{u,+}(0,0; 0)^* \Psi(0), B V^c(y) \rangle dy = -\int_0^\infty \langle \Psi(y), B V^c(y) \rangle dy \\
&= \int_0^\infty q(y) v^c(y) dy = \frac{1}{2}\int_{-\infty}^\infty q(y) v^c(y) dy ,
\end{align*}
since the last component of $\Psi(y)$ is $-q(y)$, and both $v^c(y)$ and $q(y)$ are even functions. Since $\calL(q)v^c = 1$ and $\calL(q)$ is self-adjoint, it follows from \cref{Lkernel} that
\begin{align*}
\langle \Psi(0), \partial_\lambda \tilde{V}^+(x; 0) \rangle
&= -\frac{1}{2} \int_{-\infty}^\infty \calL(q) \partial_c q(y) v^c(y) dy 
= -\frac{1}{2} \int_{-\infty}^\infty \partial_c q(y) \calL(q) v^c(y) dy \\
&= -\frac{1}{2} \int_{-\infty}^\infty \partial_c q(y) dy = -\frac{1}{2}M_c.
\end{align*}
Expanding $V^+(0; \lambda)$ in a Taylor series about $\lambda = 0$ and using \cref{VtildeVderiv}, \cref{W0IPwithtildeV}, and $\langle \Psi(0), V^c(0)\rangle = 0$, we obtain the second equation in \cref{VpmIPs}. The formulas for $V^-(0; \lambda)$ are similarly obtained.
\end{proof}
\end{lemma}

\subsection{Piecewise formulation}

As in \cite{Sandstede1998}, we will write the eigenvalue problem \cref{PDEeigsystemper3} as a piecewise system of equations. From \cref{th:perexist} and \cref{solvewithjumps}, the periodic $n$-pulse $Q_n(x)$ can be written piecewise as
\begin{equation}\label{Qnppiece}
\begin{aligned}
Q_i^-(x) &= Q^-(x; \beta_i^-) + \tilde{Q}_i^-(x) \qquad\qquad x \in [-X_{i-1}, 0] \\
Q_i^+(x) &= Q^+(x; \beta_i^+) + \tilde{Q}_i^+(x) \qquad\qquad x \in [0, X_i],
\end{aligned}
\end{equation}
where $Q_i^-: [-X_{i-1}, 0] \rightarrow \R^{2m+1}$ and $Q_i^+: [0, X_i] \rightarrow \R^{2m+1}$ are continuous, and the pieces are joined together end-to-end in a loop. We extend $Q_i^-(x)$ smoothly to $(-\infty, 0]$ and $Q_i^+(x)$ smoothly to $[0, \infty)$, so that $|Q_i^\pm(x)| \leq C e^{-\alpha_0 |x|}$. Next, we use the conjugation lemma to simplify \cref{PDEeigsystemper3} and to construct our piecewise ansatz. We will apply the conjugation lemma on $\R^\pm$ to the equation
\begin{equation}\label{Veqgeneral}
V(x)' = ( A(U(x)) + \lambda B)  V(x),
\end{equation}
where $U(x) \in C_b([0, \infty), \R^{2m+1})$ or $U(x) \in C_b((-\infty, 0], \R^{2m+1})$. Let $\Lambda = (U(x), \lambda)$ and $\Lambda_0 = (Q(x), 0)$. Using the conjugation lemma, there exists $\delta_4 \leq \delta_3$ and invertible linear transformations $P^\pm(x; U(x), \lambda) = I + \Theta^\pm(x; U(x), \lambda)$ such that for all $|\lambda| < \delta_4$ and $U(x)$ with $\| U(x) - Q(x) \| < \delta_4$, the change of coordinates $V = P^\pm(x; U(x), \lambda)$ on $\R^\pm$ conjugates \cref{Veqgeneral} into the constant-coefficient equation $Z'(x) = A(0; \lambda) Z(x)$. By \cref{solvewithjumps}, $\| Q_i^\pm(x) - Q(x) \| \leq C e^{-\alpha_0 X^*}$, thus we can choose $X^*$ sufficiently large so that $\| Q_i^\pm(x) - Q(x) \| \leq \delta_4$ for $i = 0, \dots, n-1$. Define
\begin{equation}\label{defPipm}
P_i^\pm(x; \lambda) = I + \Theta_i^\pm(x; \lambda) = I + \Theta^\pm(x; Q_i^\pm(x), \lambda).
\end{equation}
Expanding \cref{defPipm} in a Taylor series about $Q(x)$, and using \cref{conjthetadecay} and the estimates from \cref{solvewithjumps}, $\Theta_i^\pm(x; \lambda) = \Theta^\pm(x; \lambda) + \mathcal{O}(e^{-2 \alpha_0 X^*})$, thus we have the uniform estimates
\begin{equation}\label{Pipmest}
P_i^\pm(x; \lambda) = P^\pm(x; \lambda) + \mathcal{O}(e^{-2 \alpha_0 X^*}),
\end{equation}
where $P^\pm(x; \lambda)$ are the conjugation operators \cref{defPplus} and \cref{defPminus}.

Next, we note that when $\lambda = 0$, $W_0$ is a constant solution to $W'(x) = -A(Q_i^\pm(x); 0)^* W(x)$ for all $i$, thus by \cref{corr:adjconj}, $[P_i^\pm(x; \lambda)^{-1}]^* W_0$ is a solution to the constant coefficient equation $Z'(x) = -A(0)^* Z(x)$ for all $i$. Since $W_0$ is also a solution to this equation, and $[P_i^\pm(x; \lambda)^{-1}]^* W_0 \rightarrow W_0$ as $x \rightarrow \pm \infty$,
\begin{equation}\label{W0conjeq}
[P_i^\pm(x; \lambda)^{-1}]^* W_0 = W_0
\end{equation}
for $i = 0, \dots, n-1$ and all $x \in \R^\pm$.

Similarly to \cref{lemma:Vpm}, define
\begin{align}\label{defVipm}
V_i^\pm(x; \lambda) &= e^{\nu(\lambda)x} P_i^\pm(x; \lambda)V_0(\lambda) && x \in \R^\pm,
\end{align}
so that $V_i^\pm(x; \lambda)$ solves the equation
\begin{align}\label{PDEeigcenter}
[V_i^\pm(x; \lambda)]'(x) &= A(Q_i^\pm(x); \lambda) V_i^\pm(x; \lambda) && x \in \R^\pm.
\end{align}
Using \cref{defVipm} and \cref{Pipmest}, we have the estimate
\begin{align}\label{Vipmest}
V_i^\pm(x; \lambda) = V^\pm(x; \lambda) + \mathcal{O}(e^{-2 \alpha_0 X^*}).
\end{align}
Finally, for $i = 0, \dots, n-1$ define the constants
\begin{equation}\label{defklambda}
k_i^\pm(\lambda) = \frac{\langle Q'(0), V_i^\pm(0, \lambda) \rangle}{\langle Q'(0), Q'(0) \rangle},
\end{equation}
which are chosen so that for $i = 0, \dots, n-1$ and all $|\lambda| < \delta_4$,
\begin{equation}\label{klambdaIP}
\langle Q'(0), V_i^+(0; \lambda) - k_i^+(\lambda)Q'(0) \rangle = 
\langle Q'(0), V_i^-(0; \lambda) - k_i^-(\lambda)Q'(0) \rangle = 0.
\end{equation}

We can now construct our ansatz. It follows from \cref{Ekernel} and \cref{defAB} that
\begin{equation}\label{PDEeigkernel}
\begin{aligned}
{[\partial_x Q_n]}'(x) &= A(Q_n(x))\partial_x Q_n(x) \\
{[\partial_c Q_n]}'(x) &= A(Q_n(x))\partial_c Q_n(x) - B \partial_x Q_n(x).
\end{aligned}
\end{equation}
To exploit \cref{PDEeigkernel} and \cref{PDEeigcenter}, we take the piecewise ansatz for the eigenfunction $V(x)$
\begin{equation}\label{Vpiecewise}
\begin{aligned}
&d_i (\partial_x Q_i^-(x) - \lambda \partial_c Q_i^-(x)) + c_{i-1} e^{\nu(\lambda)X_{i-1} }( V_i^-(x; \lambda) - k_i^-(\lambda) Q'(x) ) + W_i^-(x) && x \in [-X_{i-1}, 0] \\
&d_i (\partial_x Q_i^+(x) - \lambda \partial_c Q_i^+(x)) + c_i e^{-\nu(\lambda)X_i }( V_i^+(x; \lambda) - k_i^+(\lambda) Q'(x) ) + W_i^+(x) && x \in [0, X_i],
\end{aligned}
\end{equation}
for $i = 0, \dots, n-1$, where $W_i^-(x) \in C([-X_{i-1}, 0], \C^{2m+1})$, $W_i^+(x) \in C([0, X_i], \C^{2m+1})$, and $c_i, d_i \in \C$. The subscripts are taken $\Mod n$, and the pieces are joined together end-to-end as in \cite{Sandstede1998}. The factors $e^{-\nu(\lambda)X_i}$ and $e^{\nu(\lambda)X_{i-1}}$ are chosen to facilitate joining the pieces at $\pm X_i$. Next, we substitute \cref{Vpiecewise} into \cref{PDEeigsystemper3}. Since the eigenfunction $V(x)$ must be continuous, the $2n$ pieces \cref{Vpiecewise} must satisfy $n$ matching conditions at $x = \pm X_i$ and $n$ matching conditions at $x = 0$. The remainder functions $W_i^\pm(x)$ must then satisfy the system of equations
\begin{equation}\label{eigsystem0}
\begin{aligned}
&(W_i^-)'(x) = A(Q_i^-(x); \lambda) W_i^-(x) + c_{i-1} e^{\nu(\lambda)X_{i-1}}  G_i^-(x; \lambda) + d_i \lambda^2 \tilde{H}_i^-(x) \\
&(W_i^+)'(x) = A(Q_i^+(x); \lambda) W_i^+(x) + c_i e^{-\nu(\lambda)X_i} G_i^+(x; \lambda) + d_i \lambda^2 \tilde{H}_i^+(x) \\
&W_i^+(X_i) - W_{i+1}^-(-X_i) = D_i d + C_i c \\
&W_i^+(0) - W_i^+(0) + c_i e^{-\nu(\lambda)X_i}(V_i^+(0; \lambda) - k_i^+(\lambda) Q'(0) )\\
&\qquad\qquad\qquad\qquad - c_{i-1} e^{\nu(\lambda)X_{i-1}} ( V_i^-(0; \lambda) - k_i^-(\lambda) Q'(0) ) = 0,
\end{aligned}
\end{equation}
for $i = 0, \dots, n-1$, where
\begin{equation}\label{defGH}
\begin{aligned}
G_i^\pm(x; \lambda) &= k_i^\pm(\lambda)\left( A(Q_i^\pm(x)) - A(Q(x)) - \lambda B \right) Q'(x) \\
\tilde{H}_i^\pm(x) &= -B \partial_c Q_i^\pm(x)\\
H(x) &= -B \partial_c Q(x)
\end{aligned}
\end{equation}
and
\begin{align}
D_i d &= d_{i+1}\left(\partial_x Q_{i+1}^-(-X_i) - \lambda \partial_c Q_{i+1}^-(-X_i)\right) - d_i \left( \partial_x Q_i^+(X_i) - \lambda \partial_c Q_i^+(X_i) \right) \label{defDid} \\
C_i c &= c_i \left( e^{\nu(\lambda) X_i} (V_{i+1}^-(-X_i; \lambda) - k_i^-(\lambda)Q'(-X_i)) - e^{-\nu(\lambda) X_i} (V_i^+(X_i; \lambda) - k_i^+(\lambda) Q'(X_i)) \right). \label{defCic}
\end{align}
As in \cite{Sandstede1998}, we will not be able to find a solution to this system for arbitrary $\lambda$. We will instead consider the system
\begin{equation}\label{eigsystem}
\begin{aligned}
(W_i^-)'(x) &= A(Q_i^-(x); \lambda) W_i^-(x) + c_{i-1} e^{\nu(\lambda)X_{i-1}}  G_i^-(x; \lambda) + d_i \lambda^2 \tilde{H}_i^-(x) \\
(W_i^+)'(x) &= A(Q_i^+(x); \lambda) W_i^+(x) + c_i e^{-\nu(\lambda)X_i}  G_i^+(x; \lambda) + d_i \lambda^2 \tilde{H}_i^+(x) \\
W_i^+(X_i) &- W_{i+1}^-(-X_i) = D_i d + C_i c \\
W_i^\pm(0) &\in \C \Psi(0) \oplus \C W_0 \oplus Y^+ \oplus Y^- \\
W_i^+(0) - W_i^-(0) &+ c_i e^{-\nu(\lambda)X_i}V_i^+(0; \lambda) - c_{i-1} e^{\nu(\lambda)X_{i-1}} V_i^-(0; \lambda) \in \C \Psi(0) \oplus \C W_0,
\end{aligned}
\end{equation}
for $i = 0, \dots, n-1$. The fourth equation states that the remainder functions $W_i^\pm(x)$ have no component in $\C Q'(0)$. The other terms in the ansatz do not appear in this equation by \cref{klambdaIP}. The final equation states that the jumps can only be in the directions of $W_0$ and $\Psi(0)$. The terms involving $Q'(0)$ in the ansatz do not appear in this equation since $Q'(0) \perp W_0\oplus\Psi(0)$. A solution to \cref{eigsystem} solves \cref{eigsystem0} if and only if $n$ jump conditions at $x = 0$ in the direction of $\C \Psi(0) \oplus \C W_0$ are satisfied, i.e.  
\begin{equation}\label{jumpxi1}
\begin{aligned}
\xi_i &= \langle \Psi(0), W_i^+(0) - W_i^-(0) + c_i e^{-\nu(\lambda)X_i}V_i^+(0; \lambda) - c_{i-1} e^{\nu(\lambda)X_{i-1}} V_i^-(0; \lambda) \rangle = 0  \\
\xi_i^c &= \langle W_0, W_i^+(0) - W_i^-(0) + c_i e^{-\nu(\lambda)X_i}V_i^+(0; \lambda) - c_{i-1} e^{\nu(\lambda)X_{i-1}} V_i^-(0; \lambda) \rangle = 0,
\end{aligned}
\end{equation}
for $i = 0, \dots, n-1$, where the terms involving $Q'(0)$ in the ansatz again do not appear.

Finally, we apply the conjugation operators \cref{defPipm} to the system \cref{eigsystem}. Making the substitution $W_i^\pm(x) = P_i^\pm(x; \lambda) Z_i^\pm(x)$, we obtain the equations
\begin{equation}\label{systemZ}
\begin{aligned}
(Z_i^-)'(x) &= A(0; \lambda) Z_i^-(x) + c_{i-1} e^{\nu(\lambda)X_{i-1}}  P_i^-(x; \lambda)^{-1} G_i^-(x; \lambda) + \lambda^2 d_i P_i^-(x; \lambda)^{-1} \tilde{H}_i^-(x) \\
(Z_i^+)'(x) &= A(0; \lambda) Z_i^+(x) + c_i e^{-\nu(\lambda)X_i}  P_i^+(x; \lambda)^{-1} G_i^+(x; \lambda) + \lambda^2 d_i P_i^+(x; \lambda)^{-1} \tilde{H}_i^+(x),
\end{aligned}
\end{equation}
with matching conditions at $x = \pm X_i$
\begin{equation}\label{systemmiddle}
P_i^+(X_i; \lambda) Z_i^+(X_i) - P_{i+1}^-(-X_i; \lambda) Z_{i+1}^-(-X_i; \lambda) = D_i d + C_i c, 
\end{equation}
and matching conditions at $x = 0$
\begin{equation}\label{systemcenter}
\begin{aligned}
&P_i^\pm(0; \lambda) Z_i^\pm(0) \in Y^+ \oplus Y^- \oplus \C \Psi(0) \oplus \C W_0 \\
&P_i^+(0; \lambda)Z_i^+(0) - P_i^-(0; \lambda) Z_i^-(0) \\
&\qquad+ c_i e^{-\nu(\lambda)X_i}V_i^+(0; \lambda) - c_{i-1} e^{\nu(\lambda)X_{i-1}} V_i^-(0; \lambda) \in \C \Psi(0) \oplus \C W_0.
\end{aligned}
\end{equation}
The jump conditions become
\begin{equation}\label{jumpcondZ}
\begin{aligned}
\xi_i &= \langle \Psi(0), P_i^+(0; \lambda) Z_i^+(0) - P_i^-(0; \lambda) Z_i^-(0) \\
&\qquad + c_i e^{-\nu(\lambda)X_i}V_i^+(0; \lambda) - c_{i-1} e^{\nu(\lambda)X_{i-1}} V_i^-(0; \lambda) \rangle = 0  \\
\xi_i^c &= \langle W_0, P_i^+(0; \lambda) Z_i^+(0) - P_i^-(0; \lambda) Z_i^-(0) \\
&\qquad + c_i e^{-\nu(\lambda)X_i}V_i^+(0; \lambda) - c_{i-1} e^{\nu(\lambda)X_{i-1}} V_i^-(0; \lambda) \rangle = 0. 
\end{aligned}
\end{equation}

To conclude this section, we collect some important estimates in the following lemma.

\begin{lemma}\label{stabestimateslemma}We have the estimates
\begin{enumerate}[(i)]
\item $|H(x)|, |\tilde{H}_i^\pm(x)| \leq C e^{-\alpha_0 |x|}$.
\item $|\tilde{H}_i^-(x) - H(x)| \leq C e^{-\alpha_0 X_{i-1}} e^{-\alpha_0(X_{i-1} + x) } + e^{-2 \alpha_0 X_i} e^{\alpha_0 x}$.
\item $|\tilde{H}_i^+(x) - H(x)| \leq C e^{-\alpha_0 X_i} e^{-\alpha_0(X_i - x) } + e^{-2 \alpha_0 X_{i-1}} e^{-\alpha_0 x}$.
\item $|G_i^\pm(x; \lambda)| \leq C |\lambda|(e^{-\alpha_0 X_i}+ |\lambda|) e^{-\alpha_0 |x|}$.
\item $D_i d = ( Q'(X_i) + Q'(-X_i))(d_{i+1} - d_i ) + \mathcal{O} ( e^{-\alpha_0 X_i} (e^{-\alpha_0 X^*} + |\lambda| )|d|)$.
\item $|C_i c| \leq C (e^{-\alpha_0 X_i} + |\lambda| ) |c|$.
\end{enumerate}
\begin{proof}
Since $H(x) = -\partial_c Q(x)$, it follows from \cref{transverseint} that $|H(x)| \leq C e^{-\alpha_0|x|}$, where we can use $\alpha_0$ in place of $\alpha_0 - \epsilon$ since $DF(0)$ is hyperbolic. The result for $\tilde{H}(x) = -\partial_c Q_i^\pm(x)$ can be similarly obtained by using Lin's method as in \cite{SandstedeStrut,Sandstede1993}. The bounds (ii) and (iii) follow from Lin's method and an adaptation of \cref{solvewithjumps} to derivatives with respect to $c$. For the estimates on $G_i^\pm(x; \lambda)$, by reversibility, $\langle Q'(0), V^\pm(0, 0) \rangle = \langle Q'(0), V^c(0) \rangle = 0$, thus $k_i^\pm(\lambda) = \mathcal{O}(|\lambda|)$. Using this together with the estimates from \cref{solvewithjumps} gives us (iv). For the estimate (v), we use \cref{VQpm} and \cref{Qpmbounds} with the derivative with respect to $x$ to get
\begin{align*}
(Q_{i+1}^-)'&(-X_i) = (Q^-)'(-X_i; \beta_{i+1}^-) + (\tilde{Q}_{i+1}^-)'(-X_i) \\
&= (Q^-)'(-X_i; \beta_{i+1}^-) + (Q^+)'(X_i; \beta_i^+) + \mathcal{O}(e^{-2 \alpha_0 X_i}) \\
&= Q'(-X_i) + Q'(X_i) + \mathcal{O}(e^{-\alpha_0 X_i}e^{-\alpha_0 X^*}).
\end{align*}
Similarly, $(Q_i^+)'(X_i) = Q'(-X_i) + Q'(X_i) + \mathcal{O}(e^{-\alpha_0 X_i}e^{-\alpha_0 X^*})$. Substitute these into \cref{defDid} and use (i) to get the estimate (v). For the estimate (vi), by \cref{lemma:Vpm}, $e^{\nu(\lambda) X_i} V^-(-X_i; \lambda) = V_0(\lambda) + \mathcal{O}(e^{-\alpha_0 X_i})$ and $e^{-\nu(\lambda) X_i} V^+(X_i; \lambda) = V_0(\lambda) + \mathcal{O}(e^{-\alpha_0 X_i})$. Subtracting these and using $k_i^\pm(\lambda) = \mathcal{O}(|\lambda|)$, we obtain the estimate (vi).
\end{proof}
\end{lemma}

\subsection{Exponential trichotomy}\label{sec:trichotomy}

We will now define exponential trichotomies on $\R^\pm$ for the conjugated system. Since there is some freedom in choosing subspaces for the trichotomy, we will make a choice that allows us to best satisfy equation \cref{systemmiddle}. The range of the stable projection on $\R^+$ is unique and is given by $E^s(\lambda)$, but we can choose any complement of $E^s(\lambda)$ to be the complement of the range of the stable projection at $X_i$. For sufficiently small $\lambda$, since the eigenvectors of $A(0; \lambda)$ are smooth in $\lambda$, $E^u(0)\oplus E^c(0)$ is a complement of $E^s(\lambda)$. Since $P_i^+(X_i; \lambda) = I + \mathcal{O}(e^{-\alpha_0 X_i})$, for sufficiently small $\lambda$ and sufficiently large $X_i$, $P_i^+(X_i, \lambda)^{-1} ( E^u(0)\oplus E^c(0) ) P_i^+(X_i, \lambda)$ is a complement of $E^s(\lambda)$. A similar result holds for the complement of the unstable projection on $\R^-$ at $-X_{i-1}$. Using these complements, there exists $\delta_5 \leq \delta_4$ such that for all $|\lambda| < \delta_5$, we can decompose the evolution operator $\Phi(x,y; \lambda)$ on $\R^\pm$ as
\begin{align}\label{Phidecomp}
\Phi(x, y; \lambda) &= \Phi_i^{s,\pm}(x, y; \lambda) + \Phi_i^{u,\pm}(x, y; \lambda) + \Phi_i^{c,\pm}(x, y; \lambda) && i = 0, \dots, n-1,
\end{align}
where
\begin{equation}\label{Zevolmod}
\begin{aligned}
\Phi_i^{s,+}(x, y; \lambda) &= \Phi(x, y; \lambda) P^s(\lambda) \\
\Phi_i^{u,+}(x, y; \lambda) &= \Phi(x, X_i; \lambda) P_i^+(X_i, \lambda)^{-1}
P^u(0) P_i^+(X_i, \lambda) \Phi(X_i, y; \lambda) \\
\Phi_i^{c,+}(x, y; \lambda) &= \Phi(x, X_i; \lambda) P_i^+(X_i, \lambda)^{-1}
P^c(0) P_i^+(X_i, \lambda) \Phi(X_i, y; \lambda)  \\
\Phi_i^{s,-}(x, y; \lambda) &= \Phi(x, -X_{i-1}; \lambda) P_i^-(-X_{i-1}, \lambda)^{-1}
P^s(0) P_i^-(-X_{i-1}, \lambda) \Phi(-X_{i-1}, y; \lambda) \\
\Phi_i^{u,-}(x, y; \lambda) &= \Phi(x, y; \lambda) P^u(\lambda) \\
\Phi_i^{c,-}(x, y; \lambda) &= \Phi(x, -X_{i-1}; \lambda) P_i^-(-X_{i-1}, \lambda)^{-1}
P^c(0) P_i^-(-X_{i-1}, \lambda) \Phi(-X_{i-1}, y; \lambda),
\end{aligned}
\end{equation}
and we have the estimates
\begin{equation}\label{Zevolbounds}
\begin{aligned}
|\Phi_i^{s,\pm}(x, y; \lambda)| &\leq C e^{-\alpha(x - y)} && y \leq x \\
|\Phi_i^{u,\pm}(x, y; \lambda)| &\leq C e^{-\alpha(y - x)} && x \leq y \\
|\Phi_i^{c,\pm}(x, y; \lambda)| &\leq C e^{\eta|x - y|},
\end{aligned}
\end{equation}
which are uniform in $\lambda$, and are the same as the estimates \cref{stdtrichbounds} but with possibly different constants $C$. The stable evolution $\Phi_i^{s,+}(x, y; \lambda)$ on $\R^+$ and the unstable evolution $\Phi_i^{u,-}(x, y; \lambda)$ on $\R^-$ do not depend on $X_i$ or $i$. The evolution of the unconjugated system $W_i^\pm(x) = A(Q_i^\pm(x); \lambda) W_i^\pm(x)$ on $\R^\pm$ is given by
\begin{equation}\label{unconjevol2}
\tilde{\Phi}_i^\pm(x, y; \lambda) = P_i^\pm(y; \lambda) \Phi(y, x; \lambda) P_i^\pm(x; \lambda)^{-1}.
\end{equation}
For $i = 0, \dots, n-1$, equations \cref{Zevolmod} induce an exponential trichotomy for the unconjugated system via the evolution operators
\begin{equation}\label{trichunconj}
\tilde{\Phi}_i^{s/u/c,\pm}(x, y; \lambda) = P_i^\pm(y; \lambda) \Phi_i^{s/u/c,\pm}(y, x; \lambda) P_i^\pm(x; \lambda)^{-1}.
\end{equation}
As a consequence of \cref{Zevolmod},
\begin{align*}
\tilde{\Phi}_i^{u,+}(X_i, X_i ; \lambda) &= P^u(0), \quad
\tilde{\Phi}_i^{c,+}(X_i, X_i ; \lambda) = P^c(0) \\
\tilde{\Phi}_i^{s,-}(-X_{i-1}, X_{i-1} ; \lambda) &= P^s(0), \quad
\tilde{\Phi}_i^{c,-}(-X_{i-1}, X_{i-1} ; \lambda) = P^c(0),
\end{align*}
which are the eigenprojections for $A(0)$ and are independent of $\lambda$.

\subsection{Inversion}

We will now solve the system \cref{eigsystem}. This follows the outline of the proof of \cite[Theorem 2]{Sandstede1998}, the main differences being the presence of a center subspace and the fact that we are on a periodic domain. Choose $\delta \leq \delta_5$. The result will hold for $|\lambda| < \delta$, and we may need to decrease $\delta$ as we proceed. Define the spaces
\begin{align*}
V_a &= \bigoplus_{i=0}^{n-1} E^u(0) \oplus E^s(0) \oplus E^c(0) \\
V_b &= \bigoplus_{i=0}^{n-1} (\C Q'(0)\oplus Y^-)\oplus(\C Q'(0)\oplus Y^+)\\
V_c &= \bigoplus_{i=0}^{n-1} \C \\
V_d &= \bigoplus_{i=0}^{n-1} \C \\
V_\lambda &= B_\delta(0) \subset \C,
\end{align*}
where the subscripts are taken $\Mod n$, since we are on a periodic domain, and the product spaces are endowed with the maximum norm. Using the variation of constants formula and splitting the evolution operator via the exponential trichotomy \cref{Zevolmod}, we write \cref{systemZ} in integrated form as
\begin{equation}\label{Zfpeq}
\begin{aligned}
Z_i^-(x) = \Phi_i^{s,-}&(x, -X_{i-1}; \lambda) P_i^-(-X_{i-1}; \lambda)^{-1} a_{i-1}^- + \Phi_i^{u,-}(x, 0; \lambda) P_i^-(0; \lambda)^{-1} b_i^- \\
&- \Phi_i^{c,-}(x, -X_{i-1}; \lambda) P_i^-(-X_{i-1}; \lambda)^{-1} a_{i-1}^c \\
&+ \int_0^x \Phi_i^{u,-}(x, y; \lambda)P_i^-(y; \lambda)^{-1}(c_{i-1} e^{\nu(\lambda)X_{i-1}} G_i^-(y; \lambda) + \lambda^2 d_i \tilde{H}_i^-(y)) dy \\
&+ \int_{-X_{i-1}}^x \Phi_i^{s,-}(x, y; \lambda) P_i^-(y; \lambda)^{-1} (c_{i-1} e^{\nu(\lambda)X_{i-1}} G_i^-(y; \lambda) + \lambda^2 d_i \tilde{H}_i^-(y)) dy \\
&+ \int_{-X_{i-1}}^x \Phi_i^{c,-}(x, y; \lambda) P_i^-(y; \lambda)^{-1} (c_{i-1} e^{\nu(\lambda)X_{i-1}} G_i^-(y; \lambda) + \lambda^2 d_i \tilde{H}_i^-(y)) dy \\ 
Z_i^+(x) = \Phi_i^{u,+}&(x, X_i; \lambda) P_i^+(X_i; \lambda)^{-1} a_i^+ + \Phi_i^{s,+}(x, 0; \lambda) P_i^+(0; \lambda)^{-1} b_i^+ \\
&+ \Phi_i^{c,+}(x, X_i; \lambda) P_i^+(X_i; \lambda)^{-1} a_i^c \\
&+ \int_0^x \Phi_i^{s,+}(x, y; \lambda) P_i^+(y; \lambda)^{-1} (c_i e^{-\nu(\lambda)X_i} G_i^+(y; \lambda) + \lambda^2 d_i \tilde{H}_i^+(y)) dy \\
&+ \int_{X_i}^x \Phi_i^{u,+}(x, y; \lambda) P_i^+(y; \lambda)^{-1}( c_i e^{-\nu(\lambda)X_i} G_i^+(y; \lambda) + \lambda^2 d_i \tilde{H}_i^+(y)) dy \\
&+ \int_{X_i}^x \Phi_i^{c,+}(x, y; \lambda) P_i^+(y; \lambda)^{-1}( c_i e^{-\nu(\lambda)X_i} G_i^+(y; \lambda) + \lambda^2 d_i \tilde{H}_i^+(y)) dy.
\end{aligned}
\end{equation}
As in \cite{Sandstede1998}, we will solve the eigenvalue problem in a series of inversion steps. Since the RHS of the fixed point equations \cref{Zfpeq} does not involve $Z_i^\pm$, these equations solve equation \cref{systemZ}. In the next lemma, we solve \cref{systemmiddle}, which are the matching conditions at the tails.

\begin{lemma}\label{Zinv1}
For $i = 0, \dots, n-1$, there is a unique set of initial conditions $(a_i^+, a_i^-, a_i^c)$ such that \cref{systemmiddle} is satisfied for any $(b, c, d)$ and $\lambda$. These are given by
\begin{equation}\label{aipmexp1}
\begin{aligned}
a_i^+ &= P_0^u(\lambda) D_i d + A_2(\lambda)_i^+(b, c, d) \\
a_i^- &= -P_0^s(\lambda) D_i d + A_2(\lambda)_i^-(b, c, d) \\
a_i^c &= \tilde{A}_2(\lambda)_i(b, c, d),
\end{aligned}
\end{equation}
where $A_2$ and $\tilde{A}_2$ are analytic in $\lambda$, linear in $(b, c, d)$, and have bounds
\begin{align}
|A_2(\lambda)_i(b, c, d)|
&\leq C \left(e^{-\alpha X_i}|b| + (e^{-\alpha_0 X_i} + |\lambda|)|c_i| + |\lambda|^2|d| \right) \label{A2bound} \\
|\tilde{A}_2(\lambda)_i(b, c, d) &\leq C |\lambda| \left( e^{-\alpha X_i} |b| + e^{-\alpha_0 X_i} |c_i| + e^{-\alpha_0 X_i} |\lambda|^2 |d| \right). \label{tildeA2bound}
\end{align}

\begin{proof}
It follows from \cref{Zevolmod} that
\begin{equation}\label{PhiXiProj}
\begin{aligned}
P_{i+1}^-(-X_i, \lambda) &\Phi_{i+1}^{s/c,-}(-X_i, -X_i; \lambda) P_{i+1}^-(-X_i, \lambda)^{-1} = P^{s/c}(0) \\
P_i^+(X_i, \lambda) &\Phi^{u/c,+}_{i}(X_i, X_i; \lambda) P_i^+(X_i, \lambda)^{-1} = P^{u/c}(0).
\end{aligned}
\end{equation}
Substituting \cref{Zfpeq} into \cref{systemmiddle} and using \cref{PhiXiProj}, we have
\begin{equation}\label{Didexpansion}
\begin{aligned}
D_i &d + C_i c = a_i^+ - a_i^- + 2 a_i^c \\
&+ P_i^+(X_i; \lambda)\Phi_i^{s,+}(X_i, 0; \lambda) P_i^+(0; \lambda)^{-1} b_i^+ - P_{i+1}^-(-X_i; \lambda)\Phi_i^{u,-}(-X_i, 0; \lambda) P_i^-(0; \lambda)^{-1} b_{i+1}^- \\
&+ P_i^+(X_i; \lambda) \int_0^{X_i} \Phi_i^{s,+}(X_i, y; \lambda)P_i^+(y; \lambda)^{-1} (c_i e^{-\nu(\lambda)X_i} G_i^+(y; \lambda) + \lambda^2 d_i \tilde{H}_i^+(y)) dy \\ 
&- P_{i+1}^-(-X_i; \lambda) \int_0^{-X_i} \Phi_i^{u,-}(-X_i, y; \lambda) P_{i+1}^-(y; \lambda)^{-1}(c_i e^{\nu(\lambda)X_i} G_{i+1}^-(y; \lambda) + \lambda^2 d_{i+1} \tilde{H}_{i+1}^-(y)) dy,
\end{aligned}
\end{equation}
which is of the form
\begin{align}\label{Dideq1}
D_i d + C_i c &= a_i^+ - a_i^- + 2 a_i^c + L_3(\lambda)_i(b, d).
\end{align}
$L_3(\lambda)_i(b, c, d)$ is defined by the RHS of \cref{Didexpansion}, is linear in $(b,c,d)$ and analytic in $\lambda$, and is independent of $a$. To obtain a bound on $L_3$, we will bound the individual terms involved. For the terms involving $b$, we use trichotomy estimates \cref{Zevolbounds} to get
\[
| P_i^+(X_i; \lambda)\Phi^s(X_i, 0; \lambda) P_i^+(0; \lambda)^{-1} P_i^+(0; \lambda)^{-1} b_i^+| \leq C e^{-\alpha X^*} |b|.
\]
The term involving $b_{i+1}^-$ is similar. For the integral terms, using \cref{Zevolbounds}, the bound \cref{expnubound}, and the estimates from \cref{stabestimateslemma},
\begin{align*}
&\left| P_i^+(X_i; \lambda) \int_0^{X_i} \Phi_i^{s,+}(X_i, y; \lambda)P_i^+(y; \lambda)^{-1} (c_i e^{-\nu(\lambda)X_i} G_i^+(y; \lambda) + \lambda^2 d_i \tilde{H}_i^+(y)) dy \right| \\
&\qquad \leq C \left( |\lambda|(e^{-\alpha_0 X^*}+ |\lambda|) |c| + |\lambda|^2 |d| \right) \int_0^{X_i} e^{-\alpha(X_i - y)} e^{-\alpha_0 y} dy \\
&\qquad \leq C e^{-\alpha_0 X_i}\left( |\lambda|(e^{-\alpha_0 X^*}+ |\lambda|) |c| + |\lambda|^2 |d| \right).
\end{align*}
The other integral is similar. Combining these, we have the bound for $L_3$
\begin{equation}\label{L3bound}
|L_3(\lambda)_i(b, c, d)| \leq C \left( e^{-\alpha X_i} |b| + e^{-\alpha_0 X_i} |\lambda| (e^{-\alpha_0 X_i } + |\lambda|) |c| + e^{-\alpha_0 X_i} |\lambda|^2 |d| \right).
\end{equation}
To solve for $a_i^+$, $a_i^-$, and $a_i^c$, we apply the projections on $P^{s/u/c}(0)$ on the eigenspaces $E^{s/u/c}(0)$ to \cref{Dideq1}. For $a_i^\pm$, using the bound \cref{L3bound} for $L_3(\lambda)(b, c, d)$ and the estimate for $C_i c$ from \cref{stabestimateslemma},
\begin{align*}
a_i^+ &= P_0^u(0) D_i d + A_2(\lambda)_i^+(b, c, d) \\
a_i^- &= -P_0^s(0) D_i d + A_2(\lambda)_i^-(b, c, d),
\end{align*}
where $A_2(\lambda)_i^\pm(b, c, d)$ has bound \cref{A2bound}. For $a_i^c$, we apply the center projection $P^c(0) = \langle W_0, \cdot \rangle$ to \cref{Dideq1}, and note that $\langle W_0, D_i d\rangle = 0$. Since $k_i^\pm(\lambda)Q'(\pm X_i) = \mathcal{O}(|\lambda|e^{-\alpha_0 X_i})$, and by reversibility
\begin{align*}
\langle W_0, &e^{\nu(\lambda) X_i} V^-(-X_i; \lambda) - e^{-\nu(\lambda) X_i} V^+(X_i; \lambda) \rangle = \left( P^-(-X_i; \lambda) - P^+(X_i; \lambda) \right) V_0(\lambda) \\
&= \langle W_0, R \Theta^+(X_i; 0) R - \Theta^+(X_i; 0), V_0 \rangle + \mathcal{O}(|\lambda|e^{-\alpha_0 X_i}) = \mathcal{O}(|\lambda|e^{-\alpha_0 X_i}),
\end{align*}
we have $|\langle W_0, C_i c\rangle| \leq |\lambda|e^{-\alpha_0 X_i}|c_i|$. For a bound on $\langle W_0, L_3(\lambda)_i(b,c,d) \rangle$, using \cref{W0conjeq} we obtain
\begin{align*}
\langle W_0, &P_i^+(X_i; \lambda)\Phi^s(X_i, 0; \lambda) P_i^+(0; \lambda)^{-1} b_i^+ \rangle
\\ 
&= \langle [P_i^+(X_i; 0)^{-1}]^* W_0, P_i^+(X_i; 0)\Phi^s(X_i, 0; 0) P_i^+(0; \lambda)^{-1} b_i^+ \rangle + \mathcal{O}(e^{-\alpha X_i} |\lambda||b_i^+|) \\
&= \langle W_0, \Phi^s(X_i, 0; 0)  P_i^+(0; \lambda)^{-1} b_i^+ \rangle + \mathcal{O}(e^{-\alpha X_i} |\lambda||b_i^+|) = \mathcal{O}(e^{-\alpha X_i} |\lambda||b_i^+|).
\end{align*}
Bounding the other terms in $L_3(\lambda)_i(b, c, d)$ in a similar fashion, we obtain the estimate
\begin{equation}\label{W0L3bound}
|\langle W_0, L_3(\lambda)_i(b, c, d) \rangle| \leq C |\lambda| \left( e^{-\alpha X_i} |b| + e^{-\alpha_0 X_i} |\lambda|(e^{-\alpha_0 X_i} + |\lambda|) |c| + e^{-\alpha_0 X_i} |\lambda|^2 |d| \right).
\end{equation}
Combining the above bounds and dividing by 2, $a_i^c = \tilde{A}_2(\lambda)_i(b, c, d)$, where $\tilde{A}_2(\lambda)_i(b, c, d)$ has bound \cref{tildeA2bound}.
\end{proof}
\end{lemma}

In the next lemma, we solve equations \cref{systemcenter}, which are matching conditions at $x = 0$ in the directions other than $\C \Psi(0) \oplus \C W_0$. Using the decomposition \cref{DSdecomp}, equations \cref{systemcenter} are equivalent to the three projections
\begin{equation}\label{centercond2}
\begin{aligned}
&P(\C Q'(0) ) P_i^-(0; \lambda) Z_i^-(0) = 0 \\
&P(\C Q'(0) ) P_i^+(0; \lambda) Z_i^+(0) = 0 \\
&P(Y^+ \oplus Y^-) ( P_i^+(0; \lambda)Z_i^+(0) - P_i^-(0; \lambda) Z_i^-(0) \\
&\qquad\qquad\qquad + c_i e^{-\nu(\lambda)X_i}V_i^+(0; \lambda) - c_{i-1} e^{\nu(\lambda)X_{i-1}} V_i^-(0; \lambda) ) = 0,
\end{aligned}
\end{equation}
where the kernel of each projection is the remaining spaces in the direct sum decomposition \cref{DSdecomp}. We do not need to include $\C Q'(0)$ in the third equation of \cref{centercond2} by \cref{klambdaIP} and since we eliminated any component of $P_i^\pm(0; \lambda) Z_i^\pm(0)$ in $\C Q'(0)$ in the first two equations.

\begin{lemma}\label{Zinv2}
There is a unique set of initial conditions $(a_i^+, a_i^-, a_i^c)$ for $i = 0, \dots, n-1$, and an operator $B_1: V_\lambda \times V_c \times V_d \rightarrow V_b$ such that for $a = (a_i^+, a_i^-, a_i^c)$ and $b = B_1(\lambda)(c,d)$, equations \cref{systemmiddle} and \cref{systemcenter} are satisfied for any $(c, d)$. $B_1(\lambda)(c, d)$ is analytic in $\lambda$, linear in $(c, d)$, and has uniform bound
\begin{align}
|B_1(\lambda)(c, d)| &\leq C\left( (e^{-\alpha_0 X^*} + |\lambda|)|c| + (|\lambda| + e^{-\alpha_0 X^*})^2 |d| \right). \label{B1bound}
\end{align} 
The initial conditions $a_i^+$, $a_i^-$, and $a_i^c$ are given by
\begin{align*}
a_i^+ &= P_0^u(0) D_i d + A_4(\lambda)_i^+(c, d) \\
a_i^- &= -P_0^s(0) D_i d + A_4(\lambda)_i^-(c, d) \\
a_i^c &= \tilde{A}_4(\lambda)_i(c, d),
\end{align*}
where $A_2$ and $\tilde{A}_2$ are analytic in $\lambda$, linear in $(c, d)$, and have bounds
\begin{align}
|A_4(\lambda)_i(b, c, d)|
&\leq C \left( (e^{-\alpha_0 X_i} + |\lambda|)|c_i| + |\lambda|^2|d| \right) \label{A4bound} \\
|\tilde{A}_4(\lambda)_i(b, c, d) &\leq C |\lambda| \left( e^{-\alpha_0 X_i} |c_i| + e^{-\alpha_0 X_i} |\lambda|^2 |d| \right). \label{tildeA4bound}
\end{align}

\begin{proof}
Using the decomposition \cref{TQ0decomp}, we can write $b_i^\pm$ uniquely as $b_i^\pm = x_i^\pm + y_i^\pm$, where $x_i^\pm \in \C Q'(0)$ and $y_i^\pm \in Y^\pm$. Using this together with \cref{Zevolmod}, we can write $P_i^\pm(0; \lambda) Z_i^\pm(0)$ as
\begin{equation}\label{PipmZ0}
\begin{aligned}
P_i^-&(0; \lambda) Z_i^-(0) = x_i^- + y_i^- + R_i^-(\lambda) b_i^- + P_i^-(0; \lambda) \Phi_i^{s,-}(0, -X_{i-1}; \lambda) P_i^-(-X_{i-1}, \lambda)^{-1} a_{i-1}^- \\
&- P_i^-(0; \lambda) \Phi_i^{c,-}(0, -X_{i-1}; \lambda) P_i^-(-X_{i-1}, \lambda)^{-1} a_{i-1}^c \\
&+ P_i^-(0; \lambda) \int_{-X_{i-1}}^0 \Phi_i^{s,-}(0, y; \lambda) P_i^-(y; \lambda)^{-1} (c_{i-1} e^{\nu(\lambda)X_{i-1}}
 G_i^-(y; \lambda) + \lambda^2 d_i \tilde{H}_i^-(y)) dy \\
&+ P_i^-(0; \lambda) \int_{-X_{i-1}}^0 \Phi_i^{c,-}(0, y; \lambda) P_i^-(y; \lambda)^{-1} (c_{i-1} e^{\nu(\lambda)X_{i-1}}
 G_i^-(y; \lambda) + \lambda^2 d_i \tilde{H}_i^-(y)) dy  \\ 
P_i^+&(0; \lambda) Z_i^+(0) = x_i^+ + y_i^+ + R_i^+(\lambda) b_i^+ + P_i^+(0; \lambda) \Phi_i^{u,+}(0, X_i; \lambda) P_i^+(X_i, \lambda)^{-1} a_i^+ \\
&+ P_i^+(0; \lambda) \Phi_i^{c,+}(0, X_i; \lambda) P_i^+(X_i, \lambda)^{-1} a_i^c \\
&+ P_i^+(0; \lambda) \int_{X_i}^0 \Phi_i^{u,+}(0, y; \lambda) P_i^+(y; \lambda)^{-1}( c_i e^{-\nu(\lambda)X_i} G_i^+(y; \lambda) + \lambda^2 d_i \tilde{H}_i^+(y)) dy \\
&+ P_i^+(0; \lambda) \int_{X_i}^0 \Phi_i^{c,+}(0, y; \lambda) P_i^+(y; \lambda)^{-1}( c_i e^{-\nu(\lambda)X_i} G_i^+(y; \lambda) + \lambda^2 d_i \tilde{H}_i^+(y)) dy,
\end{aligned}
\end{equation}
where 
\begin{equation}\label{Ripmbound}
\begin{aligned}
R_i^-(\lambda) &= P_i^-(0; \lambda) P^u(\lambda) P_i^-(0; \lambda)^{-1} 
- P^-(0; 0) P^u(0) P^-(0; 0)^{-1} = \mathcal{O}(|\lambda| + e^{-2 \alpha X^*}) \\
R_i^+(\lambda) &= P_i^+(0; \lambda) P^s(\lambda) P_i^+(0; \lambda)^{-1} 
- P^+(0; 0) P^s(0) P^+(0; 0)^{-1} = \mathcal{O}(|\lambda| + e^{-2 \alpha X^*}).
\end{aligned}
\end{equation}
Applying the projections in \cref{centercond2}, we obtain an expression of the form
\begin{equation}\label{projxy}
\begin{pmatrix}x_i^- \\ x_i^+ \\ 
y_i^+ - y_i^- \end{pmatrix} 
+ L_4(\lambda)_i(b, c, d) = 0,
\end{equation}
where $L_4(\lambda)_i(b, c, d)$ consists of the projections in \cref{centercond2} applied to the remaining terms in \cref{PipmZ0} as well as the term $P(Y^+ \oplus Y^-) ( c_i e^{-\nu(\lambda)X_i}V_i^+(0; \lambda) - c_{i-1} e^{\nu(\lambda)X_i} V_i^-(0; \lambda) )$ in the final component. To obtain a bound on $L_4$, we will bound the individual terms involved. For the $a_i^\pm$ terms, we use the expression from \cref{Zinv1}, the estimate \cref{A2bound}, and the trichotomy bounds \cref{Zevolbounds} to get
\begin{align*}
&|P_i^+(0; \lambda) \Phi_i^{u,+}(0, X_i; \lambda) P_i^+(X_i, \lambda)^{-1} a_i^+| \\
&\qquad\leq C e^{-\alpha X_i} \left( e^{-\alpha X^*} |b| + (e^{-\alpha_0 X^*} + |\lambda|)|c| + (|\lambda|^2 + |D|)|d| \right).
\end{align*}
The $a_i^-$ term is similar. For the $a_i^c$ terms, we use \cref{Zevolbounds}, the estimate \cref{expnubound}, and the expression from \cref{Zinv1} to get
\begin{align*}
| P_i^+(0; \lambda) \Phi_i^{c,+}(0, X_i; \lambda) P_i^+(X_i, \lambda)^{-1} a_i^c | &\leq C |\lambda| \left( e^{-\alpha X_i} |b| + e^{-\alpha X_i} |c_i| +|\lambda|^2 |d| \right).
\end{align*}
For the remainder terms involving $b_i$, we use \cref{Ripmbound} to get
\begin{align*}
|R_i^+(\lambda) b_i^+|\leq C \left(|\lambda| + e^{-2 \alpha X^*}\right)|b|.
\end{align*}
The term involving $b_i^-$ is similar. For the terms involving $c$ in the third equation in \cref{centercond2}, it follows from \cref{Vipmest} that $V_i^\pm(0; \lambda) = V^c(0) + \mathcal{O}(|\lambda| + e^{-2 \alpha_0 X^*})$. From the discussion following \cref{lemma:Vpm}, $V^c(0)$ contains no component in $Y^+ \oplus Y^-$, thus we have
\begin{align*}
|P(Y^+ \oplus Y^-) ( c_i e^{-\nu(\lambda)X_i}V_i^+(0; \lambda) - c_{i-1} e^{\nu(\lambda)X_{i-1}} V_i^-(0; \lambda) ) | \leq C (|\lambda| + e^{-2 \alpha_0 X^*})|c|,
\end{align*}
where we also used the estimate \cref{expnubound}. The bound on the integral terms is determined by the integral involving the center subspace, since the other integral has a stronger bound. Using the bounds from \cref{stabestimateslemma} together with the trichotomy bound \cref{Zevolbounds} and the bound \cref{expnubound}, 
\begin{align*}
&\left| P_i^+(0; \lambda) \int_{X_i}^0 \Phi_i^{c,+}(0, y; \lambda) P_i^+(y; \lambda)^{-1}( c_i e^{-\nu(\lambda)X_i} G_i^+(y; \lambda) + \lambda^2 d_i \tilde{H}_i^+(y)) dy  \right| \\
&\qquad \leq C \left( |\lambda|(e^{-\alpha_0 X^*} + |\lambda|)|c| + |\lambda|^2 |d| \right).
\end{align*}
The integral terms from $P_i^-(0; \lambda) Z_i^-(0)$ have similar bounds. Combining these and simplifying, we obtain the bound
\begin{align*}
|L_4(\lambda)_i(b, c, d)| \leq 
C\Big( (|\lambda| + e^{-\alpha X^*})|b| + (|\lambda|+e^{-\alpha_0 X^*})|c| + (|\lambda| + e^{-\alpha_0 X^*})^2 |d|  \Big).
\end{align*}
Since $|\lambda| < \delta$ and we can choose $X^*$ sufficiently large so that $e^{-\alpha X^*} < \delta$, this becomes
\begin{align*}
L_4(\lambda)(b, c, d) \leq 
C\Big( \delta |b| + (|\lambda|+e^{-\alpha_0 X^*})|c| + (|\lambda| + e^{-\alpha_0 X^*})^2 |d| \Big), 
\end{align*}
which is uniform in $|b|$. Define the map
\[
J_2: \left( \bigoplus_{j=1}^n \C Q'(0) \oplus \C Q'(0)  \right) \oplus
\left( \bigoplus_{j=1}^n Y^+ \oplus Y^- \right) 
\rightarrow \bigoplus_{j=1}^n \C Q'(0) \oplus \C Q'(0) \oplus (Y^+ \oplus Y^-)
\]
by $J_2( (x_i^+, x_i^-),(y_i^+, y_i^-))_i = ( x_i^+, x_i^-, y_i^+ - y_i^- )$, which is an isomorphism by \cref{DSdecomp}. Since $b_i = (x_i^- + y_i^-, x_i^+ + y_i^+)$, we can write \cref{projxy} as
\begin{equation}\label{projxy2}
J_2( (x_i^+, x_i^-),(y_i^+, y_i^-))_i 
+ L_4(\lambda)_i(b_i, 0, 0) + L_4(\lambda)_i(0, c, d) = 0.
\end{equation}
Let $S_2(b)_i = J_2( (x_i^+, x_i^-),(y_i^+, y_i^-))_i + L_4(\lambda)_i(b_i, 0, 0)$. Substituting this into \cref{projxy2}, we obtain the equation $S_2(b) = -L_4(\lambda)(0, c, d)$. Decreasing $\delta$ if necessary, the operator $S_2(b)$ is invertible, thus we can solve for $b$ by using
\begin{align}
b = B_1(\lambda)(c,d) 
= -S_2^{-1} L_4(\lambda)(0, c, d),
\end{align}
which has bound given by \cref{B1bound}. Substituting \cref{B1bound} into the bounds \cref{A2bound}, and \cref{tildeA2bound} and using the estimates for $C_i c$ and $D_i d$ from \cref{stabestimateslemma}, we obtain the bounds \cref{A4bound} and \cref{tildeA4bound}.
\end{proof}
\end{lemma}

\subsection{Jump Conditions}

We have constructed a unique solution to \cref{eigsystem0} which will have $n$ jumps in the directions of $\Psi(0)$ and $W_0$. For this solution to be an eigenfunction, all $n$ jumps must be 0. In the next two lemmas, we compute the jumps in the direction of $W_0$ and $\Psi(0)$. 

\begin{lemma}\label{jumpcenteradj}
The jumps in the direction of $W_0$ are given 
\begin{equation}\label{jumpW0}
\begin{aligned}
\xi^c_i &= e^{-\nu(\lambda)X_i}c_i\left( 1 - \frac{1}{2}\lambda \tilde{M} \right) 
- e^{\nu(\lambda)X_{i-1}}c_{i-1}\left( 1 + \frac{1}{2}\lambda \tilde{M} \right) \\
&\qquad + \mathcal{O}\Big( (e^{-\alpha_0 X^*} + |\lambda|)^2 |c| + |\lambda| (|\lambda| + e^{-\alpha_0 X^*})^2 |d| \Big),
\end{aligned}
\end{equation}
for $i = 0, \dots, n-1$, where
\[
\tilde{M} = \int_{-\infty}^{\infty} \left(v^c(y) - \frac{1}{c}\right) dy,
\]
and the remainder terms are analytic in $\lambda$.

\begin{proof}
The jumps in the center direction are given by
\[
\xi_i^c = 
\langle W_0, P_i^+(0; \lambda) Z_i^+(0) - P_i^-(0; \lambda) Z_i^-(0) + c_i e^{-\nu(\lambda)X_i}V_i^+(0; \lambda) - c_{i-1} e^{\nu(\lambda)X_{i-1}} V_i^-(0; \lambda) \rangle.
\]
Using \cref{lemma:VpmPsiIP} and the bound \cref{expnubound},
\begin{align*}
\langle W_0, &c_i e^{-\nu(\lambda)X_i}V_i^+(0; \lambda) - c_{i-1} e^{\nu(\lambda)X_{i-1}} V_i^-(0; \lambda) \rangle \\
&= e^{-\nu(\lambda)X_i}c_i\left( 1 - \frac{1}{2}\lambda \tilde{M} \right) 
- e^{\nu(\lambda)X_{i-1}}c_{i-1}\left( 1 + \frac{1}{2}\lambda \tilde{M} \right) + \mathcal{O}\left( (e^{-\alpha_0 X^*} +|\lambda|)^2 |c| \right).
\end{align*}
The terms $P_i^\pm(0; \lambda) Z_i^\pm(0)$ are given by \cref{PipmZ0}. The only leading order term involves the integral of $\tilde{H}_i^\pm$ in the center subspace. Using \cref{Zevolmod}, \cref{Pipmest}, \cref{deftildephi}, and \cref{lemma:Afacts},
\begin{equation*}
\begin{aligned}
&\left\langle W_0, P_i^-(0; \lambda) \lambda^2 d_i \int_{-X_{i-1}}^0 \Phi_i^{c,-}(0, y; \lambda) P_i^-(y; \lambda)^{-1} \tilde{H}_i^-(y) dy \right\rangle \\
&\qquad = \lambda^2 d_i \int_{-X_{i-1}}^0 \langle W_0, P_i^-(0; \lambda) \Phi(0, -X_{i-1}; \lambda) P_i^-(-X_{i-1}, \lambda)^{-1} \\
&\qquad\qquad\qquad
P^c(0) P_i^-(-X_{i-1}, \lambda) \Phi(-X_{i-1}, y; \lambda) P_i^-(y; \lambda)^{-1} \tilde{H}_i^-(y) \rangle dy \\
&\qquad=\lambda^2 d_i \int_{-X_{i-1}}^0 \langle \tilde{\Phi}^-(y, -X_{i-1}; 0)^* P^c(0)^* \tilde{\Phi}^-(-X_{i-1}, 0; 0)^* W_0,
 H(y) \rangle dy \\
& \qquad\qquad\qquad+ \mathcal{O}(|\lambda|^2 (|\lambda| + e^{- \alpha_0 X^*}) |d|).
\end{aligned}
\end{equation*}
Since $W_0$ is a constant solution to \cref{adjvareq2}, $P^c(0)^* W_0 = W_0$, and $H(y) = -B \partial_c Q(y)$. 
\begin{equation*}
\begin{aligned}
&\left\langle W_0, P_i^-(0; \lambda) \lambda^2 d_i \int_{-X_{i-1}}^0 \Phi_i^{c,-}(0, y; \lambda) P_i^-(y; \lambda)^{-1} \tilde{H}_i^-(y) dy \right\rangle \\
&\qquad=\lambda^2 d_i \int_{-\infty}^0 \langle W_0, H(y) \rangle dy + \mathcal{O}(|\lambda|^2 (|\lambda| + e^{- \alpha_0 X^*}) |d|) \\
&\qquad= -\lambda^2 d_i \int_{-\infty}^0 \partial_c q(y)dy + \mathcal{O}(|\lambda|^2 (|\lambda| + e^{- \alpha_0 X^*}) |d|).
\end{aligned}
\end{equation*}
Similarly, 
\begin{align*}
&\left\langle W_0, P_i^-(0; \lambda) \int_{-X_{i-1}}^0 \Phi_i^{c,-}(0, y; \lambda) P_i^-(y; \lambda)^{-1} \lambda^2 d_i \tilde{H}_i^-(y) dy  \right\rangle \\
&\qquad= \lambda^2 d_i \int_0^{\infty} \partial_c q(y)dy + \mathcal{O}(|\lambda|^2 (|\lambda| + e^{- \alpha_0 X^*}) |d|).
\end{align*}

The rest of the terms are higher order. For the terms involving $a_i^\pm$, we use \cref{Zinv2}, \cref{Zevolmod}, \cref{Zevolbounds}, and \cref{A4bound} to get
\begin{align*}
\langle W_0, &P_i^+(0; \lambda) \Phi_i^{u,+}(0, X_i; \lambda) P_i^+(X_i, \lambda)^{-1} a_i^+ \rangle = \langle W_0, \tilde{\Phi}_i^+(0, X_i; \lambda) 
P^u(0) \tilde{\Phi}_i^+(X_i, X_i; \lambda)a_i^+ \rangle \\
&= \langle \tilde{\Phi}_i^+(X_i, 0; 0)^* W_0, P^u(0) a_i^+ \rangle + \mathcal{O}(|\lambda|e^{-\alpha X^*}|a_i^+|) = \langle W_0, P^u(0) a_i^+ \rangle + \mathcal{O}(|\lambda|e^{-\alpha X^*}|a_i^+|) \\
&= \mathcal{O}\left(|\lambda|e^{-\alpha X^*} \left( (e^{-\alpha_0 X^*} + |\lambda|^2) |c| + |\lambda|^2 |d| + |D||d| \right) \right),
\end{align*}
since $W_0$ is a constant solution to $W'(x) = -A(Q_i^\pm(x))^* W(x)$. For the terms involving $a_i^c$, we use \cref{Zinv2}, \cref{Zevolmod}, \cref{Zevolbounds}, \cref{tildeA4bound}, and \cref{expnubound} to get
\begin{align*}
\langle W_0, &P_i^+(0; \lambda) \Phi_i^{c,+}(0, X_i; \lambda) P_i^+(X_i, \lambda)^{-1} a_i^c \rangle = \langle W_0, \tilde{\Phi}_{i,+}(0, X_i; \lambda) 
P^c(0) \tilde{\Phi}_i^{+}(X_i, X_i; \lambda)a_i^c \rangle \\
&= \langle \tilde{\Phi}_i^+(X_i, 0; 0)^* W_0, P^c(0) a_i^c \rangle + \mathcal{O}(|\lambda||a_i^c|) = \mathcal{O}\left( |\lambda| (e^{-\alpha_0 X^*} + |\lambda|) |c| +|\lambda|^3 |d| \right).
\end{align*}
For the terms involving $b$, $b_i^\pm = x_i^\pm + y_i^\pm$ vanishes when we take the inner product with $W_0$. For the remaining terms, we use the estimate \cref{B1bound} to get
\begin{align*}
|\langle W_0, R_i^+(\lambda) b_i^+ \rangle |\leq C \left(|\lambda| + e^{-2 \alpha_0 X^*}\right)\left((|\lambda|+e^{-\alpha_0 X^*})|c| + (|\lambda| + e^{-\alpha_0 X^*})^2 |d| \right).
\end{align*}
For the center integral involving $G_i^\pm(y, \lambda)$, we use \cref{stabestimateslemma} and \cref{expnubound} to get
\begin{align*}
\left| \left\langle W_0, P_i^+(0; \lambda) \int_{X_i}^0 \Phi_i^{c,+}(0, y; \lambda) P_i^+(y; \lambda)^{-1} c_i e^{-\nu(\lambda)X_i} G_i^+(y; \lambda)  dy \right\rangle \right| \leq C |\lambda| (e^{-\alpha_0 X^*} + |\lambda|) |c|.
\end{align*}
The non-center integral involving $G_i^\pm(y, \lambda)$ has a stronger bound. For the non-center integral involving $\tilde{H}_i^\pm$, we follow the same procedure as above, replacing $P^c(0)$ by $P^u(0)$, to get
\begin{align*}
\left| \left\langle W_0, P_i^+(0; \lambda) \int_{X_i}^0 \Phi_i^{u,+}(0, y; \lambda) P_i^+(y; \lambda)^{-1} \lambda^2 d_i \tilde{H}_i^+(y) dy \right\rangle \right| \leq C |\lambda|^2 (|\lambda| + e^{- \alpha_0 X^*}) |d|.
\end{align*}
Bounds for the terms from $P_i^-(0; \lambda) Z_i^-(0)$ are similar. Combining all of these terms and simplifying, we obtain the center jump expressions \cref{jumpW0}.
\end{proof}
\end{lemma}

\begin{lemma}\label{jumpadj}
The jumps in the direction of $\Psi(0)$ are given 
\begin{equation}\label{jumpPsi0}
\begin{aligned}
\xi_i = \langle \Psi&(X_i), Q'(-X_i) \rangle (d_{i+1} - d_i ) - \langle \Psi(X_{i-1}), Q'(-X_{i-1}) \rangle (d_i - d_{i-1} ) - \lambda^2 d_i M  \\
&-\frac{1}{2}\lambda M_c \left( e^{-\nu(\lambda)X_i}c_i + e^{\nu(\lambda)X_{i-1}}c_{i-1} \right) + \mathcal{O}\left( (e^{-\alpha_0 X^*} + |\lambda|)^2 |c| + (e^{-\alpha_0 X^*} + |\lambda|)^3 |d| \right),
\end{aligned}
\end{equation}
for $i = 0, \dots, n-1$, where
\[
M = \int_{-\infty}^\infty q(y) \partial_c q(y) dy, \qquad
M_c = \int_{-\infty}^\infty \partial_c q(y) dy,
\]
and the remainder terms are analytic in $\lambda$.
\begin{proof}
The jumps in the direction of $\Psi(0)$ are given by
\[
\xi_i = 
\langle \Psi(0), P_i^+(0; \lambda) Z_i^+(0) - P_i^-(0; \lambda) Z_i^-(0) + c_i e^{-\nu(\lambda)X_i}V_i^+(0; \lambda) - c_{i-1} e^{\nu(\lambda)X_{i-1}} V_i^-(0; \lambda) \rangle.
\]
Using \cref{lemma:VpmPsiIP} and the bound \cref{expnubound},
\begin{align*}
\langle \Psi(0), &c_i e^{-\nu(\lambda)X_i}V_i^+(0; \lambda) - c_{i-1} e^{\nu(\lambda)X_i} V_i^-(0; \lambda) \rangle \\
&= -\frac{1}{2}\lambda M_c e^{-\nu(\lambda)X_i}c_i - \frac{1}{2}\lambda M_c e^{\nu(\lambda)X_{i-1}}c_{i-1} + \mathcal{O}\left( (e^{-\alpha X^*} +|\lambda|)^2 |c| \right). 
\end{align*}
The terms $P_i^\pm(0; \lambda) Z_i^\pm(0)$ are given by \cref{PipmZ0}. As in \cref{jumpcenteradj}, we begin by computing the leading order terms. For the non-center integral involving $\tilde{H}_i^\pm$, following the same procedure as in \cref{jumpcenteradj},
\begin{align*}
&\left\langle \Psi(0), P_i^-(0; \lambda) \int_{-X_{i-1}}^0 \Phi_i^{s,-}(0, y; \lambda) P_i^-(y; \lambda)^{-1} \lambda^2 d_i \tilde{H}_i^-(y) dy \right\rangle \\
&\qquad = \lambda^2 d_i \int_{-X_{i-1}}^0 \langle \tilde{\Phi}^-(y, -X_{i-1}; 0)^* P^s(0)^* \tilde{\Phi}^-(-X_{i-1}, 0; 0)^* \Psi(0), H(y) \rangle dy \\
&\qquad \qquad \qquad+ \mathcal{O}(|\lambda|^2( |\lambda| + {e^{-\alpha_0 X^*}})|d_i|).
\end{align*}
By \cref{lemma:Afacts}, $P^s(0)^* = P^{u,*}(0)$, and $P^{u,*}(0) \Psi(-X_{i-1}) = \Psi(-X_{i-1}) + \mathcal{O}(e^{-2\alpha_0 X_{i-1}})$. Since $\Psi(x)$ is a solution to \cref{adjvareq2}, $H(y) = -B \partial_c Q(y)$, and the last component of $\Psi(y)$ is $-q(y)$, this becomes
\begin{align*}
&\left\langle \Psi(0), P_i^-(0; \lambda) \int_{-X_{i-1}}^0 \Phi_i^{s,-}(0, y; \lambda) P_i^-(y; \lambda)^{-1} \lambda^2 d_i \tilde{H}_i^-(y) dy \right\rangle \\
&\qquad\qquad = \lambda^2 d_i \int_{-\infty}^0 q(y) \partial_c q(y) dy + \mathcal{O}(|\lambda|^2( |\lambda| + {e^{-\alpha_0 X^*}})|d_i|).
\end{align*}
Similarly,
\begin{align*}
&\left\langle \Psi(0), P_i^+(0; \lambda) \int_{X_i}^0 \Phi_i^{u,+}(0, y; \lambda) P_i^+(y; \lambda)^{-1} \lambda^2 d_i \tilde{H}_i^+(y) dy \right\rangle \\
&\qquad\qquad = -\lambda^2 d_i \int_0^\infty q(y) \partial_c q(y) dy + \mathcal{O}(|\lambda|^2( |\lambda| + {e^{-\alpha_0 X^*}})|d_i|).
\end{align*}
For the terms involving $a_i^\pm$, we use \cref{Zinv2}, \cref{Zevolmod}, and \cref{A4bound} to get
\begin{align*}
\langle \Psi(0), &P_i^+(0; \lambda) \Phi_i^{u,+}(0, X_i; \lambda) P_i^+(X_i; \lambda)^{-1} a_i^+ \rangle \\
&= \langle \Psi(0), P_i^+(0; \lambda) \Phi^+(0, X_i; \lambda) P_i^+(X_i; \lambda)^{-1}
P^u(0) a_i^+ \rangle \\
&= \langle \Psi(0), \tilde{\Phi}^+(0, X_i; 0) P^u(0) a_i^+ \rangle + \mathcal{O}(e^{-\alpha_0 X^*}(|\lambda| + e^{-2 \alpha_0 X^*})|a_i^+|) \\
&= \langle \Psi(X_i), P^u(0) D_i d \rangle + \mathcal{O}\left(e^{-\alpha_0 X^*} (e^{-\alpha_0 X^*} + |\lambda|)|c| + e^{-\alpha_0 X^*}(|\lambda| + e^{-\alpha_0 X^*})^2 |d|  \right).
\end{align*}
By \cref{stabestimateslemma}, $P^u(0) D_i d = Q'(-X_i)(d_{i+1} - d_i ) + \mathcal{O}( e^{-\alpha_0 X_i} (e^{-\alpha_0 X^*} + |\lambda| )|d|)$, thus we have
\begin{align*}
\langle \Psi(0), &P_i^+(0; \lambda) \Phi_i^{u,+}(0, X_i; \lambda) P_i^+(X_i, \lambda)^{-1} a_i^+ \rangle = \langle \Psi(X_i), Q'(-X_i) \rangle (d_{i+1} - d_i ) \\
&+ \mathcal{O}\left(e^{-\alpha_0 X^*} (e^{-\alpha_0 X^*} + |\lambda|)|c| + e^{-\alpha_0 X^*}(|\lambda| + e^{-\alpha_0 X^*})^2 |d|  \right).
\end{align*}
Similarly,
\begin{align*}
\langle \Psi(0), &P_i^-(0; \lambda) \Phi_i^{s,-}(0, -X_{i-1}; \lambda) P_i^-(-X_{i-1}, \lambda)^{-1} a_{i-1}^- \rangle = -\langle \Psi(-X_{i-1}), Q'(X_{i-1}) \rangle (d_i - d_{i-1} ) \\
&+ \mathcal{O}\left(e^{-\alpha_0 X^*} (e^{-\alpha_0 X^*} + |\lambda|)|c| + e^{-\alpha_0 X^*}(|\lambda| + e^{-\alpha_0 X^*})^2 |d| \right).
\end{align*}

The remaining terms will be higher order. For the terms involving $a_i^c$, similar to \cref{jumpcenteradj},
\begin{align*}
\langle \Psi(0), &P_i^+(0; \lambda) \Phi_i^{c,+}(0, X_i; \lambda) P_i^+(X_i, \lambda)^{-1} a_i^c \rangle = \mathcal{O}\left( |\lambda| (e^{-\alpha_0 X^*} + |\lambda|)|c| + |\lambda^3 |d| \right).
\end{align*}
For the terms involving $b$, $b_i^\pm = x_i^\pm + y_i^\pm$ vanishes when we take the inner product with $\Psi(0)$. For the remaining terms, similar to \cref{jumpcenteradj},
\begin{align*}
|\langle \Psi(0), R_i^+(\lambda) b_i^+ \rangle |\leq C \left(|\lambda| + e^{-2 \alpha X^*}\right)\left((|\lambda|+e^{-\alpha_0 X^*})|c| + (|\lambda| + e^{-\alpha X^*})^2 |d| \right).
\end{align*}
For the center integral term, using the bounds from \cref{stabestimateslemma},
\begin{align*}
&\left\langle \Psi(0), P_i^+(0; \lambda) \int_{X_i}^0 \Phi_i^{c,+}(0, y; \lambda) P_i^+(y; \lambda)^{-1}( c_i e^{-\nu(\lambda)X_i} G_i^+(y; \lambda) + \lambda^2 d_i \tilde{H}_i^+(y)) dy \right\rangle \\
&\qquad= \int_{-X_{i-1}}^0 \langle \Psi(-X_{i-1}), P^c(0) \tilde{\Phi}^-(-X_{i-1}, y; 0) 
 (c_i e^{-\nu(\lambda)X_i} G_i^+(y; \lambda) + \lambda^2 d_i H(y)) \rangle dy \\
&\qquad\qquad +\mathcal{O}\left( (|\lambda| + e^{-2\alpha_0 X^*}) \left(|\lambda|(|\lambda + e^{-\alpha_0 X^*})|c| + |\lambda|^2 |d| \right) \right) \\
&\qquad= \mathcal{O}\left( |\lambda|(|\lambda| + e^{-\alpha_0 X^*})^2|c| + |\lambda|^2 (|\lambda| + e^{-\alpha_0 X^*}) |d| \right).
\end{align*}
The non-center integral involving $G_i^\pm(x; \lambda)$ has a stronger bound. Bounds for the terms from $P_i^-(0; \lambda) Z_i^-(0)$ are similar. Combining all of these terms, using the reversibility relation
\[
\langle \Psi(-X_{i-1}), Q'(X_{i-1}) \rangle = \langle R \Psi(-X_{i-1}), R Q'(X_{i-1}) \rangle = 
-\langle \Psi(X_{i-1}), Q'(-X_{i-1}) \rangle,
\]
and simplifying, we obtain the jump expressions \cref{jumpPsi0}.
\end{proof}
\end{lemma}

\cref{blockmatrixtheorem} combines the jump conditions from \cref{jumpcenteradj} and \cref{jumpadj} into a single block matrix equation $E(\lambda) = 0$ which is analytic in $\lambda$. A nontrivial solution exists if and only if $\det E(\lambda) = 0$.

\section{Proof of results for periodic single pulse }\label{sec:singlepulse}

\subsection{Proof of Lemma \ref{lemma:1blockmatrix} }

For the periodic single pulse, there is one length parameter $X_0 = X$, and for the ansatz \cref{Vpiecewise}, there is only one parameter $c$ and one parameter $d$. The system of equations \cref{eigsystem} becomes
\begin{equation}\label{eigsystemper1p}
\begin{aligned}
&(W^\pm)'(x) = A(Q_0^\pm(x); \lambda) W^\pm(x) + c e^{\mp \nu(\lambda) X} G_0^\pm(x; \lambda) + d \lambda^2 \tilde{H}_0^\pm(x)  \\
&W^+(X) - W^-(-X) = C_0 c \\
&W^\pm(0) \in \R \Psi(0) \oplus \R W_0 \oplus Y^+ \oplus Y^- \\
&W^+(0) - W^-(0) + c ( e^{-\nu(\lambda)X}V_0^+(0; \lambda) - e^{\nu(\lambda)X} V_0^-(0; \lambda)) \in \C \Psi(0) \oplus \C W_0.
\end{aligned}
\end{equation}
Equation \cref{1pblockmatrix} follows directly from \cref{blockmatrixtheorem}. The bounds for the remainder terms follow from \cref{blockmatrixtheorem} and the fact that $d$ only appears as $\lambda^2 d$ in \cref{eigsystemper1p}. 

To show symmetry, suppose that $(W^+(x; c, d, \lambda), W^-(x; c, d, \lambda))$ is a solution to \cref{eigsystemper1p}. Since the periodic single pulse is symmetric, $Q_0^-(x) = R Q_0^+(-x)$, from which it follows that $G_0^-(x; \lambda) = -R G_0^+(-x; -\lambda)R$. Since $\nu(-\lambda) = -\nu(\lambda)$, $V_0^-(x; \lambda) = R V_0^+(-x; -\lambda)$, $Y^+ = R Y^-$, $R \Psi(0) = \Psi(0)$, and $R W_0 = W_0$, if we replace $(c, d, \lambda)$ by $(-c, d, -\lambda)$, equations \cref{eigsystemper1p} are satisfied by $(-RW^-(-x; c, d, \lambda), -RW^+(-x; c, d, \lambda))$. Since \cref{eigsystemper1p} has a unique solution, $\left(W^+(x; c, d, -\lambda), W^-(x; c, d, -\lambda)\right)
= -\left(RW^-(-x; -c, d, \lambda), RW^+(-x; -c, d, \lambda)\right)$. The jump conditions then become
\begin{align*}
\langle W_0, W^+(0; c, d, -\lambda) - W^-(0; c, d, -\lambda) \rangle &= \langle W_0, W^+(0; -c, d, \lambda) - W^-(0; -c, d, \lambda) \rangle \\
\langle \Psi(0), W^+(0; c, d, -\lambda) - W^-(0; c, d, -\lambda) \rangle &= \langle \Psi(0), W^+(0; -c, d, \lambda) - W^-(0; -c, d, \lambda) \rangle .
\end{align*}
Since $W(x; c, d; \lambda)$ is linear in $(c, d)$,
\begin{align*}
E(-\lambda) &= K E(\lambda), \quad K = \begin{pmatrix}-1 & 0 \\ 0 & 1 \end{pmatrix},
\end{align*}
thus $\det E(-\lambda) = \det K \det E(\lambda) = -\det E(\lambda)$. We compute the determinant directly with the aid of Wolfram Mathematica to get \cref{1pblockmatrixdet}.

\subsection{Proof of Theorem \ref{theorem:1pess} }

First, we make a change of variables to simplify the problem. Since $\nu(0) = 0$ and $\nu'(0) = 1/c$, $\nu(\lambda)$ is invertible near 0. Let $\lambda = \nu^{-1}(\mu)$. Expanding in a Taylor series about $\mu = 0$,
\begin{equation}\label{lambdamu}
\lambda = \nu^{-1}(\mu) = c \mu + \mathcal{O}(\mu^3).
\end{equation}
Substituting this into \cref{1pblockmatrixdet}, dividing by $c \mu^2$ since we are looking for the nonzero essential spectrum eigenvalues, and simplifying, we wish to solve
\begin{equation}\label{1pdetmu2}
2 M \sinh(\mu X) + c K \mu \cosh(\mu X) + \mathcal{O}\left(|\mu|(|\mu| + r^{1/2} X )\right) = 0,
\end{equation}
where $K = M \tilde{M} + M_c^2$. Choose any positive integer $m \leq N$, and take the ansatz
\begin{equation}\label{singlemu}
\mu = \frac{m \pi i}{X + c \frac{K}{2 M}} + \frac{h}{X}.
\end{equation}
For sufficiently large $X$, expand the denominator of \cref{singlemu} in a Taylor series to get
\begin{align}\label{singlemuTaylor}
\mu &= \frac{m \pi i}{X\left(1  + c \frac{K}{2 M X} \right) } + \frac{h}{X}
= \frac{m \pi i}{X}\left( 1 - c \frac{K}{2 M X} + c^2 \frac{K^2}{4 M^2 X^2} + \mathcal{O}\left(\frac{1}{X^3}\right) \right) + \frac{h}{X}.
\end{align}
Substituting this into \cref{1pdetmu2}, expanding the $\sinh$ and $\cosh$ terms in a Taylor series about $m \pi i$, and simplifying, equation \cref{1pdetmu2} is equivalent to 
\begin{align*}
\left( 2M + \frac{cK}{X} \right)h + \mathcal{O}\left( \frac{m+h}{X}\left(\frac{m+h}{X} + r^{1/2} X \right) \right) = 0.
\end{align*}
Since $X = X(r) = \mathcal{O}(|\log r|)$ and $|\log r|^{-1}$ is lower order than $r^{1/2} |\log r|$, we wish to solve
\begin{align*}
G_m(h, r) = \left( 2M + \mathcal{O}\left( \frac{1}{|\log r|} \right)\right) h + \mathcal{O}\left( \frac{(m+h)^2}{|\log r|^2} \right) = 0.
\end{align*}
Since $G_m(0,0) = 0$ and $\partial_h G_m(0,0) = 2M \neq 0$, by the implicit function theorem, there exists $r_1^m \leq r_*$ and a unique smooth function $h_m(r)$ with $h_m(0) = 0$ such that for all $r \leq r_1^m$, $G_m(h_m(r),r) = 0$. Expanding $h_m(r)$ in a Taylor series about $r = 0$,
$h_m(r) = \mathcal{O}\left( \frac{m^2}{|\log r|^2} \right)$.
Let $r_1 = \min\{ r_1^1, \dots, r_1^N \}$. Substituting $h_m(r)$ into \cref{singlemu}, the essential spectrum eigenvalues are located at
\begin{align*}
\mu_m(r) &= \frac{m \pi i}{X + c \frac{K}{2 M}} + \mathcal{O}\left( \frac{m^3}{|\log r|^3} \right) & m = 1, \dots, N.
\end{align*}
Changing variables back to $\lambda$ and using \cref{lambdamu},
\[
\lambda_m^{\text{ess}}(r) = c  \frac{m \pi i}{X + c \frac{M \tilde{M} + M_c^2}{2 M}} + \mathcal{O}\left( \frac{m^3}{|\log r|^3} \right),
\]
from which we obtain \cref{1pess} by factoring out $X$ from the denominator. By Hamiltonian symmetry (and the symmetry of $E(\lambda)$), eigenvalues must come in quartets. Since there is nothing else above the real axis with similar magnitude, $\lambda_m^{\text{ess}}(r)$ is on the imaginary axis, and there is another essential spectrum eigenvalue at $-\lambda_m^{\text{ess}}(r)$. 

\section{Proof of results for periodic double pulse}\label{sec:doublepulse}

\subsection{Proof of Lemma \ref{lemma:2blockmatrix} }

For the periodic double pulse, we have the symmetry relation $Q_i^-(x) = R Q_{i-1}^+(-x)$, where the subscript $i$ is taken $\Mod 2$. Let $(W_0^-(x), W_0^+(x), W_1^-(x), W_1^+(x))$ be the unique solution to \cref{eigsystem0} for given $(c_1, c_0, d_1, d_0, \lambda)$. Following the same procedure as in the proof of \cref{lemma:1blockmatrix}, when $(c_1, c_0, d_1, d_0, \lambda) \mapsto (-c_1, -c_0, d_0, d_1, -\lambda)$, equations \cref{eigsystem0} are satisfied by $(-R W_1^+(-x), -R W_1^-(-x), -R W_0^+(-x), -R W_1^-(-x))$, which must be the same as the original solution by uniqueness. We then compute the jump conditions to get $E(-\lambda) = K_1 E(\lambda) K_2$, where
\begin{align*}
K_1 = \begin{pmatrix}T & 0 \\ 0 & T \end{pmatrix}, \quad
K_2 = \begin{pmatrix}-I & 0 \\ 0 & T \end{pmatrix}, \quad
T = \begin{pmatrix} 0 & 1 \\ 1 & 0 \end{pmatrix},
\end{align*}
from which follows that $\det E(-\lambda) = \det K_1 \det E(\lambda) \det K_2 = -\det E(\lambda)$. The form of the remainder matrix follows from this symmetry together with the fact that $\lambda = 0$ is an eigenvalue with at least algebraic multiplicity 3. The determinant of $E(\lambda)$ is then computed directly with the aid of Wolfram Mathematica. 

\subsection{Change of variables}

As in the proof of \cref{theorem:1pess}, we make the change of variables $\lambda = \nu^{-1}(\mu)$ so that $\det E(\lambda)$ becomes 
\begin{equation}\label{2pdetmu}
\begin{aligned}
\det E(&\mu, r) = -2 \mu^2 c^2 (2a + \mu^2 c^2 + R_1) \left( M \sinh(\mu X) + c \mu (M \tilde{M} + M_c^2 ) \cosh(\mu X) \right) \\
&+4 a c^3 \mu^3 M_c^2 \sinh(\mu X_1)\sinh(\mu X_0) 
+ R_2 \mu^2 \sinh(\mu(X_1 - X_0)) \\
&+ R_3 \mu^2 \sinh(\mu X) + \mu^3 R_4.
\end{aligned}
\end{equation}
The remainder terms have the same bounds as in \cref{lemma:2blockmatrix} with $\lambda$ replaced by $\mu$. Define
\begin{align}
\mu^*_{int} &= \sqrt{-\frac{2 a}{M c^2}} \label{defmustarint} \\
\mu^*_m &= \frac{m \pi i}{X + c \frac{K}{M}}  && m = 1, \dots, N, \label{defmustaress}
\end{align}
where $K = M \tilde{M} + M_c^2$, and $a$ is defined in \cref{2pa}. For \cref{theorem:2peigsassym}, we will ensure that the interaction eigenvalues and essential spectrum eigenvalues do not interfere. Since $\mu^*_{int} = \mathcal{O}(r^{1/2})$ and $\mu^*_1 = \mathcal{O}(|\log r|^{-1})$, choose $r_0 \leq r_*$ sufficiently small so that $|\mu^*_{int}| \leq \frac{1}{2} |\mu^*_1|$ for all $r \leq r_0$. We can then simplify $\det E(\mu, r)$ to obtain 
\begin{equation}\label{2pdetmu2}
\begin{aligned}
\det E(\mu, r) &= -2 c^2 \mu^2 (2a + c^2 \mu^2 M) \left( M \sinh(\mu X) + c \mu (M \tilde{M} + M_c^2 ) \cosh(\mu X) \right) \\
&\qquad + \tilde{R} \mu^2 \sinh(\mu X) + R \mu^3,
\end{aligned}
\end{equation}
where $|\tilde{R}|, |R| \leq C(r^{1/2} + |\mu|)^3$. 

\subsection{Essential spectrum eigenvalues}

Choose any positive integer $m \leq N$, and take the ansatz $\mu = \mu^*_m + \frac{h}{X}$. As in the proof of \cref{theorem:1pess}, we expand the denominator of $\mu^*_m$ in a Taylor series, substitute the result into \cref{2pdetmu2}, and simplify to obtain the equation
\begin{align*}
G_m(h, r) = \left( M + \mathcal{O}\left( \frac{1}{|\log r|} \right)\right) h + \mathcal{O}\left( \frac{(m+h)^2}{|\log r|^2} \right) = 0.
\end{align*}
Since $G_m(0,0) = 0$ and $\partial_h G_m(0,0) = M \neq 0$, the result follows from the implicit function theorem and Hamiltonian symmetry as in the proof of \cref{theorem:1pess}.

\subsection{Interaction eigenvalues}\label{sec:2pinteigs}

The interaction eigenvalues will be at approximately $\mu = \pm \mu_{int}^*$, which is defined in \cref{defmustarint}. The interaction pattern will thus be determined by the sign of $a$. The first step is to characterize $a$. We start with the following lemma.

\begin{lemma}\label{lemma:Hoverlaplemma}
Define the functions $H: \R^+ \times \R^+ \rightarrow \R$ and $H_1: \R^+ \times \R^+ \rightarrow \R$ by 
\begin{align}
H(b_0, b_1) &= b_0 \sin \left( -\rho \log b_0 \right) - b_1 \sin \left( -\rho \log b_1 \right) \label{perdefH} \\
H_1(b_0, b_1) &= b_0 \left[ \rho \cos \left( -\rho \log b_0 \right) - \sin \left( -\rho \log b_0 \right) \right] + b_1 \left[ \rho \cos \left( -\rho \log b_1 \right) - \sin \left( -\rho \log b_1 \right) \right]. \label{perdefH1}
\end{align}
Then the zero sets of $H(b_0, b_1)$ and $H_1(b_0, b_1)$ intersect only at the discrete set of points $(b_0, b_1) = (b_k^*, b_k^*)$ for $k \in \Z$, where $b^*_k = e^{-\frac{1}{\rho} (k \pi + p^*) }$ are the pitchfork bifurcation points from \cref{2pulsebifurcation}.
\begin{proof}
Let $f(x) = x \sin \left( -\rho \log x \right)$. Then $H(b_0, b_1) = f(b_0) - f(b_1)$ and $H_1(b_0, b_1) = -b_0 f'(b_0) - b_1 f'(b_1)$. If $H(b_0, b_1) = 0$ then $f'(b_0) = f'(b_1)$, thus since $H_1(b_0, b_1) = 0$, $(b_0 + b_1) f'(b_0) = 0$. If $f'(b_0) = 0$, then by \cref{pitchforkH}, $b_0 = b_k^*$, and so $b_1 = b_k^*$ as well. Otherwise, $b_0 = -b_1$, which is not possible since $b_0$ and $b_1$ are both positive.
\end{proof}
\end{lemma}

We can now characterize $a$ in terms of the parameterization of the periodic double pulse. 

\begin{lemma}\label{lemma:chara}
Let $r_*$ be as in \cref{2pulsebifurcation}. Then for any $r \in \mathcal{R}$ with $r \leq r_*$,
\begin{enumerate}[(i)]
	\item For a symmetric periodic 2-pulse $\tilde{Q}_2(x; m_0, s_0, r)$, $a = r \tilde{a}(r; m_0, s_0)$, where $\tilde{a}(r; m_0, s_0)$ is continuous in $r$. Furthermore $\tilde{a}(0; m_0, s_0) = 0$ if and only if $s_0 = p^*$. For $s_0 \neq p^*$,
	\begin{equation}\label{symmsigns}
	\begin{aligned}
	\tilde{a}(0; 0, s_0) > 0 \text{ and } \tilde{a}(0; 1, s_0) < 0 && \text{if }s_0 > p^* \\
	\tilde{a}(0; 0, s_0) < 0 \text{ and } \tilde{a}(0; 1, s_0) > 0 && \text{if }s_0 < p^*.
	\end{aligned}
	\end{equation}

	\item For an asymmetric periodic 2-pulse $Q_2(x; m_0, s_1, r)$, $a = r \tilde{a}(r; m_0, s_1)$, where $\tilde{a}(r; m_0, s_0)$ is continuous in $r$. The sign of $\tilde{a}(0; m_0, s_0)$ is completely determined by $m_0$ and is given for all $s_1 > p^*$ by
	\begin{equation}\label{asymmsigns}
	\begin{aligned}
	\tilde{a}(0; m_0, s_1) &< 0 && \text{if }m_0 = 0 \\
	\tilde{a}(0; m_0, s_1) &> 0 && \text{if }m_0 = 1.
	\end{aligned}
	\end{equation}
\end{enumerate}	
\begin{proof}
Using \cite[Lemma 6.1(ii)]{Sandstede1998}, rescaling as in \cref{jumplemma3}, and simplifying, 
\begin{align*}
\langle \Psi(X_i), Q'(-X_i) \rangle &=
-p_0 \alpha_0 e^{\alpha_0 \phi/\beta_0} r b_i \left( \rho \cos\left(-\rho \log b_i \right) - \sin \left(-\rho \log b_i \right) \right) + \mathcal{O}(r^{1+\gamma/2\alpha_0}),
\end{align*}
from which it follows that $a = -p_0 \alpha_0 e^{\alpha_0 \phi/\beta_0} r H_1(b_0, b_1)$, where $H_1$ is defined in \cref{perdefH1}. For part (i), $a = r \tilde{a}(r; m_0, s_0)$, where
\begin{align*}
\tilde{a}(r; m_0, s_0) &= -p_0 \alpha_0 e^{\alpha_0 \phi/\beta_0} r H_1(b_0, b_1) \\
&= -2 p_0 \alpha_0 e^{\alpha \phi/\beta_0} e^{-\frac{1}{\rho}(m_0 \pi + s_0)} (-1)^{m_0} \left( \rho \cos s_0 - \sin s_0 \right) + \mathcal{O}(r^{\gamma/2\alpha}).
\end{align*}
When $r = 0$, $\tilde{a}(0; m_0, s_0) = 0$ if and only if $s_0 = p^* = \arctan \rho$, and the sign of $\tilde{a}(0; m_0, s_0)$ is given by \cref{symmsigns} for $s_0 \neq p^*$. 

For part (ii), $a = r \tilde{a}(r; m_0, s_1)$, where
\begin{align*}
\tilde{a}(r; m_0, s_1) &= -p_0 \alpha_0 e^{\alpha \phi/\beta_0} H_1( b_0(m_0, s_1), b_1(s_1) ) + \mathcal{O}(r^{\gamma/2\alpha}),
\end{align*}
and the parameterization $(b_0(m_0, s_1), b_1(s_1))$ is defined in \cref{armpersists}. Since $(b_0(m_0, s_1), b_1(s_1))$ is contained entirely within the zero set of $H(b_0, b_1)$, where $H$ is defined in \cref{perdefH}, it follows from \cref{lemma:Hoverlaplemma} that $\tilde{a}(0; m_0, s_1) = 0$ if and only if $s_1 = p^*$. By continuity in $s_1$, $\tilde{a}(0; m_0, s_1)$ must have the same sign for all $s_1 > p^*$. Sending $s_1 \rightarrow \infty$,
\begin{align*}
\lim_{s_1 \rightarrow \infty}
\tilde{a}(0; m_0, s_1) = -p_0 \beta_0 e^{\alpha \phi/\beta_0} e^{-\frac{1}{\rho} m_0 \pi} \cos(m_0 \pi) = (-1)^{m_0+1} p_0 \beta_0 e^{\alpha \phi/\beta_0} e^{-\frac{1}{\rho} m_0 \pi},
\end{align*}
thus since $p_0 > 0$ and $\beta_0 > 0$, the sign of $\tilde{a}(0; m_0, s_1)$ is given by \cref{asymmsigns} for $s_1 > p^*$.
\end{proof}
\end{lemma}

We can now find the interaction eigenvalues. To rescale equation \cref{2pdetmu2}, let $a = r \tilde{a}(r)$, where $\tilde{a}(r) = \tilde{a}(r; m_0, s_1)$ is defined in \cref{lemma:chara} and $\tilde{a}(0) \neq 0$ since $s_1 > p^*$. Let $\mu^*_{int}(r) = r^{1/2} \tilde{\mu}^*(r)$, where $\tilde{\mu}^*(0) \neq 0$, and take the ansatz $\mu = r^{1/2} \tilde{\mu}^*(r) + r^{1/2}$. Substituting this into \cref{2pdetmu2}, dividing by $-2 c^2 r^{5/2} X$, using the fact that $\mu = \mathcal{O}(r^{1/2})$ and $X = \mathcal{O}(|\log r|)$, and simplifying, we obtain the equation
\[
G(h, r) = M c^2 h (\tilde{\mu}^*(r) + h)^3 (2 \tilde{\mu}^*(r) + h)
\left( M + \mathcal{O}\left(\frac{1}{|\log r|} \right) \right)
 + \mathcal{O}\left( r^{1/2}(\tilde{\mu}^*(r) + h) \right) = 0.
\]
Since $G(0, 0) = 0$ and $\partial_h G(0, 0) = 2 M^2 c^2 (\tilde{\mu}^*(0))^4 \neq 0$, by the implicit function theorem there exists $r_1 \leq r_0$ and a unique smooth function $h(r)$ with $h(0) = 0$ such that for all $r \leq r_1$, $G(h(r), r) = 0$. Expanding $h(r)$ in a Taylor series about $r = 0$, $h(r) = \mathcal{O}(r^{1/2})$. Undoing the rescaling and changing variables back to $\lambda$, there is an interaction eigenvalue located at 
\begin{align*}
\lambda^{\text{int}}(r) = \sqrt{-\frac{2 \tilde{a}(r)}{M}}r^{1/2} + \mathcal{O}(r) = \sqrt{ -\frac{2 a}{M} } + \mathcal{O}(r).
\end{align*}
By Hamiltonian symmetry, there is also an eigenvalue at $-\lambda(r)$. Since eigenvalues must come in quartets, and there only two eigenvalues of this magnitude, we conclude that for $r \leq r_1$, there is a pair of interaction eigenvalues given by $\lambda = \pm \lambda^{\text{int}}(r)$, which is either real or purely imaginary. Since $M = d''(c) > 0$ by \cref{hyp:dccpos}, these are real if $\tilde{a}(0) < 0$ and purely imaginary if $\tilde{a}(0) > 0$, which depends only on $m_0$ by \cref{lemma:chara}.

\subsection{Eigenvalues at 0}\label{sec:asymmeigcount0}

In this section, we use Rouch\'{e}'s theorem to show that there are exactly three eigenvalues at 0. As in the previous section, let $\mu^*_{int} = r^{1/2} \tilde{\mu}^*(r)$, and let $\xi = \frac{1}{2}|\tilde{\mu}^*(0)|$. Let $\mu = r^{1/2} \tilde{\mu}$, and take $\tilde{\mu}$ on the circle $|\tilde{\mu}| = \xi$. Making these substitutions, dividing by $-2 c^2 r^{5/2} X$, and simplifying, equation \cref{2pdetmu2} is equivalent to
\begin{equation}\label{def:eig0G}
G(\tilde{\mu}, r) = M c^2 \tilde{\mu}^3 (\tilde{\mu} + \tilde{\mu}^*(r))(\tilde{\mu} - \tilde{\mu}^*(r))
\left( M + \mathcal{O}\left(\frac{1}{|\log r|} \right) \right)
 + \mathcal{O}\left( r^{1/2} \right) = 0.
\end{equation}
Let $G(\tilde{\mu}, r) = G(\tilde{\mu}, 0) + G_1(\tilde{\mu}, r)$, where $G_1(\tilde{\mu}, r) = G(\tilde{\mu}, r) - G(\tilde{\mu}, 0)$. On the circle $|\tilde{\mu}| = \xi$, $|G(\tilde{\mu}, 0)| \geq \frac{Mc^2}{16}|\tilde{\mu}_*(0)|^5$. Since $|G_1(\tilde{\mu}, r)| \rightarrow 0$ as $r \rightarrow 0$ and $\tilde{\mu}_*(0) \neq 0$, there exists $r_2 \leq r_0$ such that for $r \leq r_2$, $|G_1(\tilde{\mu}, r)| < |G(\tilde{\mu}, 0)|$ on the circle $|\tilde{\mu}| = \tilde{\xi}$. By Rouch\'{e}'s theorem $G(\tilde{\mu}, r)$ and $G(\tilde{\mu}, 0)$ have the same number of zeros (counted with multiplicity) inside the circle of radius $\xi$. By the choice of $\tilde{\xi}$, $G_1(\tilde{\mu})$ has exactly 3 zeros inside the circle, thus $G(\tilde{\mu}, r)$ does as well. Undoing the scaling and changing variables back to $\lambda$, we obtain the result.

\subsection{Proof of Theorem \ref{theorem:2peigssym}}

Let $\tilde{a}(r, s_0) = \tilde{a}(r; m_0, s_0)$ as in \cref{lemma:chara}, so that $\tilde{a}(0, p^*) = 0$ and $\tilde{a}(r, s_0) \rightarrow 0$ as $r \rightarrow 0$ and $s_0 \rightarrow p^*$. Using Rouch\'{e}'s as in the previous section, for sufficiently small $r$, there are five zeros of $E(\lambda, r)$ in a small ball around the origin. Since two of these eigenvalues correspond to the kernel eigenfunctions $\partial_x \tilde{Q}_2(x)$ and $\partial_c \tilde{Q}_2(x)$, and there is third kernel eigenfunction by \cref{qnkernel}, this leaves two eigenvalues unaccounted for. Let $s_0 = p^* + h$. Using the scaling $\mu = r^{1/2} \tilde{\mu}$, following the same steps as in the previous section, factoring out $\tilde{\mu}^3$ from \cref{def:eig0G} since we have already accounted for three eigenvalues at 0, and simplifying, the remaining eigenvalues must satisfy the equation
\begin{equation}\label{def:symmG}
G(\tilde{\mu}, r, h) = \left( \tilde{\mu}^2 + \frac{2 \tilde{a}(r, p^* + h)}{M c^2} \right)
\left( M + \mathcal{O}\left(\frac{1}{|\log r|} \right) \right)
 + \mathcal{O}\left( r^{1/2} \right) = 0.
\end{equation}
When $r = 0$ and $h = 0$, $G(\tilde{\mu}, 0, 0)$ has a double root at $\tilde{\mu} = 0$. We will show that for small $r$, the double root persists at $h = h(r)$. For $\tilde{\mu}$ to be a double root, it must satisfy $K(\tilde{\mu}, r, h) = 0$, where
\begin{equation*}
K(\tilde{\mu}, r, h) = 
\begin{pmatrix}G(\tilde{\mu}, r, h) \\ \partial_{\tilde{\mu}}G(\tilde{\mu}, r, h) \end{pmatrix}. 
\end{equation*}
For $r = 0$, $K(0, 0, 0) = 0$. Using \cref{lemma:chara}, $\tilde{a}(0, p^* + h) = c_0 h + \mathcal{O}(|h|^2)$ for some constant $c_0 \neq 0$, thus we have
\begin{equation*}
D_{(\tilde{\mu}, h)}K(0, 0, 0) = 
\begin{pmatrix}
0 & \frac{2 c_0}{c^2} \\
2M & 0
\end{pmatrix},
\end{equation*}
which is nonsingular. Using the implicit function theorem, there exists $r_1 \leq r_0$ and a unique smooth function $(\tilde{\mu}(r), h(r))$ with $(\tilde{\mu}(0), h(0)) = (0, 0)$ such that for all $r \leq r_1$, $K(\tilde{\mu}(r), h(r), r) = 0$. In other words, $\tilde{\mu}(r)$ is a double root of $G(\tilde{\mu}, r, h)$ when $h = h(r)$. By Hamiltonian symmetry, we must have $\tilde{\mu}(r) = 0$, thus for $r \leq r_1$, there are two more eigenvalues at 0 when $s_0 = p^* + h(r)$, which brings the total to five. Since the pitchfork bifurcation in the family of periodic 2-pulses occurs near $p^*$ at $p^*(r)$, it follows from standard PDE bifurcation theory that this quintuple zero must occur at the pitchfork bifurcation point.

\subsection{Proof of Theorem \ref{th:Kreinbubble}}

Choose $m_0 = 1$, let $r_*$ be as in \cref{2pulsebifurcation}, and let $Q_2(x; s_1, r)$ be the family of periodic 2-pulse solutions parameterized by $s_1$. By \cref{theorem:2peigsassym}, there is a pair of purely imaginary interaction eigenvalues for sufficiently small $r$. Define
\begin{align*}
\mu_1(s_1, r) &= \frac{\pi i}{X(s_1,r) + c \frac{K}{M}} \\
\mu_*(s_1, r) &= \sqrt{-\frac{ 2a(s_1,r) - R_1(r) }{M c^2}} = \sqrt{-\frac{2a(s_1,r)}{M c^2}} + \mathcal{O}(r),
\end{align*}
both of which lie on the imaginary axis. Note that we have included the remainder term $R_1(r)$, which independent of $\mu$, in the definition of $\mu_*$. As $s_1 \rightarrow \infty$, $X_1(r, s_1) \rightarrow \infty$ and $X_0(r, s_1) \rightarrow \frac{1}{2 \alpha_0}|\log r| + \frac{\pi}{2\beta_0} + \tilde{L}$, which is a nonzero constant. By the definition of $a$ in \cref{2pa}, $a(s_1,r)$ approaches a nonzero constant as $s_1 \rightarrow \infty$. Thus there exists $r_0 \leq r_*$ such that for all $r \leq r_0$, we can find $s_*(r)$ such that $\mu_1(s_*(r), r) = \mu_*(s_*(r), r)$. Let
\begin{equation}\label{XstarKrein}
X_*(r) = X(s_*(r), r) = \frac{\pi i}{\mu_*(s_*(r), r)} - c \frac{K}{M}.
\end{equation}
For $s_1$ close to $s_*(r)$, define the parameter $k(s_1, r)$ by 
\[
2 k(s_1, r) i = \mu_*(s_1, r) - \mu_1(s_1, r).
\]
This measures the distance on the imaginary axis between $\mu_1$ and $\mu_*$, so that $k(s_*(r), r) = 0$. From this point forward, we will drop the dependence on $s_1$ and $r$ for convenience of notation.

Let $\mu = \mu_* + h$, so that $\mu = \mu_1 + h + 2 k i$. As in \cref{singlemuTaylor}, we expand $\mu_1$ in a Taylor series to get
\begin{align}\label{muTaylor2}
\mu &= \frac{m \pi i}{X}\left( 1 - c \frac{K}{2 M X} + c^2 \frac{K^2}{4 M^2 X^2} + \mathcal{O}\left(\frac{1}{X^3}\right) \right) + h + 2 k i.
\end{align}
Anticipating what will follow, we take the rescaling
\begin{align*}
\mu_* &= r^{1/2} \tilde{\mu}, \quad h = r^{5/4} X_0 \tilde{h}, \quad k = r^{5/4} X_0 \tilde{k}.
\end{align*}
Using this rescaling and the fact that $1/X(r) = \mathcal{O}(r^{1/2})$ in this case, we substitute \cref{muTaylor2} into $\sinh \mu X$ and $\cosh \mu X$ to get
\begin{equation}\label{sinhcoshmuX}
\begin{aligned}
\sinh(\mu X) &= (-1)\left( - c \frac{K \pi i }{MX} + c^2 \frac{K^2 \pi i }{M^2 X^2} + (h + 2 k i)X \right) + \mathcal{O}\left( r^{3/2} \right) \\
\quad \cosh(\mu X) &= (-1) + \mathcal{O}\left( r^{3/2} \right).
\end{aligned}
\end{equation}
It follows from \cref{sinhcoshmuX} that $\sinh(\mu X) =  \mathcal{O}(r^{1/2})$. Furthermore, 
\begin{align*}
\sinh(\mu X_1)\sinh(\mu X_0) &= \left( \sinh(\mu X) \cosh( \mu X_0)  - \cosh(\mu X)\sinh(\mu X_0) \right)\sinh(\mu X_0) \\
&= \mu_*^2 X_0^2 + \mathcal{O}(r^{1/2})
\end{align*}
and
\begin{align*}
&\sinh(\mu(X_1 - X_0)) = \sinh(\mu X) \cosh(2 \mu X_0) - \cosh(\mu X)\sinh(2 \mu X_0)
= \mathcal{O}\left( r^{1/2}|\log r| \right).
\end{align*}
Substituting all of these into \cref{2pdetmu}, dividing by $-2 \mu^2 c^2$ since $\mu \neq 0$, using the expression $a = -M c^2 \mu_*^2/2 + \mathcal{O}(r^{3/2})$, and simplifying, we obtain the equation
\begin{equation}\label{det3}
\begin{aligned}
&M \tilde{\mu} \tilde{h} X_0 r^{7/4}
\left( \left( M + \frac{c K}{X} \right) X ( \tilde{h} + 2 \tilde{k} i) X_0 r^{5/4} + \mathcal{O}(r^{3/2}) \right) \\
&\qquad\qquad\qquad\qquad\qquad- 2 M M_c^2 c \tilde{\mu}^5 r^{5/2} X_0^2 + \mathcal{O}( r^{5/2} X_0 ) = 0.
\end{aligned}
\end{equation}
Since
\[
\left( M + \frac{c K}{X} \right) X  = MX + c K = \frac{M \pi i}{c \mu_1} = \frac{M \pi i}{c(\mu_* - 2 k i)}
= \frac{M \pi i}{c \tilde{\mu}r^{1/2}} + \mathcal{O}(r^{1/4}),
\]
equation \cref{det3} simplifies to
\begin{equation}\label{det4}
\begin{aligned}
2 \pi c M^2 X_0^2 r^{5/2} \tilde{h}( \tilde{h} + 2 \tilde{k} i) - 2 M M_c^2 c^3 | \tilde{\mu}|^5 r^{5/2} X_0^2 + \mathcal{O}( r^{5/2} |\log r| ) = 0.
\end{aligned}
\end{equation}
Dividing by $X_0^2 r^{5/2}$ and noting that $X_0 = \mathcal{O}(|\log r|)$, we obtain the equation
\begin{equation}\label{det5}
\begin{aligned}
G(\tilde{h}, \tilde{k}, r) = \tilde{h}^2 + 2 \tilde{k}  i \tilde{h}
- T + \mathcal{O}\left( \frac{1}{|\log r|} \right) = 0, \quad T = \frac{M_c^2 }{2 \pi M^2 } c^2 |\tilde{\mu}|^5 > 0.
\end{aligned}
\end{equation}
When $r = 0$, $G(\tilde{h}, \tilde{k}, 0) = 0$ if and only if $\tilde{h} = \tilde{h}_0^\pm(\tilde{k}) = -\tilde{k}i \pm \sqrt{T - \tilde{k}^2}$. When $|k| \leq \sqrt{T}$, this is a circle in the complex plane of radius $\sqrt{T}$. We will show this persists for small $r$. For $\tilde{h} = \tilde{h}_0^\pm(\tilde{k})$, $G(\tilde{h}_0^\pm(\tilde{k}), \tilde{k}, 0) = 0$ and $\partial_{\tilde{h}}G(\tilde{h}_0^\pm(\tilde{k}), \tilde{k}, 0) = 2 \sqrt{T - \tilde{k}^2}$, which is nonzero as long as $\tilde{k} \neq \pm \sqrt{T}$. Choose any $\epsilon > 0$ and define
\[
K_\epsilon = \left[-2 \sqrt{T} -\sqrt{T} - \epsilon\right]
\cup \left[-\sqrt{T} + \epsilon, \sqrt{T} - \epsilon\right]
\cup \left[\sqrt{T} + \epsilon, 2\sqrt{T}\right].
\]
Then for $\tilde{k} \in K_\epsilon$, $G_{\tilde{h}_0}(\tilde{h}_0(\tilde{k}), \tilde{k}, 0)$ is bounded away from 0, with bound dependent on $\epsilon$. By the uniform contraction mapping principle, there exists $r_1 \leq r_0$ and smooth functions $\tilde{h}_*^\pm(\tilde{k}, r)$ such that $h_*^\pm(\tilde{k}, 0) = \tilde{h}_0^\pm(\tilde{k})$, $h_*^\pm(\tilde{k}, r) = \tilde{h}_0^\pm(\tilde{k}) + \mathcal{O}\left(\frac{1}{|\log r|}\right)$, and for all $r \leq r_1$ and $\tilde{k} \in K_\epsilon$, $G(\tilde{h}_*^\pm(\tilde{k},r),\tilde{k},r) = 0$. For $k \in [-\sqrt{T} + \epsilon, \sqrt{T} - \epsilon]$, $\tilde{h}_*^\pm(\tilde{k}, r)$ is real and is symmetric across the imaginary axis by Hamiltonian symmetry. All that remains is show there is a double zero on the imaginary axis when $\tilde{k}$ is close to $\pm \sqrt{T}$. Similar to the proof of \cref{theorem:2peigssym}, for $\tilde{h}$ to be a double root of \cref{det5}, it must satisfy $H(\tilde{h}, \tilde{k}, r)$, where
\[
H(\tilde{h}, \tilde{k}, r) = \begin{pmatrix} G(\tilde{h}, \tilde{k}, r) \\
\partial_{\tilde{h}} G(\tilde{h}, \tilde{k}, r) \end{pmatrix}.
\]
For $r = 0$, $H(-\sqrt{T}i, \sqrt{T}, 0) = 0$, and
\begin{equation*}
D_{(\tilde{h}, \tilde{k})}H(-\sqrt{T}i, \sqrt{T}, 0) = 
\begin{pmatrix}
0 & 2 \sqrt{T} \\
2 & 2i
\end{pmatrix},
\end{equation*}
which is nonsingular. Using the implicit function theorem, there exists $r_2 \leq r_1$ and a unique smooth function $(\tilde{h}_1(r), \tilde{k}_1(r))$ with $(\tilde{h}_1(0), \tilde{k}_1(0)) = (-\sqrt{T}i, \sqrt{T})$ and $(\tilde{h}_1(r), \tilde{k}_1(r)) = (-\sqrt{T}i, \sqrt{T}) + \mathcal{O}\left(\frac{1}{|\log r|}\right)$ such that for all $r \leq r_2$, $H(\tilde{h}_1(r), \tilde{k}_1(r), r) = 0$. In other words, $\tilde{h}_1(r)$ is a double root of $G(\tilde{h}, \tilde{k}, r)$ when $\tilde{k} = \tilde{k}(r)$. By Hamiltonian symmetry, $\tilde{h}(r)$ must lie on the imaginary axis. We have a similar result for $(\tilde{h}, \tilde{k}) = (\sqrt{T}i, -\sqrt{T})$.

Undoing the scaling and change of variables and letting $s = k/c$, for $r \leq r_2$ there is a pair of eigenvalues
\begin{align*}
\lambda(r) = \lambda_*(r) - s i \pm \sqrt{ T_1 -  s^2} + \mathcal{O}(r^{5/4}|\log r|^{1/2}), \quad T_1 = -\frac{M_c^2}{2 \pi M^2 c } |\lambda_*|^5 X_0^2,
\end{align*}
which is approximately a circle centered at $\lambda_*$ in the complex plane of radius $\sqrt{T_1}$. For $s = s_\pm(r) = \pm \sqrt{T_1}\left( 1 + \mathcal{O}\left(\frac{1}{|\log r|} \right) \right)$, there are double eigenvalues at $\lambda = \lambda_* + s_\pm(r) + \mathcal{O}\left(\frac{1}{|\log r|} \right)$, which is on the imaginary axis by Hamiltonian symmetry. Again by Hamiltonian symmetry, the eigenvalues $\lambda$ are symmetric about the imaginary axis for $s_-(r) < s < s_+(r)$ and are purely imaginary for $s < s_-(r)$ and $s > s_+(r)$. The maximum real part of $\lambda$ occurs when $s = 0$, and is approximately $\sqrt{T_1}$. Expanding $\mu_1$ in a Taylor series about $X = X_*$, this occurs when $X(r) = X_*(r) + \frac{2 c \pi s}{\lambda_*^2} + \mathcal{O}(s^2)$. In particular, the Krein collision (and reverse collision) occur when $X(r)$ is given by \cref{KreinDeltaX}.

\section{Conclusions}\label{sec:conclusions}

In this paper, we use Lin's method to construct periodic multi-pulse solutions to KdV5 and to determine the spectrum near the origin associated with these solutions. The results hold more generally for Hamiltonian systems which are both reversible and translation invariant,
including higher order KdV equations.
This technique allows us to compute both the eigenvalues resulting from interactions between neighboring pulses in the periodic multi-pulse structure and the essential spectrum eigenvalues resulting from the background state. Of note, we show that as the domain size is increased, brief instability bubbles form when interaction eigenvalues and essential spectrum eigenvalues of opposite Krein signature collide on the imaginary axis. These Krein bubbles can be found numerically, and their location and size agree with the theoretical results. Numerical timestepping experiments show that these Krein bubbles correspond to slowly growing, oscillatory instabilities.

In \cref{th:perexist}, we prove the existence of periodic multi-pulse solutions, but have to exclude certain periodic parameterizations to avoid bifurcation points. We only demonstrate the complete bifurcation structure for periodic 2-pulses (\cref{2pulsebifurcation}). It is likely that bifurcation structures exist for arbitrary periodic $n$-pulses, but that these become more complex as $n$ increases. As a first step to elucidating these bifurcations, it would be useful to study periodic 3-pulses by using AUTO for parameter continuation not only in the domain size $X$ but also in the parameters $c$ and $p$. This could suggest a theoretical result as well as a generalization to higher $n$. We could also extend the spectral results to arbitrary multi-pulses. As long as any Krein collisions are avoided, it should be possible to determine the interaction eigenvalue patten solely from the geometry of the periodic multi-pulse, as is the case with the periodic double pulse. It would also be worth looking at whether Krein bubbles occur for higher order periodic multi-pulses. As an initial step, we could construct periodic 3-pulses with two pairs of imaginary interaction eigenvalues as well as ``mixed'' periodic 3-pulses with one pair each of real and imaginary interaction eigenvalues. We could then determine numerically if Krein bubbles occur for either of these two periodic 3-pulses. We also note that we obtained better error estimates for the Krein bubble size in the numerical computations in \cref{sec:numerics} than were predicted by \cref{th:Kreinbubble}, thus it might be possible to obtain sharper error estimates analytically. 

\begin{figure}
\begin{center}
\includegraphics[width=8cm]{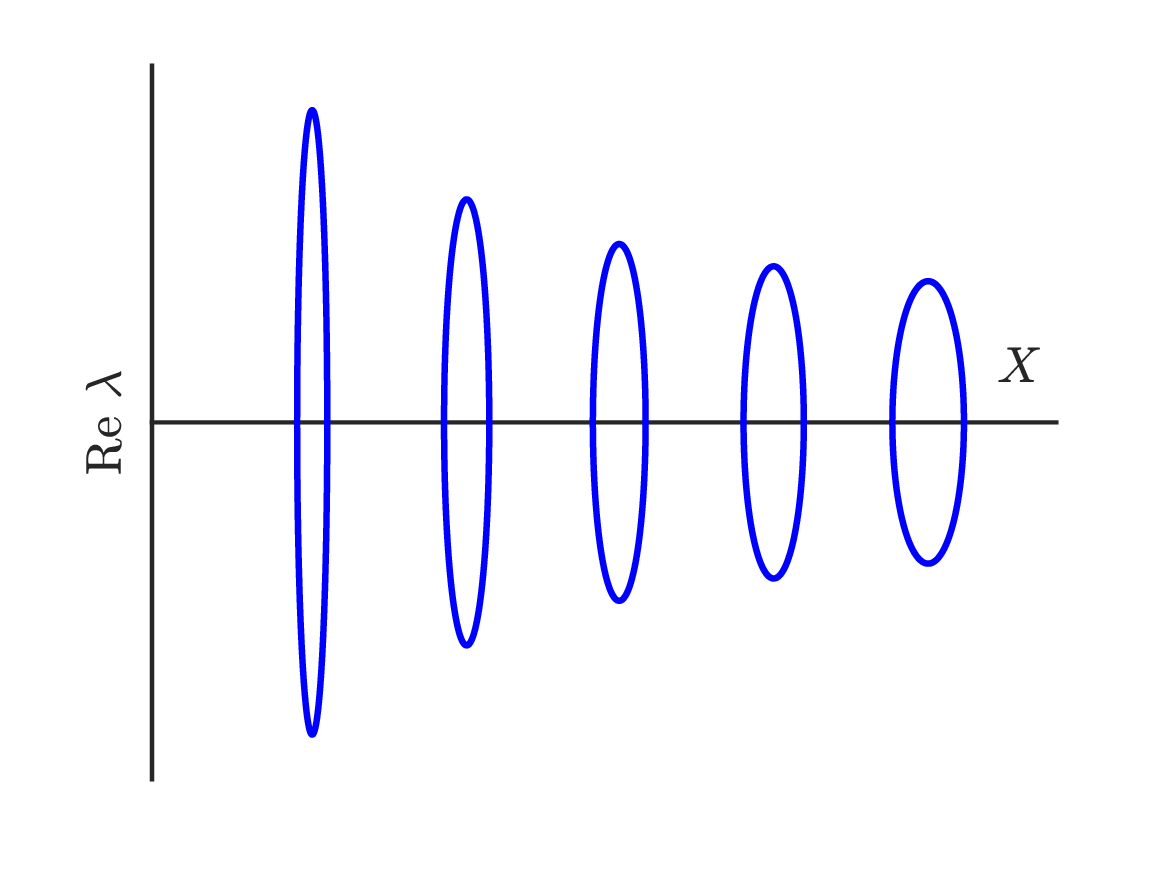}
\includegraphics[width=8cm]{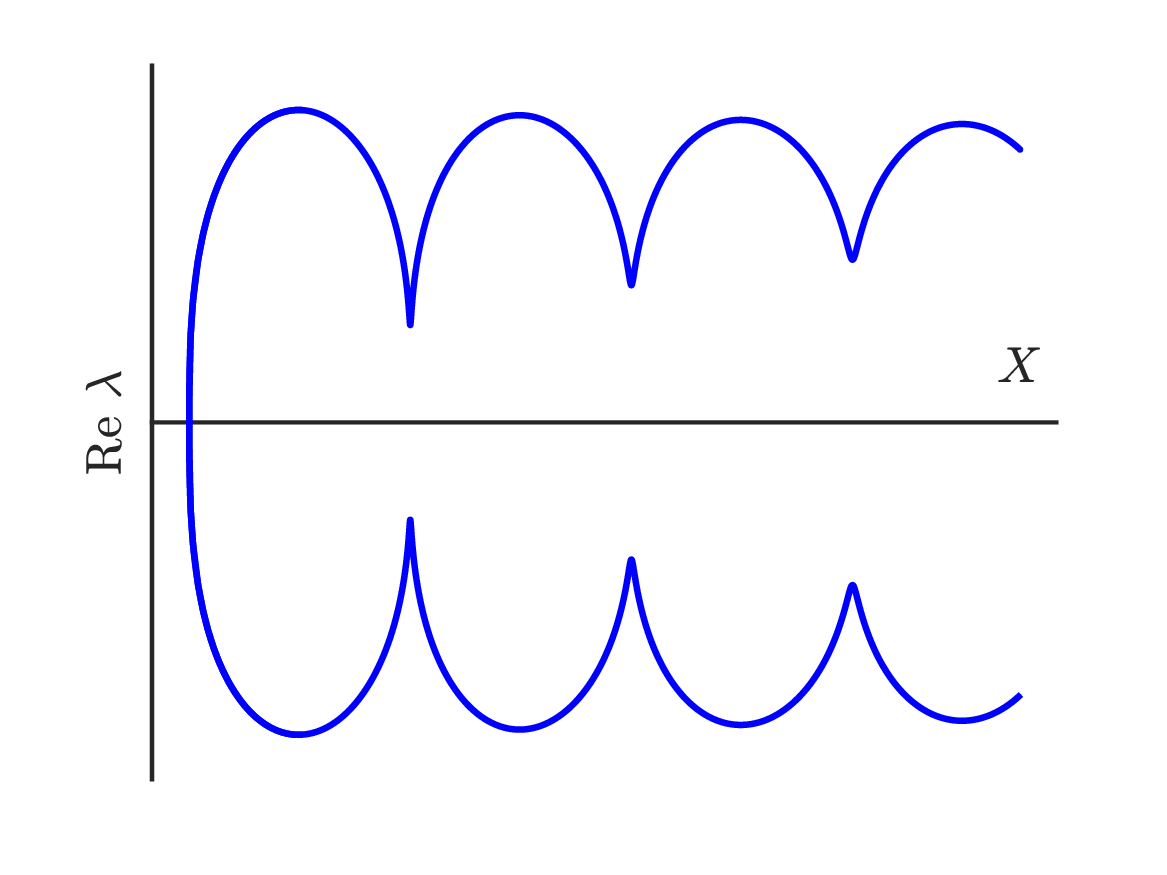}
\end{center}
\caption{Real part of $\lambda$ vs. $X$ for first five Krein bubbles as $X$ is increased (left). Real part of $\lambda$ vs. $X$ showing overlapping Krein bubbles as $X$ is further increased.
}
\label{fig:KreinBubbleCollision}
\end{figure}

Finally, we note from \cref{remark:kreinbubbles} that as the domain size $X$ increases, the centers of subsequent Krein bubbles (indexed by increasing $m$) are approximately equally spaced in $X$. However, their widths in $X$ scale with $\sqrt{m}$, so that they grow with each subsequent bubble (\cref{fig:KreinBubbleCollision}, left panel). At some critical value $X = X_c$, the Krein bubbles will start overlapping (\cref{fig:KreinBubbleCollision}, right panel). After this occurs, we expect that there will always be an eigenvalue with positive real part, thus the periodic double pulses should all be unstable for $X > X_c$. As $X$ is further increased, we expect that more than two Krein bubbles will interact, and that the eigenvalue behavior will become increasing complicated. (This is very difficult to simulate numerically, as it involves extremely large domain sizes.) At the same time, the radius in the complex plane of subsequent Krein bubbles scales with $1/\sqrt{m}$ (\cref{fig:KreinBubbleCollision}, left panel), which suggests that the maximum real part of the corresponding interaction eigenvalue approaches 0 as $X\rightarrow \infty$. Although we cannot generalize the periodic case to obtain the behavior on $\R$ by taking $X\rightarrow \infty$, this suggests that the spectrum of the double pulse on the real line, which is the formal limit of the periodic double pulse as $X\rightarrow \infty$, may in fact be purely imaginary.

\paragraph{Acknowledgments}

This material is based upon work supported by the U.S. National Science Foundation under grants DMS-1148284 (R.P. and B.S.), DMS-1840260 (R.P.), and DMS-1714429 (B.S.).

\bibliographystyle{elsarticle-num.bst}
\bibliography{kdv5periodic}

\end{document}